\documentclass[11pt]{amsart}
\usepackage{amssymb}
\usepackage{amscd}
\usepackage{lscape}
\usepackage[arrow,matrix,curve]{xy}
\usepackage{times}
\input epsf.tex

\newtheorem{Thm}{Theorem}[section]
\newtheorem{Lem}[Thm]{Lemma}
\newtheorem{Cor}[Thm]{Corollary}
\newtheorem{Prop}[Thm]{Proposition}

\newtheorem{example}{Example}

\newcommand{\A}{\mathbb{A}}
\newcommand{\D}{\mathbb{D}}
\newcommand{\E}{\mathbb{E}}
\newcommand{\Z}{\mathbb{Z}}
\newcommand{\Q}{\mathbb{Q}}

\newcommand{\N}{\mathbb{N}}
\newcommand{\C}{\mathbb{C}}

\newcommand{\M}{{\mathcal M}}

\newcommand{\rk}{\operatorname{rk}}
\newcommand{\iso}{\operatorname{iso}}
\newcommand{\ql}{\operatorname{ql}}  
\newcommand{\md}{\operatorname{mod}} 
\newcommand{\smd}{\underline{\md}}
\newcommand{\ind}{\operatorname{ind}}
\newcommand{\irr}{\operatorname{Irr}}

\newcommand{\Hom}{\operatorname{Hom}}
\newcommand{\sHom}{\underline{\Hom}} 
\newcommand{\Ext}{\operatorname{Ext}}
\newcommand{\End}{\operatorname{End}}
\newcommand{\ext}{\operatorname{ext}} 

\newcommand{\rep}{\operatorname{rep}}
\newcommand{\dm}{\operatorname{dim}}  
\newcommand{\GL}{\operatorname{GL}} 

\newcommand{\dimv}[1]{\mathbf{dim}(#1)} % or \operatorname{\mathbf{dim}}#1}
\newcommand{\dva}[1]{\mathbf{#1}}  % dimensionvectos for \md(\tLam)
\newcommand{\dvp}[1]{\dva{p}^\phg_{#1}}
\newcommand{\cla}[1]{\mathbf{#1}}  % elements of $K_0(\Da)=K_0(\smd(\tLam))

\newcommand{\supp}{\operatorname{supp}}
\newcommand{\field}{\C}
\newcommand{\bil}[1]{\langle #1\rangle}

\newcommand{\Lam}{\Lambda}
\newcommand{\lam}{\lambda}
\newcommand{\tLam}{{\widetilde{\Lam}}}

\newcommand{\Da}{\Delta}
\newcommand{\Gam}{\Gamma}
\newcommand{\sGam}{\underline{\Gam}}

\newcommand{\tQ}{\widetilde{Q}}
\newcommand{\tQa}[1]{\tQ^{[#1]}}
\newcommand{\Daa}[1]{\Da^{[#1]}}
\newcommand{\qda}[1]{q^{[#1]}}
\newcommand{\bia}[2]{\bil{#2}^{[#1]}}
\newcommand{\bha}[1]{\dva{h}^{[#1]}}

\newcommand{\cMa}[1]{\mathcal{M}^{[#1]}}
\newcommand{\cTa}[1]{\Tub^{[#1]}}
\newcommand{\cHe}{\mathcal{H}}

\newcommand{\Idv}{I}  % Set of dim vectors of indec \tLam modules

\newcommand{\ENn}{\N}
\newcommand{\ENp}{\Z^+}

\newcommand{\udel}{\underline{\delta}}

\newcommand{\phg}{{\phantom{g}}}
\newcommand{\phX}{\phantom{\bsm 1\\1\\1\\1\\1\esm\hspace{0.20cm}}}

\newcommand{\phZ}{\phantom{\bsm 1\\1\\1\esm}\hspace{0.16cm}}
\newcommand{\bsm}{\begin{smallmatrix}}
\newcommand{\esm}{\end{smallmatrix}}

\newcommand{\orb}{\mathcal O}
\newcommand{\Tub}{\mathcal{T}}

\newcommand{\msm}{\mathfrak{m}}  % multisegment m
\newcommand{\msp}{\mathfrak{p}}  % multisegment p
 
\newcommand{\X}{\mathbb{X}}
\newcommand{\coh}{\operatorname{coh}}
\newcommand{\rad}{\operatorname{rad}}

\def\<{\langle\,}
\def\>{\,\rangle}

\def\eg{{\em e.g. }}
\def\SG{\mathfrak S}

\def\F{{\mathcal F}}

\def\B{{\mathcal B}}
\def\U{{\mathbf U}}
\def\V{{\mathbf V}}

\def\g{\mathfrak g}

\def\<{\langle}
\def\>{\rangle}
\def\n{{\mathfrak n}}
\def\si{\sigma}
\def\lra{\longrightarrow}
\def\ra{\rightarrow}

\def\1{\mathbf 1}

\def\res{{\rm Res}}
\def\cc{{\mathbf c}}
\def\ii{{\mathbf i}}
\def\jj{{\mathbf j}}

\def\f{{\mathfrak f}}
\def\soc{{\rm Socle}}

\def\iu{\underline{i}}
\def\ju{\underline{j}}
\def\ku{\underline{k}}

\def\MM{{\mathbf M}}
\newdimen\Squaresize \Squaresize=14pt
\newdimen\Thickness \Thickness=0.5pt
 
\def\Square#1{\hbox{\vrule width \Thickness
   \vbox to \Squaresize{\hrule height \Thickness\vss
      \hbox to \Squaresize{\hss#1\hss}
   \vss\hrule height\Thickness}
\unskip\vrule width \Thickness}
\kern-\Thickness}
 
\def\Vsquare#1{\vbox{\Square{$#1$}}\kern-\Thickness}

\def\shuff#1#2{\mathbin{
\hbox{\vbox{ \hbox{\vrule \hskip#2 \vrule height#1 width 0pt}%
\hrule}%
\vbox{ \hbox{\vrule \hskip#2 \vrule height#1 width 0pt
\vrule}%
\hrule}%
}}}

\def\SHUF{{\mathchoice{\shuff{7pt}{3.5pt}}%
{\shuff{6pt}{3pt}}%
{\shuff{4pt}{2pt}}%
{\shuff{3pt}{1.5pt}}}}%

\def\shuffle{\,\SHUF\,\,}

\begin{document}
%%%%%%%%%%%%%%%%

%\today
\CompileMatrices

\bigskip
\title{Semicanonical bases and preprojective algebras}

\author{Christof Gei{ss}}
\address{Christof Gei{ss}\newline
Instituto de Matem\'aticas, UNAM\newline
Ciudad Universitaria\newline
04510 Mexico D.F.\newline
Mexico}
\email{christof@math.unam.mx}

\author{Bernard Leclerc}
\address{Bernard Leclerc\newline
Laboratoire LMNO\newline
Universit\'e de Caen\newline
F-14032 Caen Cedex\newline
France}
\email{leclerc@math.unicaen.fr}

\author{Jan Schr\"oer}
\address{Jan Schr\"oer\newline
Department of Pure Mathematics\newline
University of Leeds\newline
Leeds LS2 9JT\newline
England}
\email{jschroer@maths.leeds.ac.uk}

\dedicatory{
Dedicated to Claus Michael Ringel on the occasion of his sixtieth
birthday.}

\thanks{
C. Geiss acknowledges support from DGAPA-UNAM
for a sabbatical stay at the University of Leeds.
J. Schr\"oer was supported by a research fellowship from 
the DFG (Deutsche Forschungsgemeinschaft).
He also thanks the Laboratoire LMNO (Caen) for an invitation
in Spring 2003 during which this work was started. 
B. Leclerc is grateful to the GDR 2432
(Alg\`ebre non commutative) 
and the GDR 2249 (Groupes, g\'eom\'etrie et repr\'esentations) 
for their support.
}

\bigskip
\begin{abstract}
We study the multiplicative properties of the dual of 
Lusztig's semicanonical basis.
The elements of this basis are naturally indexed by the
irreducible components of Lusztig's nilpotent varieties, which
can be interpreted as varieties of modules over preprojective algebras.
We prove that the product of two dual semicanonical basis vectors
$\rho_{Z'}$ and $\rho_{Z''}$ 
is again a dual semicanonical basis vector provided the closure of
the direct sum of the
corresponding two irreducible components $Z'$ and $Z''$ 
is again an irreducible component.
It follows that the semicanonical basis and the canonical basis
coincide if and only if we are in Dynkin type $\A_n$ with $n \leq 4$.
Finally, we provide a detailed study of the varieties of 
modules over the preprojective
algebra of type $\A_5$.
We show that in this case the multiplicative properties of
the dual semicanonical basis are controlled by the Ringel form of a 
certain tubular algebra of type $(6,3,2)$ and by the
corresponding elliptic root system of type $\E_8^{(1,1)}$. 

\bigskip\noindent
{\sc R\'esum\'e.}  Nous \'etudions les propri\'et\'es multiplicatives
de la base duale de la base semi-canonique de Lusztig.
Les \'el\'ements de cette base sont naturellement param\'etr\'es
par les composantes irr\'educ\-ti\-bles des vari\'et\'es nilpotentes
de Lusztig, qui peuvent \^etre interpr\'et\'ees comme vari\'et\'es
de modules sur les alg\`ebres pr\'eprojectives.
Nous d\'emontrons que le produit de deux vecteurs $\rho_{Z'}$ et
$\rho_{Z''}$ de la base semi-canonique duale est encore un vecteur
de la base semi-canonique duale si la somme directe des composantes
irr\'eductibles $Z'$ et $Z''$ est encore une composante
irr\'eductible.
Il en r\'esulte que les bases canonique et semi-canonique 
ne co\"\i ncident que pour le type de Dynkin $\A_n$ avec $n\le 4$.
Finalement, nous \'etudions en d\'etail les vari\'et\'es de modules
sur l'alg\`ebre pr\'eprojective de type $\A_5$.
Nous montrons que dans ce cas les propri\'et\'es multiplicatives
de la base semi-canonique duale sont controll\'ees par la forme
de Ringel d'une alg\`ebre tubulaire de type $(6,3,2)$ et par le
syst\`eme de racines elliptique de type $\E_8^{(1,1)}$ qui lui
est associ\'e. 
\end{abstract}

%%%%%%%%%%
\maketitle
%%%%%%%%%%

\bigskip
\setcounter{tocdepth}{1}
\tableofcontents

%%%%%%%%%%%%%%%%%%%%%%%%%%%%%%%%%%%%%%%%%%%%%%

\section{Introduction}

%%%%%%%%%%%%%%%%%%%%%%%%%%%%%%%%%%%%%%%%%%%%%%

\subsection{}
Let $\g$ be a simple Lie algebra of simply-laced type $\A, \D, \E$, and 
let $\n$ be a maximal nilpotent subalgebra. 
Let $\B_q$ be the canonical basis  
of the quantum enveloping algebra $U_q(\n)$ \cite{K,Lu90}
and $\B_q^*$ the basis dual to $\B_q$.
When $q$ tends to 1, these two bases specialize to bases 
$\B$ and $\B^*$ of $U(\n)$ and $\C[N]$, respectively.
Here $N$ stands for a maximal unipotent subgroup of a complex simple 
Lie group $G$ with Lie algebra $\g$.

Let $I$ denote an indexing set for the simple roots of $\g$.
Given a finite-dimensional $I$-graded vector space~$\V$ with graded dimension
$|\V|$, we denote by $\Lam_\V$ 
the corresponding nilpotent variety, see \cite[\S 12]{Lu91}.
This variety can be seen as the variety of modules over the
preprojective algebra $\Lam$ attached to the Dynkin diagram of $\g$,
with underlying vector space $\V$ \cite{Ri1}.

For a variety $X$ let $\irr(X)$ be the set of irreducible components of
$X$.
Lusztig has shown that there are natural bijections
\[
\irr(\Lam_\V) \longrightarrow \B_q(|\V|)  \text{ (resp. $\B_q^*(|\V|)$)}
\]
\[
Z \mapsto b_Z \text{ (resp. $b_Z^*$)}
\]
where $B_q(|\V|)$ (resp. $B_q^*(|\V|)$) is the 
subset of $\B_q$ (resp. $\B^*_q$) consisting of the elements of degree $|\V|$.
Kashiwara and Saito \cite{KSa} proved that the crystal basis
of $U_q(\n)$ can be constructed geometrically in terms of these
irreducible components (this was a conjecture of Lusztig).
  
This paper is motivated by several problems about the bases $\B_q$
and $\B_q^*$ and their relations with the varieties $\Lam_\V$ and
the preprojective algebra $\Lam$.

\subsection{}
One problem, which was first considered by Berenstein and Zelevinsky
\cite{BZ}, is to study the multiplicative structure of the basis $\B_q^*$.
Two elements $b^*_1$ and $b^*_2$ of $\B_q^*$ are called 
{\em multiplicative} if their product belongs to $\B_q^*$
up to a power of $q$.
It was conjectured in \cite{BZ} that $b^*_1$ and $b^*_2$
are multiplicative if and only if they $q$-commute.
We refer to this as the BZ-conjecture.
The conjecture was proved for types $\A_2$ and $\A_3$ 
\cite{BZ}, and it also holds for $\A_4$ \cite{Z}.

More recently, Marsh and Reineke observed a strong relationship
between the multiplicative structure of $\B_q^*$ and
properties of the irreducible components of the varieties
$\Lam_\V$.
They checked \cite{MR} that for $\g$ of type $\A_n\ (n\le 3)$,
if the irreducible components $Z_1\subset\Lam_{\V_1}$ and 
$Z_2\subset\Lam_{\V_2}$ 
are the closures of the 
isomorphism classes of two indecomposable $\Lam$-modules $x_1$ and
$x_2$, then $b^*_{Z_1}$ and $b^*_{Z_2}$ are multiplicative if and only if
$\Ext^1_\Lam(x_1,x_2)=0$.
This was verified by a case-by-case calculation,
using the fact that for type $\A_n\ (n\le 3)$ 
the preprojective algebra is of finite representation
type, that is, it has only a finite number of isomorphism classes of 
indecomposable modules \cite{DlRi}.
They also calculated many examples in type $\A_4$ 
and conjectured that this property still holds in this case
(note that $\Lam$ is again representation-finite for $\A_4$).
But a conceptual explanation was still missing.

Let $Z_1 \oplus Z_2$ denote the
subset of $\Lam_{\V_1\oplus\V_2}$ consisting of all $\Lam$-modules $x$
isomorphic to $y_1\oplus y_2$ with $y_1\in Z_1$ and $y_2\in Z_2$.
It follows from a general decomposition theory for irreducible
components of varieties of modules developed in \cite{CBSc}
that the condition $\Ext^1_\Lam(x_1,x_2) = 0$ for some 
$(x_1,x_2) \in Z_1 \times Z_2$ is equivalent to 
$\overline{Z_1\oplus Z_2}$ being an irreducible component of 
$\Lam_{\V_1\oplus\V_2}$.

In \cite{L} counterexamples to the BZ-conjecture were
found for all types other than $\A_n$ with $n\le 4$.
In particular in type~$\A_5$, for a certain $\V$ of dimension 8
one can find an irreducible component $Z$ of $\Lam_\V$
such that
\begin{equation}\label{CE1}
(b^*_Z)^2 = q^{-2}(b^*_{Z'} + b^*_{Z''})
\end{equation}
where $Z'=\overline{Z\oplus Z}$ and $Z''$ are two irreducible
components of $\Lam_{\V\oplus\V}$, see also \cite{GS2}.
This seems to be the smallest counterexample to the BZ-conjecture
in type $\A$.
Moreover, it also shows that the result of Marsh and Reineke does not 
generalize to $\A_5$.
Note however that the BZ-conjecture was proved for large
families of elements of $\B^*_q$ \cite{C1,C2,CalM,LNT}. 
For example, in type $\A$ it holds for quantum flag minors,
and the reformulation in terms of direct sums of irreducible
components is also valid \cite{Sch}. 

So one would like to get a better understanding of the relationship between 
multiplicativity of elements of $\B_q^*$ and direct sum
decompositions of irreducible components of varieties of $\Lam$-modules.

\subsection{}
Another interesting problem concerns the singular supports of 
the simple perverse sheaves used by Lusztig \cite{Lu90} to define the 
canonical basis $\B_q$.
Let $Q$ be a Dynkin quiver, which is obtained from the Dynkin diagram of $\g$
by choosing an orientation.
Let $\rep(Q,\V)$ be the affine space of representations
of $Q$ with underlying finite-dimensional $I$-graded vector space $\V$. 
This is a finite union of isomorphism classes (or orbits) $\orb$.
In Lusztig's geometric construction, the elements of $\B_q(|\V|)$ 
are given by the perverse extensions $L_\orb$
of the constant sheaves $\C_\orb$ on the orbits $\orb$.
In \cite{Lu91} Lusztig considered the singular supports $SS(L_\orb)$
of these sheaves and showed that they are unions of irreducible 
components of $\Lam_\V$ (independent of the chosen orientation of the 
Dynkin diagram of $\g$).
He conjectured that in fact each $SS(L_\orb)$ is irreducible,
equal to the closure $\Lam_\orb$ of the conormal bundle of $\orb$.
Unexpectedly, Kashiwara and Saito \cite{KSa} produced a counterexample
to this conjecture. 
They exhibited two orbits $\orb', \orb''$ for type $\A_5$ such that 
\[
SS(L_{\orb''})=\Lam_{\orb'} \cup \Lam_{\orb''}.
\]
The corresponding vectors $b_{\orb'}$ and $b_{\orb''}$ of $\B_q$ have principal
degree 16, and apparently this is the smallest counterexample
in type $\A$.
   
It turns out that this counterexample is dual to the
counterexample above for $\B_q^*$, in the sense that $\Lam_{\orb'}=Z'$
and $\Lam_{\orb''}=Z''$, see \cite[Remark 1]{L}.
One motivation for this paper was to find an explanation
for this coincidence.

\subsection{}
What makes these problems difficult is that, although the canonical basis
reflects by definition the geometry of the varieties 
$\overline{\orb} \subseteq \rep(Q,\V)$, we want to relate it to the geometry 
of some other varieties, namely the irreducible components of the 
nilpotent varieties $\Lam_\V$.
It is natural to think of an intermediate object, that is, 
a basis reasonably close to the canonical basis, but directly 
defined in terms of the varieties $\Lam_\V$. 
Lusztig \cite{Lu00} has constructed such a basis ${\mathcal S} = \{f_Z\}$ 
and called it the {\em semicanonical basis}. 
This is a basis of $U(\n)$ (not of the $q$-deformation $U_q(\n)$)
which gives rise, like $\B$, to a basis in each irreducible highest 
weight $U(\g)$-module.
Let ${\mathcal S}^*=\{\rho_Z\}$ denote the basis of $\C[N]$ dual 
to ${\mathcal S}$.
Our first main result is the following:

\begin{Thm}\label{ThmA}
If $Z_1\subset \Lam_{\V_1}$ 
and $Z_2\subset\Lam_{\V_2}$ are irreducible components 
such that 
\[
Z=\overline{Z_1\oplus Z_2}
\]
is an irreducible component of $\Lam_{\V_1\oplus\V_2}$, then 
$\rho_{Z_1}\rho_{Z_2}=\rho_Z$.
\end{Thm}
In other words, the dual semicanonical basis ${\mathcal S}^*$ satisfies the
multiplicative property which was expected to hold for the
dual canonical basis $\B^*$.

An irreducible component $Z \in \irr(\Lam_\V)$ is called
{\it indecomposable} if $Z$ contains a dense subset of indecomposable
$\Lam$-modules.
By \cite{CBSc}, every irreducible component $Z$ of $\Lam_\V$
has a {\it canonical decomposition} 
\[
Z = \overline{Z_1\oplus \cdots \oplus Z_m}
\]
where the $Z_i\subset\Lam_{\V_i}$ are indecomposable 
irreducible components.
Our theorem implies that 
\[
\rho_Z = \rho_{Z_1}\cdots\rho_{Z_m}.
\]
Hence ${\mathcal S}^*$ has a natural description as a collection of families
of monomials in the elements indexed by indecomposable
irreducible components.
Such a description of ${\mathcal S}^*$ resembles 
the description of $\B^*$ for type $\A_n\ (n\le 4)$ obtained 
by Berenstein and Zelevinsky.

\subsection{}
So a natural question is how close are the bases ${\mathcal S}^*$ 
and $\B^*$?
In type $\A$, Berenstein and Zelevinsky \cite{BZ} proved that
all minors of the triangular matrix of coordinate functions
on $N$ belong to $\B^*$.
We prove that they also belong to ${\mathcal S}^*$. 
Hence using \cite{LNT, Sch}, it follows that
${\mathcal S}^*\cap\B^*$ contains all 
multiplicative products of flag minors.
However the two bases differ in general. More precisely we have:

\begin{Thm}\label{ThmB}
The bases ${\mathcal S}^*$ and $\B^*$ coincide if and only if $\g$ is of 
type $\A_n$ with $n \le 4$.
\end{Thm}

For example in type $\A_5$, we deduce from 
Equation~(\ref{CE1}) and Theorem~\ref{ThmA} that
\begin{equation}
\rho_{Z'} = b^*_{Z'} + b^*_{Z''}
\end{equation}
(where for simplicity we use the same notation $b^*_{Z'}$ and $b^*_{Z''}$
for the specializations at $q=1$).
Nevertheless, since ${\mathcal S}^*$ and $\B^*$ have lots
of elements in common, we get an explanation why the 
BZ-conjecture (or rather its reformulation in terms
of irreducible components of varieties of $\Lam$-modules)
holds for large families of elements of $\B^*$.

Of course, by duality, these results also allow to compare
the bases ${\mathcal S}$ and $\B$.
In particular, returning to the example of \cite{KSa},
we can check that  
\begin{equation}
b_{\orb''} = f_{\Lam_{\orb'}} + f_{\Lam_{\orb''}},
\end{equation}
and this is probably the smallest example in type $\A$
for which the canonical and semi\-canonical bases differ.
One may conjecture that, in general, the elements $f_Z$ occurring in the 
${\mathcal S}$-expansion of $b_\orb\in\B$ are indexed by the irreducible 
components $Z$
of $SS(L_\orb)$, so that $SS(L_\orb)$ is irreducible
if and only if $b_\orb=f_{\Lam_\orb}$.
(There is a similar conjecture of Lusztig \cite{Lu97} for the ``semicanonical
basis'' of the group algebra of a Weyl group obtained from
the irreducible components of the Steinberg variety.)
Assuming this conjecture we get an explanation of the
relationship between the counterexamples to the conjectures of
Berenstein-Zelevinsky and Lusztig.   

\subsection{}
In the last part of the paper, we consider the first case
which is not well understood, namely type $\A_5$.
In this case, the preprojective algebra $\Lam$ is
representation-infinite, 
but it is still of tame representation type \cite{DlRi}.
Motivated by our description of ${\mathcal S}^*$ in terms of indecomposable
irreducible components of varieties of $\Lam$-modules, we 
give a classification of the indecomposable irreducible components 
for the case $\A_5$.
We also give an explicit criterion to decide when the closure of the 
direct sum of 
two such components is again an irreducible component.
These results are deduced from \cite{GSc}, in which a general
classification of irreducible components of varieties of modules
over tubular algebras is developed. 
They are naturally formulated in terms of the Ringel bilinear
form $\<-,-\>$ of a convex subalgebra $\Da$ of a Galois covering 
of $\Lam$.
The algebra $\Da$ is a tubular algebra of type $(6,3,2)$
and the corresponding 10-dimensional infinite root system~$R$
is an elliptic root system of type $\E_8^{(1,1)}$
in the classification of Saito \cite{Sa},
with a 2-dimensional lattice of imaginary roots. 
Note that the irreducible component $Z$ of Equation~(\ref{CE1})
corresponds to a generator of this lattice.
(This is an a posteriori justification for calling $b^*_Z$ 
an imaginary vector in \cite{L}.)
The Ringel form $\<-,-\>$ allows to define a distinguished
Coxeter matrix $\Phi$ of order $6$ acting on $R$.
We prove the following:

\begin{Thm}\label{ThmC}
There is a one-to-one correspondence 
$r\mapsto Z(r)$ between the set 
of Schur roots of $R$ and the set of indecomposable irreducible
components of the nilpotent varieties of type $\A_5$ which
do not contain an indecomposable projective $\Lam$-module.
Moreover $\overline{Z(r_1)\oplus Z(r_2)}$ is an
irreducible component if and only if the Schur roots
$r_1$ and $r_2$ satisfy certain conditions which are
all expressible in terms of $\<-,-\>$ and $\Phi$.
\end{Thm}
We also explain how to translate from the language of 
roots to the language of multisegments,
which form a natural indexing set of canonical and semicanonical
bases in type~$\A$.

\subsection{}
The paper is organized as follows.
In Section~\ref{var} we recall the general theory of varieties of
modules.
We explain a general decomposition theory for irreducible components
of such varieties.
This is followed in Section~\ref{prepro} by a short introduction to
preprojective algebras. 
Then we recall the concept of a constructible function in Section 
\ref{constf}.
Following Lusztig~\cite{Lu00}, we review in Section~\ref{semi}
the definition of the semicanonical basis of $U(\n)$, 
which is obtained by realizing $U(\n)$ as an algebra $\M$ of
constructible functions on the nilpotent varieties.
In order to study the dual semicanonical basis and its
multiplicative properties we also need to describe  
the natural comultiplication of $U(\n)$ in terms of $\M$.
This was not done in \cite{Lu00}, so we provide this 
description in Section~\ref{comult}.
In Section~\ref{multipli} we introduce the dual
semicanonical basis ${\mathcal S}^*$ of $\M^*$ and prove Theorem~\ref{ThmA}.
Note that for this theorem we do not restrict ourselves to types $\A, \D, \E$,
and only assume that $\n$ is the positive part of a symmetric
Kac-Moody Lie algebra. 
We end this section with the proof of the ``only if'' part of 
Theorem~\ref{ThmB}.
In Section~\ref{emb} we embed $\M^*$ into the shuffle algebra.
This gives a practical way of computing elements of ${\mathcal S}^*$.
We use this to prove that in type $\A$ all nonzero minors in the
coordinate functions of $N$ belong to ${\mathcal S}^*$.
In the rest of the paper we focus on the Dynkin cases $\A_n\ (n\le 5)$.
In Section~\ref{Galois} we consider a Galois covering $\tLam$ of the 
algebra $\Lam$, with Galois group $\Z$, and we use it to calculate
the Auslander-Reiten quiver of $\Lam$ for $n\le 4$.
We also introduce an algebra $\Da$ whose repetitive algebra
is isomorphic to $\tLam$. For $n\le 4$, $\Da$ has finite 
representation type, while for $n=5$ it is a tubular algebra of
tubular type $(6,3,2)$. 
In Section~\ref{indmult} we recall from \cite{GSc} 
that the indecomposable irreducible components of $\Lam$ are in 
one-to-one correspondence with the $\Z$-orbits of Schur roots of $\tLam$.
We also describe the map which associates to such a Schur root
the multisegment indexing the corresponding indecomposable
irreducible component.
The component graphs for the representation-finite cases
$\A_2$, $\A_3$ and $\A_4$ are constructed in Section~\ref{comp234},
and the corresponding graphs of prime elements of ${\mathcal B}^*$
are described in Section~\ref{prime234}.
In Section~\ref{compare} we prove the ``if'' part of Theorem~\ref{ThmB}.
All the remaining sections are devoted to the case $\A_5$.
In Section~\ref{tubweight} we relate the category of 
$\tLam$-modules to the category $\md(\Da)$ of modules over
the tubular algebra $\Da$ and to the
category $\coh(\X)$ of coherent sheaves on a weighted projective line
$\X$ of type $(6,3,2)$ in the sense of Geigle and Lenzing \cite{GL}.
In Section~\ref{roots} we consider the Grothendieck groups
$K_0(\md(\Da))\simeq K_0(\coh(\X))\simeq \Z^{10}$. 
They are naturally endowed with a (non-symmetric) bilinear 
form $\bil{-,-}$ (the Ringel form) and a Coxeter matrix.
This gives rise to an elliptic root system of type $\E_8^{(1,1)}$.
We give an explicit description of its set of positive roots
and of the subset $R^+_S$ of Schur roots.
In Section~\ref{category}, we show that $R^+_S$ naturally parametrizes
the $\Z$-orbits of Schur roots of $\tLam$, hence also
the indecomposable irreducible components of $\Lam$.
Then Section~\ref{extensions} describes the component
graph of $\Lam$ for type $\A_5$, thus making precise the
statements of Theorem~\ref{ThmC}.
Section~\ref{proof1} consists of the proof of Theorem~\ref{Thm1}.
We conclude by noting the existence of similar results
for type $\D_4$ and by pointing out some possible connections
with the theory of cluster algebras of Fomin and Zelevinsky
(Section~\ref{concl}).
Section~\ref{pictures} contains a collection of pictures and
tables to which we refer at various places in the text.

%%%%%%%%%%%%%%%%%%%%%%%%%%%%%%%%%%%%%%%%%%%%%%%%%%
\subsection{} %  Conventions
%%%%%%%%%%%%%%%%%%%%%%%%%%%%%%%%%%%%%%%%%%%%%%%%%%
Throughout, we use the following conventions.
If $f: M_1 \to M_2$ and $g: M_2 \to M_3$ are
maps, then the composition is denoted by $gf: M_1 \to M_3$.
Similarly, if $\alpha: 1 \to 2$ and $\beta: 2 \to 3$ are arrows in
a quiver, then the composition of $\alpha$ and $\beta$ is denoted
by $\beta\alpha$.

Modules are always assumed to be left modules. 

All vector spaces are over the field $\C$ of complex numbers.

We set $\Q^+ = \{ q \in \Q \mid q > 0 \}$, 
$\Q^- = \{ q \in \Q \mid q < 0 \}$ and $\Q_\infty = \Q \cup \{ \infty
\}$.
We also set $\Z^+=\{z\in\Z\mid z>0\}$ and $\ENn = \Z^+\cup\{0\}$.

%%%%%%%%%%%%%%%%%%%%%%%%%%%%%%%%%%%%%%%%%%%%%%%%%%%%%%%%%%%%%%%%%

\section{Varieties of modules}\label{var}

%%%%%%%%%%%%%%%%%%%%%%%%%%%%%%%%%%%%%%%%%%%%%%%%%%%%%%%%%%%%%%%%%
%%%%%%%%%%%%%%%%%%%%%%%%%%%%%%%%%%%%
\subsection{} % Quivers and algebras
%%%%%%%%%%%%%%%%%%%%%%%%%%%%%%%%%%%%
A {\it quiver} is a quadruple 
$Q = (I,Q_1,s,e)$
where $I$ and $Q_1$ are
sets with $I$ non-empty, and $s,e: Q_1 \to I$ are maps such that
$s^{-1}(i)$ and $e^{-1}(i)$ are finite for all $i \in I$.
We call $I$ the set of {\it vertices} and $Q_1$ the set of {\it arrows}
of $Q$.
For an arrow $\alpha \in Q_1$ one calls $s(\alpha)$ the starting vertex and
$e(\alpha)$ the end vertex of $\alpha$.

A {\it path} of length $t$ in $Q$ is a sequence 
$p = \alpha_1 \alpha_2 \cdots \alpha_t$ of arrows such that
$s(\alpha_i) = e(\alpha_{i+1})$ for $1 \leq i \leq t-1$.
Set $s(p) = s(\alpha_t)$ and $e(p) = e(\alpha_1)$.
Additionally, for each vertex $i \in I$ let $1_i$ be a 
path of length 0.
By $\C Q$ we denote the {\it path algebra} of $Q$, with
basis the set of all paths in $Q$ and product given
by concatenation.
A ${\it relation}$ for $Q$ is a linear combination
\[
\sum_{i = 1}^t \lam_i p_i
\]
where $\lam_i \in \C$ and the $p_i$ are paths of length at
least two in $Q$ with
$s(p_i) = s(p_j)$ and $e(p_i) = e(p_j)$ for all $1 \leq i,j \leq t$.
Thus, we can regard a relation as an element in $\C Q$.

An ideal $J$ in $\C Q$ is {\it admissible} if it is 
generated by a set of relations for $Q$.
Note that this differs from the usual definition of an admissible ideal,
where one also assumes that the factor algebra $\C Q/J$ is finite-dimensional.
%%%%%%%%%%%%%%%%%%%%%%%%%%%%%%%%%%%%%
\subsection{} % Representations
%%%%%%%%%%%%%%%%%%%%%%%%%%%%%%%%%%%%%
A map $d: I \to \ENn$ such that $I \setminus d^{-1}(0)$ is finite
is called a {\it dimension vector for} $Q$.
We also write $d_i$ instead of $d(i)$, and we often use the notation
$d = (d_i)_{i \in I}$.
By $\ENn^{(I)}$ we denote the semigroup of dimension vectors for $Q$.

Let ${\mathcal V}_{\rm fin}(I)$ be the category of finite-dimensional
$I$-graded vector spaces. 
Thus, the objects of ${\mathcal V}_{\rm fin}(I)$ are of the form
$\V = \bigoplus_{i \in I} V_i$ where the $V_i$ are finite-dimensional
vector spaces, and only finitely many of the $V_i$ are nonzero.
We call $|\V| = (\dm(V_i))_{i \in I}$ the {\it dimension vector} or
{\it degree} of $\V$.
The morphisms in ${\mathcal V}_{\rm fin}(I)$ are just linear maps
respecting the grading.
Direct sums in  ${\mathcal V}_{\rm fin}(I)$ are defined in the obvious way.

A {\it representation} of $Q$ with underlying vector space 
$\V \in {\mathcal V}_{\rm fin}(I)$ is
an element 
\[
x = (x_\alpha)_{\alpha \in Q_1} \in \rep(Q,\V) =
\bigoplus_{\alpha \in Q_1} \Hom_\C(V_{s(\alpha)},V_{e(\alpha)}).
\]
For a representation $x = (x_\alpha)_{\alpha \in Q_1} \in \rep(Q,\V)$ 
and a path 
$p = \alpha_1 \alpha_2 \cdots \alpha_t$ in $Q$ set
\[
x_p = x_{\alpha_1} x_{\alpha_2} \cdots x_{\alpha_t}.
\]
Then $x$ {\it satisfies a relation}
$\sum_{i = 1}^t \lam_i p_i$ 
if
$\sum_{i = 1}^t \lam_i x_{p_i} = 0$.
If $R$ is a set of relations for $Q$, then let
$
\rep(Q,R,\V)
$
be the set of all representations $x \in \rep(Q,\V)$ which
satisfy all relations in $R$.
This is a closed subvariety of $\rep(Q,\V)$.
Let $A$ be the algebra $\C Q/J$, where $J$ is the admissible
ideal generated by $R$.
Note that every element in $\rep(Q,R,\V)$ can be naturally interpreted as an 
$A$-module structure on $\V$, so we shall also write 
\[
\md(A,\V) = \rep(Q,R,\V).
\]
This is the affine variety of $A$-modules with underlying vector space $\V$.
A {\it dimension vector for} $A$ is by definition the same as a 
dimension vector for $Q$, that is, an element of $\ENn^{(I)}$.
For $x \in \md(A,\V)$ we call 
$
\dimv{x} = |\V|
$
the {\it dimension vector} of $x$.

%%%%%%%%%%%%%%%%%%%%%%%%%%%%%%%%%%%%%%%%%%%%%
\subsection{} % orbits
%%%%%%%%%%%%%%%%%%%%%%%%%%%%%%%%%%%%%%%%%%%%%
Define
$
G_\V = \prod_{i \in I} \GL(V_i). 
$
This algebraic group acts on $\md(A,\V)$ as follows.
For $g = (g_i)_{i \in I} \in G_\V$ and 
$x = (x_\alpha)_{\alpha \in Q_1} \in \md(A,\V)$ define
\[
g \cdot x = (x'_\alpha)_{\alpha \in Q_1}
\text{ where }
x'_\alpha = g_{e(\alpha)} x_\alpha g_{s(\alpha)}^{-1}.
\]
The $G_\V$-orbit of an $A$-module $x \in \md(A,\V)$ 
is denoted by $\orb(x)$.
Two $A$-modules $x,y \in \md(A,\V)$ are isomorphic if and only if they
lie in the same orbit.

For a dimension vector $d$ for $A$ set
\[
\V^d  = \bigoplus_{i \in I} \C^{d_i}, \quad
\md(A,d) = \md(A,\V^d), \quad
\GL(d) = G_{\V^d}.
\]
Thus 
$
\md(A,|\V|) \cong \md(A,\V)
$
for all $\V \in {\mathcal V}_{\rm fin}(I)$.
For this reason, we often do not distinguish between 
$\md(A,|\V|)$ and $\md(A,\V)$.
Any problems arising from this can be solved via the existence of an
isomorphism between these two varieties which respects the group
actions and the gradings.

%%%%%%%%%%%%%%%%%%%%%%%%%%%%%%%%%%%%%%%%%%%%%
\subsection{} % nilpotent orbits
%%%%%%%%%%%%%%%%%%%%%%%%%%%%%%%%%%%%%%%%%%%%%

By abuse of notation, we identify $i \in I$ with the dimension
vector mapping $i$ to 1 and $j\not = i$ to $0$.
If $\V$ is an $I$-graded vector space with $|\V| = i$, then
the variety $\md(A,\V)$ consists just of a single point and is
denoted by $Z[i]$.
The corresponding 1-dimensional $A$-module is denoted by $S_i$.

An element $x \in \md(A,\V)$ is said to be {\it nilpotent} if there exists
an $N \in \ENp$ such that for any path
$p$ of length greater than $N$ we have
$x_p = 0$.
By $\md^0(A,\V)$ we denote the closed subset of nilpotent elements in 
$\md(A,\V)$.
The nilpotent elements are exactly the $A$-modules whose composition series
contains only factors isomorphic to the simple modules $S_i$, 
$i \in I$.

%%%%%%%%%%%%%%%%%%%%%%%%%%%%%%%%%%%%%%%%%%%%%
\subsection{} % Irreducible components, notation
%%%%%%%%%%%%%%%%%%%%%%%%%%%%%%%%%%%%%%%%%%%%%
An irreducible component $Z \in \irr(\md(A,\V))$ is called 
{\it indecomposable}
if it contains a dense subset of indecomposable $A$-modules.
Let $\ind(\irr(\md(A,\V)))$ be the set of indecomposable irreducible
components of $\md(A,\V)$.
Let 
\[
\irr(A) = \bigcup_{d \in \ENn^{(I)}} \irr(\md(A,d))
\]
be the set of all irreducible components of varieties of
$A$-modules, and set
\[
\ind(\irr(A)) =  \bigcup_{d \in \ENn^{(I)}} \ind(\irr(\md(A,d))).
\]

%%%%%%%%%%%%%%%%%%%%%%%%%%%%%%%%%%%%%%%%%%%%%
\subsection{} % Irreducible components, canonical decomposition
%%%%%%%%%%%%%%%%%%%%%%%%%%%%%%%%%%%%%%%%%%%%%
Let $Z_i \in \irr(\md(A,\V^i))$, $1 \le i \le m$ be irreducible
components of $A$-modules, and set $\V = \V^1 \oplus \cdots \oplus \V^m$.
Let
$
Z_1 \oplus \cdots \oplus Z_m
$
be the set of modules in $\md(A,\V)$, which are isomorphic
to $x_1 \oplus \cdots \oplus x_m$ with $x_i \in Z_i$ for all $i$.
The closure
$
\overline{Z_1 \oplus \cdots \oplus Z_m}
$
is irreducible, but in general it is not an irreducible component.
According to \cite{CBSc}
any irreducible component $Z \in \irr(A)$
has a decomposition
\[
Z = \overline{Z_1 \oplus \cdots \oplus Z_m}
\]
for certain indecomposable irreducible components 
$Z_i \in \ind(\irr(A))$.
Moreover the components $Z_1,\ldots ,Z_m$ are uniquely 
determined up to reordering. 
This is called the {\it canonical decomposition} of $Z$.

For irreducible components $Z'$ and $Z''$ of $A$-modules define
\[
\ext^1_A(Z',Z'') = 
\min \{ \dm \Ext_A^1(x',x'')\mid (x',x'') \in Z'\times Z'' \}.
\]
This is the dimension of the extension group $\Ext_A^1(x',x'')$
for all $(x',x'')$ in a certain dense open subset of $Z'\times Z''$.
For irreducible components $Z_i \in \irr(A)\ (1 \le i \le m)$
it is known that
$
\overline{Z_1\oplus\cdots\oplus Z_m}
$
is an irreducible component if and only if
$\ext^1_A(Z_i,Z_j) = 0$ for all $i\not = j$, see \cite{CBSc}.

The {\it component graph} ${\mathcal C}(A)$ of $A$ is defined as follows.
The vertices of ${\mathcal C}(A)$ are the elements in $\ind(\irr(A))$.
There is an edge between vertices $Z'$ and $Z''$ if and only if
$\overline{Z' \oplus Z''}$ is again an irreducible component, or equivalently
if 
\[
{\rm ext}_A^1(Z',Z'') = {\rm ext}_A^1(Z'',Z') = 0. 
\]

%%%%%%%%%%%%%%%%%%%%%%%%%%%%%%%%%%%%%%%%%%%%%%%%%%%%%%%%%%%%%%%%%%

\section{Preprojective algebras}\label{prepro}

%%%%%%%%%%%%%%%%%%%%%%%%%%%%%%%%%%%%%%%%%%%%%%%%%%%%%%%%%%%%%%%%%%
%%%%%%%%%%%%%%%%%%%%%%%%%%%%%%%%%%%%%%%%%%%%%%%%%%%%%%%%%%%
\subsection{} % Definitions
%%%%%%%%%%%%%%%%%%%%%%%%%%%%%%%%%%%%%%%%%%%%%%%%%%%%%%%%%%%
Assume that $Q = (I,Q_1,s,e)$ is a finite quiver without loops.
(A {\it loop} is an arrow $\alpha$ with $s(\alpha) = e(\alpha)$,
and $Q$ is {\it finite} if $I$ is finite.
Note that this implies that $Q_1$ is finite as well.)
The {\it double quiver} 
$
\overline{Q} = (I,\overline{Q_1},\overline{s},\overline{e})
$
of $Q$ is obtained from $Q$ by adding for each arrow $\alpha \in Q_1$ 
an additional arrow $\alpha^*$.
Define $\overline{s}(\alpha) = s(\alpha)$,
$\overline{e}(\alpha) = e(\alpha)$,
$\overline{s}(\alpha^*) = e(\alpha)$ and
$\overline{e}(\alpha^*) = s(\alpha)$.
For any $i \in I$ let
\[
r_i = \sum_{\alpha \in Q_1: s(\alpha) = i} \alpha^* \alpha 
- \sum_{\alpha \in Q_1: e(\alpha) = i} \alpha \alpha^*
\]
be the {\it Gelfand-Ponomarev relation} associated to $i$.
This is a relation for $\overline{Q}$.
The {\it preprojective algebra} $P(Q)$ corresponding to $Q$ is defined as
\[
P(Q) = \C \overline{Q}/J
\]
where $J$ is generated by the relations $r_i$, $i \in I$.
These algebras were introduced and studied by Gelfand and Ponomarev,
compare also \cite{Ri1}.

%%%%%%%%%%%%%%%%%%%%%%%%%%%%%%%%%%%%%%%%%%%%%%%%%%%
\subsection{} % Nilpotent varieties
%%%%%%%%%%%%%%%%%%%%%%%%%%%%%%%%%%%%%%%%%%%%%%%%%%%
For a preprojective algebra $\Lam = P(Q)$ and some 
$\V \in {\mathcal V}_{\rm fin}(I)$ set
\[
\Lam_\V = \md^0(P(Q),\V).
\]
Lusztig proved that $\Lam_\V$ has pure dimension 
$\dm(\rep(Q,\V))$, i.e.
all irreducible components of $\Lam_\V$ have dimension
$\dm(\rep(Q,\V))$.
Usually the varieties $\Lam_\V$ are called {\it nilpotent varieties}.
If $Q$ is a Dynkin quiver, then one might call them
just {\it preprojective varieties}, since the nilpotency condition
follows automatically in these cases, as shown by the next proposition.
\begin{Prop}
For a preprojective algebra $P(Q)$ the following are equivalent:
\begin{itemize}
\item[(a)]
$P(Q)$ is finite-dimensional;
\item[(b)]
$\md^0(P(Q),\V) = \md(P(Q),\V)$ for all $\V \in {\mathcal V}_{\rm fin}(I)$; 
\item[(c)]
$Q$ is a Dynkin quiver.
\end{itemize}
\end{Prop}
The equivalence of (a) and (c) is classical (see for instance
\cite{Rei}).
The implication (c) $\Longrightarrow$ (b) is proved by Lusztig
\cite[14.2]{Lu91}, and the converse by Crawley-Boevey 
\cite{CB2}. 

%%%%%%%%%%%%%%%%%%%%%%%%%%%%%%%%%%%%%%%%%%%%%%%%%%%%%%%%%%%%%
\subsection{} %known results
%%%%%%%%%%%%%%%%%%%%%%%%%%%%%%%%%%%%%%%%%%%%%%%%%%%%%%%%%%%%%
The following remarkable property of preprojective algebras 
was proved in \cite{CB}.
\begin{Prop}\label{PropCB}
For finite-dimensional modules $X$ and $Y$ over a preprojective
algebra $\Lam$ we have
\[
\dm \Ext^1_{\Lam}(X,Y) = \dm \Ext^1_{\Lam}(Y,X).
\]
\end{Prop}

Most preprojective algebras are of wild representation type.
The following proposition lists the exceptional cases.
We refer to \cite{DlRi} and \cite{GSc} for further details.
\begin{Prop}\label{infinite_type}
Let $P(Q)$ be a preprojective algebra.
Then the following hold:
\begin{itemize}
\item[(a)]
$P(Q)$ is of finite representation type
if and only if $Q$ is of Dynkin type $\A_n$ with $n \leq 4$;
\item[(b)]
$P(Q)$ is of tame representation type 
if and only if $Q$ is of Dynkin type $\A_5$ or $\D_4$.
\end{itemize}
\end{Prop}

%%%%%%%%%%%%%%%%%%%%%%%%%%%%%%%%%%%%%%%%%%%%%%%%%%%%%%%%%%%%%%%%%

\section{Constructible functions}\label{constf}

%%%%%%%%%%%%%%%%%%%%%%%%%%%%%%%%%%%%%%%%%%%%%%%%%%%%%%%%%%%%%%%%%

%%%%%%%%%%%%%%%%%%%%%%%%%%%%%%%%%%%%%%%%%%%%%%%%%%%%
\subsection{} % constructible sets and functions
%%%%%%%%%%%%%%%%%%%%%%%%%%%%%%%%%%%%%%%%%%%%%%%%%%%%
Let $X$ be an algebraic variety over~$\C$, endowed with
its Zariski topology.
A subset $A$ of $X$ is said to be constructible if it is
a finite union of locally closed subsets.
It is easy to see that if $X$ is irreducible and if we
have a partition
$X=A_1\sqcup \cdots \sqcup A_k$ 
into a finite number of constructible subsets,
then there exists a unique $A_i$ containing a dense open subset of $X$.

A function $f: X \ra \C$ is {\it constructible} if it is a finite
$\C$-linear combination of characteristic functions $\1_A$ for 
various constructible subsets $A$.
Equivalently, $f(X)$ is finite and $f^{-1}(m)$ is a constructible
subset of $X$ for all $m\in\C$.
The set of constructible functions on $X$ is denoted by $M(X)$.
This is a $\C$-vector space. 

%%%%%%%%%%%%%%%%%%%%%%%%%%%%%%%%%%%%%%
\subsection{} % Definition of $\rho_Z$
\label{rhoZ}
%%%%%%%%%%%%%%%%%%%%%%%%%%%%%%%%%%%%%%
If $Z$ is an irreducible component of $X$ and $f\in M(X)$, then
\[
Z=\bigsqcup_{m\in\C} (Z\cap f^{-1}(m))
\]
is a finite partition into constructible subsets, 
hence there is a unique $m\in\C$
such that $Z\cap f^{-1}(m)$ contains a dense open subset
of $Z$.
In other words, a constructible function $f$ has to be
constant on a dense open subset of each irreducible 
component $Z$ of $X$. 
We denote by 
\[
\rho_Z : M(X) \ra \C
\]
the linear form associating
to $f$ its value on this dense open subset of $Z$.

%%%%%%%%%%%%%%%%%%%%%%%%%%%%%%%%%%%%%%%%%%%%%%%%%%%
\subsection{} % Integral
%%%%%%%%%%%%%%%%%%%%%%%%%%%%%%%%%%%%%%%%%%%%%%%%%%%
By $\chi(A)$ we denote the Euler characteristic of a constructible 
subset $A$.
If $A \cap B = \emptyset$, then  
$ \chi(A\sqcup B) = \chi(A) + \chi(B)$.
Hence $\chi$ can be regarded 
as a ``measure'' on the set of constructible subsets of $X$.
For $f\in M(X)$ it is then natural to define
\[
\int_{x\in X} f(x) = \sum_{m\in\C} m\, \chi(f^{-1}(m)) \in \C.
\]
This is a linear form on $M(X)$.
More generally, for a constructible subset $A$ of $X$ we write
\[
\int_{x\in A} f(x) = \sum_{m\in\C} m\, \chi(f^{-1}(m)\cap A).
\]

%%%%%%%%%%%%%%%%%%%%%%%%%%%%%%%%%%%%%%%%%%%%%%%%%%%%%%%%%%%%%%%%%

\section{Semicanonical bases}\label{semi}

%%%%%%%%%%%%%%%%%%%%%%%%%%%%%%%%%%%%%%%%%%%%%%%%%%%%%%%%%%%%%%%%%

%%%%%%%%%%%%%%%%%%%%%%%%%%%%%%%%%%%%%%%%%%%%%%%
\subsection{} % Notation
%%%%%%%%%%%%%%%%%%%%%%%%%%%%%%%%%%%%%%%%%%%%%%%
In this section we assume that $Q = (I,Q_1,s,e)$ is a finite quiver 
without loops,
and as before for $\V \in {\mathcal V}_{\rm fin}(I)$ let 
$\Lam_\V = \md^0(P(Q),\V)$ be the corresponding nilpotent variety.
We denote by $M(\Lam_\V)^{G_\V}$ the subspace of $M(\Lam_\V)$ consisting
of the constructible functions which are constant on
the orbits of $G_\V$.

%%%%%%%%%%%%%%%%%%%%%%%%%%%%%%%%%%%%%%%%%%%%%%%
\subsection{} % Definition of the star product
%%%%%%%%%%%%%%%%%%%%%%%%%%%%%%%%%%%%%%%%%%%%%%%
For $\V, \V', \V'' $ such that $|\V'|+|\V''|=|\V|$, Lusztig
\cite{Lu00} defines a bilinear map 
\[
*:\ M(\Lam_{\V'})^{G_{\V'}}\times M(\Lam_{\V''})^{G_{\V''}} \lra 
M(\Lam_\V)^{G_\V}
\]
as follows.
Let $x\in\Lam_\V$. 
Define $V_x$ to be the variety of all $I$-graded subspaces $\U$ 
of $\V$ such that $x(\U) \subseteq \U$ and $|\U|=|\V''|$.
In other words, $V_x$ is the variety of all $\Lam$-submodules
of $x$ with dimension vector $|\V''|$.
For such a $\U$ let $x'\in \Lam_{\V/\U}$ and $x''\in\Lam_\U$ 
be the elements 
induced by $x$, and let $\widetilde{x}'\in\Lam_{\V'}$ and 
$\widetilde{x}''\in\Lam_{\V''}$
be the elements obtained by transporting $x'$ and $x''$ via
some isomorphisms 
$\V/\U \stackrel{\sim}{\lra} \V'$ 
and 
$\U\stackrel{\sim}{\lra}\V''$.
For 
$f'\in M(\Lam_{\V'})^{G_{\V'}}$ 
and 
$f''\in M(\Lam_{\V''})^{G_{\V''}}$ 
define 
$
\phi_x(f',f''): V_x \ra \C
$
by
\[ 
\phi_x(f',f'')(\U) = f'(\widetilde{x}')f''(\widetilde{x}'').
\]
Following \cite{Lu00} define
\[
(f' * f'')(x) = \int_{\U\in V_x} \phi_x(f',f'')(\U).
\]

%%%%%%%%%%%%%%%%%%%%%%%%%%%%%%%%%%%%%%%%%%%%%%%
\subsection{} % star of two orbits
%%%%%%%%%%%%%%%%%%%%%%%%%%%%%%%%%%%%%%%%%%%%%%%
Let $\orb'\subset\Lam_{\V'}$ (resp. $\orb''\subset\Lam_{\V''}$) be   
a $G_{\V'}$-orbit (resp. a $G_{\V''}$-orbit). 
For $x\in\Lam_\V$ let $\F(\orb',\orb'',x)$ be the constructible
subset of $V_x$ consisting of all $\Lam$-submodules $y$ of $x$ of 
isomorphism type $\orb''$ and such that $x/y$ has isomorphism
type $\orb'$.
The above definition yields
\[
(\1_{\orb'}*\1_{\orb''})(x) = \chi(\F(\orb',\orb'',x)).
\] 
Note that in general the variety $\Lam_\V$ has infinitely
many orbits.
(Indeed, by Proposition~\ref{infinite_type} the algebra $\Lam$ has in general infinite
representation type, hence, by the validity of the second
Brauer-Thrall conjecture (see \cite{Bau}), it has in 
general an infinite number of non-isomorphic representations
of a given dimension).
Therefore the support of the function $\1_{\orb'}*\1_{\orb''}$
may consist of an infinite number of orbits.

%%%%%%%%%%%%%%%%%%%%%%%%%%%%%%%%%%%%%%%%%%%%%%%
\subsection{} % Definition of the algebra $\M$
%%%%%%%%%%%%%%%%%%%%%%%%%%%%%%%%%%%%%%%%%%%%%%%
Let 
$
\widetilde{\M} = \bigoplus_\V M(\Lam_\V)^{G_\V}, 
$
where $\V$ runs over the 
set of all isomorphism classes of vector spaces in 
${\mathcal V}_{\rm fin}(I)$.
(For example, we can take 
$\V \in \{ \V^d \mid d \in \ENn^{(I)} \}$.)
The operation $*$ defines the structure of an $\ENn^{(I)}$-graded associative 
$\C$-algebra on $\widetilde{\M}$.

For $i \in I$, we recall that $Z[i]$ denotes the variety $\Lam_\V$ where
$|\V| = i$. ($Z[i]$ is just a single point.)
Following \cite{Lu00} define
$\M$ to be the subalgebra of $(\widetilde{\M},*)$ generated by
the functions $\1_{Z[i]}$, $i \in I$.  
We set 
$
\M_\V = \M \cap M(\Lam_\V)^{G_\V}.
$

%%%%%%%%%%%%%%%%%%%%%%%%%%%%%%%%%%%%%%%%%%%%%%
\subsection{} %The KM algebra
%%%%%%%%%%%%%%%%%%%%%%%%%%%%%%%%%%%%%%%%%%%%%%
For two distinct vertices $i,j \in I$, let $a_{ij}$
denote the negative of the number of arrows $\alpha \in Q_1$
such that $\{ s(\alpha), e(\alpha) \} = \{ i,j \}$.
Set also $a_{ii}=2\ (i\in I)$.
Let $\g$ be the symmetric Kac-Moody Lie algebra over $\C$ 
with Cartan matrix $(a_{ij})_{i,j\in I}$. 
Let $\n$ be a maximal nilpotent subalgebra of $\g$, and
let $U(\n)$ be its enveloping algebra. 
We denote by $e_i\ (i\in I)$ the Chevalley generators of $U(\n)$.
The defining relations of $U(\n)$ are
\[
\sum_{k=0}^{-a_{ij}+1}
(-1)^k e_i^{(k)} e_j\, e_i^{(-a_{ij}+1-k)} =0,
\]
where $e_i^{(k)}=e_i^k/k!$.

The algebra $U(\n)$ is $\ENn^{(I)}$-graded by assigning to $e_i$
the degree $i$.
It is known that the number of irreducible components of Lusztig's nilpotent
variety $\Lam_\V$ is equal to the dimension of the homogeneous
component of $U(\n)$ of degree $|\V|$. 
This was proved by Lusztig \cite{Lu91,Lu92} when $\g$ is of finite or affine
type and by Kashiwara and Saito \cite{KSa} in general.

%%%%%%%%%%%%%%%%%%%%%%%%%%%%%%%%%%%%%%%%%%%%%%%%%%%%
\subsection{} % Isomorphism and semicanonical basis
%%%%%%%%%%%%%%%%%%%%%%%%%%%%%%%%%%%%%%%%%%%%%%%%%%%%
Lusztig has proved that there is an algebra isomorphism  
\[
\Psi : U(\n) \ra \M
\]
given by 
$\Psi(e_i)=\1_{Z[i]}$.
To do this he constructed for every $I$-graded vector space $\V$ 
a $\C$-basis 
\[
\{ f_Z \mid Z \in \irr(\Lam_\V) \}
\]
of $\M_\V$, naturally labelled by
the irreducible components of $\Lam_\V$.
Using the notation of Section~\ref{rhoZ}, it is characterized by 
\begin{equation}\label{caract}
\rho_{Z'}(f_Z) = \delta_{Z,Z'}, \qquad (Z,Z'\in\irr(\Lam_\V)). 
\end{equation}
In other words, the function $f_Z$ is the unique element of $\M_\V$
equal to $1$ on a dense
open subset of $Z$ and equal to $0$ on a dense open subset of any
other irreducible component $Z'$ of $\Lam_\V$, 
see \cite[Lemma 2.5]{Lu00}.

The basis of $U(\n)$ obtained by transporting via $\Psi^{-1}$
the collection 
\[
\bigcup_\V\,\{ f_Z \mid Z \in \irr(\Lam_\V) \},
\] 
where $\V$ ranges over the set of all isomorphism classes 
of vector spaces in ${\mathcal V}_{\rm fin}(I)$,
is called the {\it semicanonical basis} of $U(\n)$ and is
denoted by ${\mathcal S}$.

\begin{example}\label{ex1}
{\rm
Let $Q$ be the quiver with two vertices $1$ and $2$ and one arrow
$\alpha: 1 \to 2$.
Thus $Q$ is a Dynkin quiver of type $\A_2$.
 
(a)\quad 
Let $\V = V_1 \oplus V_2$ with $\dm(V_1) = \dim{V_2} =1$.
Then
\[
\Lam_\V = \{ x=(a,b) \in \C \times \C \mid ab=0 \}.
\]
The variety $\Lam_\V$ has two irreducible components
\[
Z = \{ (a,0) \mid a \in \C \}, \qquad Z' = \{ (0,b) \mid b \in \C \}.
\]
The group $G_\V = \C^* \times \C^*$ acts on $\Lam_\V$ with three orbits 
$\{(0,0)\}, Z - \{(0,0)\}, Z'- \{(0,0)\}$.
We have 
\[
f_Z = \1_Z = \1_{Z[1]} * \1_{Z[2]}, \qquad
f_{Z'} = \1_{Z'} = \1_{Z[2]} * \1_{Z[1]}.
\]  

(b)\quad Let $\V = V_1 \oplus V_2$ with $\dm(V_1) = \dm(V_2) = 2$.
Then 
\[
\rep(\overline{Q},\V) = M_2(\C)\times M_2(\C)
\] 
and 
$\Lam_\V \subset \rep(\overline{Q},\V)$ has dimension $4$. 
The variety $\Lam_\V$ has three irreducible components
\[
Z = \{x\in\Lam_\V \mid \rk(x_\alpha) \le 1,\ \rk(x_{\alpha^*}) \le 1 \},\quad
Z' = \{ x \mid  x_{\alpha^*} = 0\},\quad
Z'' = \{x \mid x_\alpha = 0\}. 
\]
We have:
\[
f_Z=\1_Z + \1_{\{(0,0)\}} 
= {\frac{1}{2}}(\1_{Z[1]} *\1_{Z[2]}*\1_{Z[2]}*\1_{Z[1]})
= {\frac{1}{2}}(\1_{Z[2]} *\1_{Z[1]}*\1_{Z[1]}*\1_{Z[2]}), 
\]
\[
f_{Z'} = \1_{Z'} 
= {\frac{1}{4}}(\1_{Z[1]} *\1_{Z[1]}*\1_{Z[2]}*\1_{Z[2]}), \quad
f_{Z''} = \1_{Z''}
= {\frac{1}{4}}(\1_{Z[2]} *\1_{Z[2]}*\1_{Z[1]}*\1_{Z[1]}).
\]
Note that $f_Z\not = \1_Z$ and $\1_Z \not\in\M$.
}
\end{example}

%%%%%%%%%%%%%%%%%%%%%%%%%%%%%%%%%%%%%%%%%%%%%%%%%%%%%
\subsection{} % Flags and composition series
%%%%%%%%%%%%%%%%%%%%%%%%%%%%%%%%%%%%%%%%%%%%%%%%%%%%%
Next, we consider composition series of modules over preprojective
algebras.
Let ${\mathcal X}$ denote the set of pairs $(\ii,\cc)$
where $\ii = (i_1,\ldots ,i_m)$ is a sequence of elements of $I$
and $\cc = (c_1,\ldots ,c_m)\in \{0,1\}^m$.
The integer $m$ is called the length of $(\ii,\cc)$.

Given $(\ii,\cc)\in{\mathcal X}$ such that $\sum_k c_ki_k = |\V|$, we 
define a {\em flag} in $\V$ of type $(\ii,\cc)$ as a sequence
\[
\f = \left(\V=\V^0 \supseteq \V^1 \supseteq \cdots \supseteq \V^m = 0\right)
\]
of graded subspaces of $\V$ such that 
\[
|\V^{k-1}/\V^k| = c_k i_k
\]
for $k=1,\ldots, m$.
Thus $\dim{\V^{k-1}/\V^k}$ is equal to $0$ or $1$.
So these are complete flags, with possible repetition
of some subspaces. 
(It will be convenient below to allow such flags with repeated 
subspaces.)
We denote by $\Phi_{\ii,\cc}$ the variety of flags of type
$(\ii,\cc)$.
When $(c_1,\ldots ,c_m) = (1,\ldots ,1)$
(flags without repetition), we simply write~$\Phi_{\ii}$.

Let $x \in \Lam_\V$. 
A flag $\f$ is $x$-{\it stable} 
if $x(\V^k)\subseteq \V^k$ for all $k$.
We denote by $\Phi_{\ii,\cc,x}$ (resp. $\Phi_{\ii,x}$) 
the variety of $x$-stable flags of type
$(\ii,\cc)$ (resp. of type $\ii$).
Note that an $x$-stable flag is the same as a composition series
of $x$ regarded as a $\Lam$-module.

%%%%%%%%%%%%%%%%%%%%%%%%%%%%%%%%%%%%%%%%%%%%%%%%%%
\subsection{} % The functions $d_\i$
%%%%%%%%%%%%%%%%%%%%%%%%%%%%%%%%%%%%%%%%%%%%%%%%%%
For $(\ii,\cc)\in{\mathcal X}$ with $\sum_k c_ki_k=|\V|$, define
\[
d_{\ii,\cc} = \1_{Z[i_1]}^{c_1} * \cdots * \1_{Z[i_m]}^{c_m} \in \M_\V.
\]
If $c_i=1$ for all $k$, we simply write $d_\ii$ instead of $d_{\ii,\cc}$.
In general, $d_{\ii,\cc}=d_\jj$ where $\jj$ is the subword of $\ii$
consisting of the letters $i_k$ for which $c_k=1$.
By definition, the functions $d_\jj$, where $\jj$ runs over
all words, span $\M$.  
The following lemma results immediately from the definition of the 
product $*$ of constructible functions.

\begin{Lem}\label{lem2}
Let $x\in\Lam_\V$. We have
\[
d_{\ii,\cc}(x) =d_{\jj}(x) =\chi(\Phi_{\jj,x}).
\]
\end{Lem}

\begin{example}
{\rm
Retain the notation of Example~\ref{ex1} (b).
Let $x=(a,b)\in \Lam_\V$ be given by the following matrices
(with respect to some fixed bases of $V_1$ and $V_2$)
\[
a = \left( \begin{matrix} 0 & 0\\ 0 & 0 \end{matrix} \right) \in
\Hom_\C(V_1,V_2),
\qquad
b = \left( \begin{matrix} 1 & 0\\ 0 & 0 \end{matrix} \right) \in
\Hom_\C(V_2,V_1).
\]
Let us calculate $d_{(2,1,2,1)}(x) = \chi(\Phi_{(2,1,2,1),x})$.
To construct a flag 
\[
\f = (\V \supset \V^1 \supset \V^2 \supset \V^3 \supset 0)
\in \Phi_{(2,1,2,1),x}
\]
we first have 
to choose a line $\V^3$ in the 2-dimensional vector space
$V_1 \cap \soc(x) = V_1$.
We may take $\V^3$ to be 

(a) the 1-dimensional image of $x$, or 

(b) any line except this one.

\noindent
In case (a) the module $x^3$ induced by $x$ in the quotient 
$\V/\V^3$ is the semisimple
module.
So we get 
\[
\chi(\Phi_{(2,1,2),x^3}) = 2 \cdot 1 \cdot 1.
\]
In case (b), $x^3=(a^3,b^3)$ where
\[
a^3 = \left( \begin{matrix} 0 \\ 0 \end{matrix} \right) \in 
\Hom_\C(V_1/\V^3,V_2),
\qquad
b^3 = \left( \begin{matrix} 1 & 0 \end{matrix} \right) \in 
\Hom_\C(V_2,V_1/\V^3).
\]  
and at the next stage $\V^2/\V^3$ must be taken as the kernel
of $x^3$ (no choice), and $\V^1$ is also completely determined.
Thus, in case (b) we get
\[
\chi(\Phi_{(2,1,2),x^3}) = 1 \cdot 1 \cdot 1.
\]
So finally, 
\[
\chi(\Phi_{(2,1,2,1),x})= 2 + 1 = 3.
\]
}
\end{example}

%%%%%%%%%%%%%%%%%%%%%%%%%%%%%%%%%%%%%%%%%%%%%%
\subsection{} % PBW bases
%%%%%%%%%%%%%%%%%%%%%%%%%%%%%%%%%%%%%%%%%%%%%%

\def\PBW{\kappa}

In this section we assume that $\g$ is a simple
finite-dimensional Lie algebra.
Equivalently, $Q$ is a Dynkin quiver. 
Then $U(\n)$ has a PBW-basis $B_Q$ associated to this
quiver $Q$.
The image $\Psi(B_Q)$ is easy to describe.
Let $\V\in{\mathcal V}_{\rm fin}(I)$.
The affine space $\rep(Q,\V)$ can be regarded as a subset
of $\Lam_\V$ by identifying it to the set of 
$x\in\Lam_\V$ with $x_\alpha=0$ for every 
$\alpha\in{\overline{Q}}_1 - Q_1$.
Clearly this is an irreducible component
of $\Lam_\V$. Moreover by our assumption
it has finitely many $G_\V$-orbits.

\begin{Lem}\label{PBW}
Let $\orb$ be a $G_\V$-orbit in $\rep(Q,\V)$.
There exists a unique $\PBW_\orb\in \M_\V$
whose restriction to $\rep(Q,\V)$ is the
characteristic function of $\orb$.
The collection of all $\PBW_\orb$ where $\orb$
runs through all $G_\V$-orbits in $\rep(Q,\V)$
is equal to $\Psi(B_Q)\cap\M_\V$.
\end{Lem}
\begin{proof}
By \cite[10.19, 12.12]{Lu91}, the map from $\M_\V$
to $M(\rep(Q,\V))^{G_\V}$ sending $f$ to its restriction
to $\rep(Q,\V)$ is a vector space isomorphism.
Moreover, the space
\[
{\mathcal H}_Q = \bigoplus_\V M(\rep(Q,\V))^{G_\V}
\]
(where $\V$ ranges over all isoclasses of vector spaces in ${\mathcal V}_{\rm fin}(I)$)
endowed with the image of the product $*$
is the geometric realization of the Hall algebra
of $Q$ over $\C$ due to Schofield (see \cite[10.19]{Lu91}). 
In this setting, the PBW-basis is the basis of ${\mathcal H}_Q$
given by the characteristic functions of the $G_\V$-orbits 
in $\rep(Q,\V)$. 
Thus the lemma is proved.
\end{proof}

%%%%%%%%%%%%%%%%%%%%%%%%%%%%%%%%%%%%%%%%%%%%%%%%%%%%%%%%%%%%%%%%%

\section{Comultiplication}\label{comult}

%%%%%%%%%%%%%%%%%%%%%%%%%%%%%%%%%%%%%%%%%%%%%%%%%%%%%%%%%%%%%%%%%

The algebra $U(\n)$ is in fact a bialgebra, the
comultiplication being defined by
\[
n \mapsto n\otimes 1 + 1 \otimes n, \qquad (n\in \n).
\]
In this section we describe the comultiplication $\Delta$ of $\M$
obtained by transporting this comultiplication via $\Psi$.

%%%%%%%%%%%%%%%%%%%%%%%%%%%%%%%%%%%%%%%%%%%%%%%%
\subsection{} % definition of Res
%%%%%%%%%%%%%%%%%%%%%%%%%%%%%%%%%%%%%%%%%%%%%%%%
For $\V, \V', \V'' \in {\mathcal V}_{\rm fin}(I)$ 
such that $|\V'|+|\V''|=|\V|$, define a linear map 
\[
\res^\V_{\V',\V''} : M(\Lam_\V)^{G_\V} \lra 
M(\Lam_{\V'}\times\Lam_{\V''})^{G_{\V'}\times G_{\V''}}
\]
as follows.
We have $\V\stackrel{\sim}{\lra} \V'\oplus\V''$, so given $x'\in\Lam_{\V'}$
and $x''\in\Lam_{\V''}$ we can regard $x'\oplus x''$ as an
element of $\Lam_\V$.
Here, $x'\oplus x''$ denotes the direct sum of $x'$ and $x''$
as endomorphisms of $\V'$ and $\V''$, or equivalently as
modules over the preprojective algebra $\Lam$.
For $f \in M_\V$, $x' \in \Lam_{\V'}$ and $x'' \in \Lam_{\V''}$
set
\[
\left(\res^\V_{\V',\V''}f\right)(x',x'') = f(x'\oplus x'').
\]
This is clearly a constructible function on 
$\Lam_{\V'}\times\Lam_{\V''}$ which is constant on  
$G_{\V'}\times G_{\V''}$-orbits.

%%%%%%%%%%%%%%%%%%%%%%%%%%%%%%%%%%%%%%%%%%%%%%%%%%
\subsection{} %Main lemma
%%%%%%%%%%%%%%%%%%%%%%%%%%%%%%%%%%%%%%%%%%%%%%%%%%
Let $\ii=(i_1,\ldots,i_m)$ with $\sum_k i_k = |\V|$.
Let $\V=\V'\oplus \V''$, $x'\in\Lam_{\V'}$, $x''\in\Lam_{\V''}$ 
and $x=x'\oplus x''\in\Lam_\V$.

\begin{Lem}\label{homo}
We have
\[
\left(\res^\V_{\V',\V''}d_\ii\right)(x',x'')
= \sum_{(\cc',\cc'')} d_{\ii,\cc'}(x')\, d_{\ii,\cc''}(x'')
\]
where the sum is over all pairs $(\cc',\cc'')$ of sequences
in $\{0,1\}^m$ such that
\begin{equation}\label{condition}
\sum_{k=0}^m c'_ki_k=|\V'|,
\qquad \sum_{k=0}^m c''_ki_k=|\V''|,
\qquad c'_k + c''_k = 1\ (0\le k\le m).
\end{equation}  
\end{Lem}

\begin{proof}
By definition we have 
\[
\left(\res^\V_{\V',\V''} d_\ii\right)(x',x'')=d_\ii(x)=\chi(\Phi_{\ii,x}).
\]
To $\f = (\V^l)_{0\le l\le m}\in\Phi_{\ii,x}$, we associate 
the pair of flags
\[
\f'=(\V^l/\V'')_{0\le l\le m},
\qquad
\f''=(\V^l \cap \V'')_{0\le l\le m}.
\]
We regard $\f''$ as a flag in $\V''$, and $\f'$ as
a flag in $\V'$ by identifying $\V'$ with $\V/\V''$.
Clearly, we have $\f'\in\Phi_{\ii,\cc',x'}$ and 
$\f''\in\Phi_{\ii,\cc'',x''}$ for some sequences 
$\cc'$, $\cc''$ in $\{0,1\}^m$ satisfying 
the conditions~(\ref{condition}).
Let $W_1$ denote the set of pairs $(\cc',\cc'')$ satisfying
(\ref{condition}).
For $(\cc',\cc'')\in W_1$, let 
$\Phi_{\ii,x}(\cc',\cc'')$ be the subset of $\Phi_{\ii,x}$
consisting of those $\f$ for which 
$(\f',\f'')\in \Phi_{\ii,\cc',x'}\times\Phi_{\ii,\cc'',x''}$.
Then clearly we have a finite partition
\[
\Phi_{\ii,x} = \bigsqcup_{(\cc',\cc'') \in W_2} \Phi_{\ii,x}(\cc',\cc''),
\]
where $W_2\subseteq W_1$ consists of the pairs $(\cc',\cc'')$ such that
$\Phi_{\ii,x}(\cc',\cc'')$ is nonempty.

Now for $(\cc',\cc'')\in W_2$, the map 
$\alpha(\cc',\cc''): 
\Phi_{\ii,x}(\cc',\cc'') \lra
\Phi_{\ii,\cc',x'}\times\Phi_{\ii,\cc'',x''}$
sending $\f$ to $(\f',\f'')$ is a vector bundle, 
see \cite[Lemma 4.4]{Lu91}.
Hence, 
\[
\chi(\Phi_{\ii,x}(\cc',\cc''))=
\chi( \Phi_{\ii,\cc',x'}\times\Phi_{\ii,\cc'',x''})
=\chi( \Phi_{\ii,\cc',x'})\chi(\Phi_{\ii,\cc'',x''})
\]
and
\[
\left(\res^\V_{\V',\V''} d_\ii\right)(x',x'')=
\chi(\Phi_{\ii,x}) = \sum_{(\cc',\cc'')\in W_2}
\chi( \Phi_{\ii,\cc',x'})\chi(\Phi_{\ii,\cc'',x''}).
\]
On the other hand, 
\[
\sum_{(\cc',\cc'')\in W_1} d_{\ii,\cc'}(x') d_{\ii,\cc''}(x'')
=
\sum_{(\cc',\cc'')\in W_1}\chi(
\Phi_{\ii,\cc',x'})\chi(\Phi_{\ii,\cc'',x''}),
\]
and it only remains to prove that $W_1=W_2$.
Clearly $W_2\subseteq W_1$, so let $(\cc',\cc'')\in W_1$.
Let $\f'\in\Phi_{\ii,\cc',x'}$ and $\f''\in\Phi_{\ii,\cc'',x''}$.
The assumption $x=x'\oplus x''$ implies that $\Phi_{\ii,x}(\cc',\cc'')$
is nonempty. 
Indeed, the flag $\f$ in $\V$ whose $k$th subspace is the 
direct sum of the $k$th subspaces of $\f'$ and $\f''$
is $x$-stable and 
by construction $\f\in\Phi_{\ii,x}(\cc',\cc'')$. 
So $(\cc',\cc'')\in W_2$ and the lemma is proved.
\end{proof}

%%%%%%%%%%%%%%%%%%%%%%%%%%%%%%%%%%%%%%%%%%%%%%%%%%%%%
\subsection{} % The comultiplication
%%%%%%%%%%%%%%%%%%%%%%%%%%%%%%%%%%%%%%%%%%%%%%%%%%%%%

By Lemma~\ref{homo}, the map $\res^\V_{\V',\V''}$ induces 
a linear map from $\M_\V$ to $\M_{\V'}\otimes\M_{\V''}$,
given by
\[
d_\ii\mapsto \sum_{(\cc',\cc'')} d_{\ii,\cc'}\otimes d_{\ii,\cc''},
\]
where the pairs $(\cc',\cc'')$ satisfy (\ref{condition}).
Taking direct sums, we obtain a linear map 
\[
\M_\V \lra \bigoplus_{\V',\V''} \M_{\V'}\otimes \M_{\V''},
\]
the sum being over all isomorphism types $\V'$ and $\V''$
of $I$-graded vector spaces such that $|\V'|+|\V''|=|\V|$.
Taking direct sums over all isomorphism types $\V$, we get
a linear map 
\[
\Delta: \M \lra \M \otimes \M.
\]
Since $(x'\oplus x'')\oplus x''' = x'\oplus(x''\oplus x''')$,
$\Delta$ is coassociative. 
Since $x'\oplus x''\simeq x''\oplus x'$, $\Delta$ is cocommutative. 

Lemma~\ref{homo} shows that $\Delta$ is multiplicative on the elements
$d_\ii$, that is, for $\ii=(i_1,\ldots,i_m)$
\begin{equation}\label{multDelta}
\Delta(d_\ii)=\Delta(d_{i_1})*\cdots *\Delta(d_{i_m})
=
(d_{i_1}\otimes 1 + 1\otimes d_{i_1})
*\cdots *
(d_{i_m}\otimes 1 + 1\otimes d_{i_m}),
\end{equation}
where the product in $\M\otimes\M$
is defined by 
\[
(f_1\otimes f_2)*(g_1\otimes g_2)=(f_1*g_1)\otimes (f_2*g_2).
\]

\begin{Prop}
Under the algebra isomorphism
$\Psi^{-1}:\M\lra U(\n)$, the map $\Delta$ gets identified with the standard
comultiplication of $U(\n)$.
\end{Prop}

\begin{proof}
Equation~(\ref{multDelta}) shows that $\Delta$
is an algebra homomorphism $(\M,*)\ra (\M\otimes\M,*)$.
Moreover the generators $d_{(i)}=\1_{Z[i]}=\Psi(e_i)$ are clearly
primitive. The result follows. \end{proof}

%%%%%%%%%%%%%%%%%%%%%%%%%%%%%%%%%%%%%%%%%%%%%%%%%%%%%%%%%%%%%%%%%%%

\section{Multiplicative properties of the dual semicanonical basis}
\label{multipli}

%%%%%%%%%%%%%%%%%%%%%%%%%%%%%%%%%%%%%%%%%%%%%%%%%%%%%%%%%%%%%%%%%%%

%%%%%%%%%%%%%%%%%%%%%%%%%%%%%%%%%%%%%%%%%%%%%%%%%%%%%%%%%%
\subsection{} % Definition of the dual semicanonical basis
%%%%%%%%%%%%%%%%%%%%%%%%%%%%%%%%%%%%%%%%%%%%%%%%%%%%%%%%%%
The vector space $\M$ is $\ENn^{(I)}$-graded, with 
finite-dimensional homogeneous components. 
Let $\M^*$ denote its graded dual.
Given an $I$-graded space $\V$ and an irreducible component
$Z$ of $\Lam_\V$, we have defined a linear form 
$\rho_Z$ on $M(\Lam_\V)$, see Section~\ref{rhoZ}.
We shall also denote by $\rho_Z$ the element of $\M^*$
obtained by restricting $\rho_Z$ to $\M_\V$ and then
extending by $0$ on all $\M_{\V'}$ with $|\V'|\not = |\V|$.  
Note that by (\ref{caract}), the basis of $\M^*$ dual to the
semicanonical basis $\{f_Z\}$ is nothing but $\{\rho_Z\}$. 
From now on ${\mathcal S}^*=\{\rho_Z\}$ will be called the 
{\it dual semicanonical basis} of $\M^*$.

\begin{Lem}\label{delta}
For $Z \in \irr(\Lam_\V)$
there exists an open dense $G_\V$-stable subset $O_Z \subset Z$ such that
for all $f\in \M_\V$ and all $x \in O_Z$ we have
$\rho_Z(f) = f(x)$.
\end{Lem}

\begin{proof}
For a given $f$, this follows from Section~\ref{rhoZ}. 
Moreover, there exists such an open set simultaneously for all
$f$ because $\M_\V$ is finite-dimensional.
\end{proof}

%%%%%%%%%%%%%%%%%%%%%%%%%%%%%%%%%%%%%%%%%%%%%%%%%%%%%%%%
\subsection{} % delta function
%%%%%%%%%%%%%%%%%%%%%%%%%%%%%%%%%%%%%%%%%%%%%%%%%%%%%%%%

For $x \in \Lam_\V$ define the delta-function $\delta_x\in \M^*$ by 
$
\delta_x(f)=f(x),\ (f\in \M).
$
We then have 
\begin{equation}
\delta_x = \rho_Z, \qquad (x\in O_Z).
\end{equation}
The next Lemma follows immediately.
\begin{Lem}\label{lemmaopen}
Let $Z\in\irr(\Lam_\V)$ and suppose that the orbit of $x\in Z$
is open dense. Then $\rho_Z=\delta_x$.
\end{Lem}

Let $\cdot$ denote the multiplication of $\M^*$ dual to the 
comultiplication $\Delta$ of $\M$.

\begin{Lem}\label{lemmaplus}
Let $x_1\in\Lam_{\V_1}$ and $x_2\in\Lam_{\V_2}$. We have
$\delta_{x_1}\cdot\delta_{x_2} = \delta_{x_1\oplus x_2}.$
\end{Lem}

\begin{proof}
For $f\in \M$, one has
\[
\left(\delta_{x_1}\cdot\delta_{x_2}\right)(f) = 
\left(\delta_{x_1}\otimes\delta_{x_2}\right)(\Delta (f))
=\Delta (f)(x_1,x_2) = f(x_1\oplus x_2)
=\delta_{x_1\oplus x_2}(f).
\]
\end{proof}

\begin{Lem}\label{lem3}
Suppose that $Z = \overline{Z_1\oplus Z_2}$ is an irreducible
component of $\Lam_\V$.
Then there exists $x\in O_Z$ of the form $x=x_1\oplus x_2$ with
$x_1\in O_{Z_1}$ and $x_2\in O_{Z_2}$.
\end{Lem}

\begin{proof} 
The direct sum $Z_1\oplus Z_2$ is the image of the morphism
\[
\theta : G_\V\times Z_1 \times Z_2 \longrightarrow Z 
\]
defined by
\[
\theta(g,x_1,x_2) = g \cdot (x_1\oplus x_2).
\]
Since $O_{Z_1}$ (resp. $O_{Z_2}$) is open dense in $Z_1$
(resp. in $Z_2$), the subset $G_\V\times O_{Z_1} \times O_{Z_2}$
is open dense in $G_\V\times Z_1 \times Z_2$.
Now, since $\theta$ is a dominant morphism between irreducible
varieties, the image under $\theta$ of 
$G_\V\times O_{Z_1} \times O_{Z_2}$ contains a dense open
subset of $Z$, hence it has a nonempty intersection with $O_Z$.
Since both $\theta(G_\V\times O_{Z_1} \times O_{Z_2})$ and $O_Z$ 
are $G_\V$-stable we can find $x$ in their intersection of the form 
$x=x_1\oplus x_2$ with $x_1\in O_{Z_1}$ and $x_2\in O_{Z_2}$. 
\end{proof}

%%%%%%%%%%%%%%%%%%%%%%%%%%%%%%%%%%%%%
\subsection{} % Proof of th1
%%%%%%%%%%%%%%%%%%%%%%%%%%%%%%%%%%%%%

We can now give the proof of Theorem~\ref{ThmA}.

\smallskip\noindent

{\it Proof of Theorem~\ref{ThmA}.}\quad
Choose $x, x_1, x_2$ as in Lemma~\ref{lem3}.
Then Lemma~\ref{delta} and Lemma~\ref{lemmaplus} yield
\[
\rho_{Z_1}\cdot\rho_{Z_2} = \delta_{x_1}\cdot\delta_{x_2}
= \delta_{x_1\oplus x_2} = \delta_x = \rho_Z.
\]
\qed

\begin{Cor}\label{decomp}
Let $Z=\overline{Z_1\oplus\cdots\oplus Z_m}$ be the canonical
decomposition of the irreducible component $Z$ of $\Lam_\V$.
The dual semicanonical basis vector $\rho_Z$ factorizes as 
\[
\rho_Z = \rho_{Z_1} \cdots \rho_{Z_m}.
\]
\end{Cor}

\begin{proof}
For $m=2$ this follows from Theorem~\ref{ThmA}.
Assume that $m>2$. 
By \cite{CBSc} 
\[
Z'=\overline{Z_1\oplus\cdots\oplus Z_{m-1}}
\]
is an irreducible component. 
Moreover 
\[
\overline{Z'\oplus Z_m} = Z,
\]
so by Theorem~\ref{ThmA} we get
$\rho_Z = \rho_{Z'}\cdot\rho_{Z_m}$.
The result follows by induction on $m$.
\end{proof}

The factorization given by Corollary~\ref{decomp} will be called
the {\em canonical factorization} of $\rho_Z$.

%%%%%%%%%%%%%%%%%%%%%%%%%%%%%%%%%%%%%%%%%%%%%%%%%%%%%
\subsection{} % Proof of the only if part of th. 1.2
%%%%%%%%%%%%%%%%%%%%%%%%%%%%%%%%%%%%%%%%%%%%%%%%%%%%%

We shall now deduce from Theorem~\ref{ThmA} the proof of the
``only if'' part of Theorem~\ref{ThmB}.

\begin{Thm}
Let $\g$ be of type $\A_n\ (n\ge 5), \D_n\ (n\ge 4), \E_6, \E_7$ or $\E_8$.
Then, the bases $\B^*$ and ${\mathcal S}^*$ do not coincide.
\end{Thm}

\begin{proof}
Assume first that $\g$ is of type $\A_5$ or $\D_4$.
Then the preprojective algebra $\Lam$ is of tame representation
type.
In this case, we have $\ext_\Lam^1(Z,Z)=0$ for any irreducible component $Z$
of $\Lam_\V$, see \cite{GSc,GS2}.
Therefore by Theorem~\ref{ThmA} and \cite{CBSc} the
square of any vector of ${\mathcal S}^*$ 
belongs to~${\mathcal S}^*$.

On the other hand, it was shown in \cite{L} that for the cases 
$\A_5$ and $\D_4$ there exist elements of $\B^*$ whose square
is not in $\B^*$. 
They are called {\it imaginary vectors} of $\B^*$.
This shows that $\B^*$ and ${\mathcal S}^*$
are different in these cases.

Now if $\g$ is not of type $\A_n$ with $n \leq 4$, then 
the Dynkin diagram of $\g$ contains a subdiagram of type $\A_5$ or $\D_4$,
and the result follows from the cases $\A_5$ and $\D_4$.
\end{proof}

In the next sections we shall prepare some material for the
proof of the ``if'' part of Theorem~\ref{ThmB}, to be given
in Section~\ref{compare}.

%%%%%%%%%%%%%%%%%%%%%%%%%%%%%%%%%%%%%%%%%%%%%%%%%%%%%%%%%%%%%%

\section{Embedding of $\M^*$ into the shuffle algebra}
\label{emb}

%%%%%%%%%%%%%%%%%%%%%%%%%%%%%%%%%%%%%%%%%%%%%%%%%%%%%%%%%%%%%%

We describe a natural embedding of $\M^*$ into the shuffle
algebra.
This is then used to describe a certain family of elements of 
${\mathcal S}^*$ in type $\A_n$.

%%%%%%%%%%%%%%%%%%%%%%%%%%%%%%%%%%%%
\subsection{} % notation
%%%%%%%%%%%%%%%%%%%%%%%%%%%%%%%%%%%%
Let $F=\C\<I\>$ be the free associative algebra over $\C$ generated by $I$.
A monomial in $F$ is called a word.
This is nothing else than a sequence 
$\ii=(i_1,\ldots ,i_k)$ in $I$.
Let 
\[
\pi : F \lra \M
\]
be the surjective algebra homomorphism
given by $\pi(i)=\1_{Z[i]}$, and more generally by
$\pi(\ii)=d_\ii$.
Let $F^*$ denote the graded dual of $F$.
We thus get an embedding of vector spaces 
\[
\pi^* : \M^* \lra F^*.
\]
Let $\{w[\ii]\}$ denote the basis of $F^*$ dual to the
basis $\{\ii\}$ of words in $F$.
Let $\si\in \M^*$.
We have 
\[
\pi^*(\si)=\sum_\ii \pi^*(\si)(\ii)\, w[\ii]
          =\sum_\ii \si(\pi(\ii))\, w[\ii]  
          =\sum_\ii \si(d_\ii)\, w[\ii].  
\]
By Lemma~\ref{lem2}, we obtain in particular
\begin{equation}\label{eqdel}
\pi^*(\delta_x)=\sum_\ii \chi(\Phi_{\ii,x})\, w[\ii].
\end{equation}

%%%%%%%%%%%%%%%%%%%%%%%%%%%%%%%%%%%%%%%%%%%%%
\subsection{} % The shuffle
%%%%%%%%%%%%%%%%%%%%%%%%%%%%%%%%%%%%%%%%%%%%%
Denote by $\shuffle$ the multiplication on $\pi^*(\M^*)$ obtained 
by pushing $\cdot$ with $\pi^*$, that is, for $\si, \tau \in \M^*$
set
\[
\pi^*(\si)\shuffle\pi^*(\tau) = \pi^*(\si\cdot\tau).
\]

\begin{Lem}\label{lemshuf}
The product $\shuffle$
is the restriction to $\pi^*(\M^*)$ of the classical shuffle product 
on $F^*$ defined by
\[
  w[i_1,\ldots ,i_m]\shuffle w[i_{m+1},\ldots ,i_{m+n}]
=
\sum_s w[i_{s(1)},\ldots ,i_{s(m+n)}], 
\]
where the sum runs over the permutations $s\in\SG_{m+n}$ 
such that
\[
s(1)<\cdots < s(m) \quad\mbox{\ and\ }\quad s(m+1)<\cdots < s(m+n).
\]
\end{Lem}

\begin{proof}
This follows easily from Lemma~\ref{homo} and the duality
of $\cdot$ and $\Delta$. Indeed,
\begin{eqnarray*}
\pi^*(\si\cdot\tau)&=& \sum_{\ii} (\si\cdot\tau)(d_\ii)\,w[\ii] \\
                   &=& \sum_{\ii} (\si\otimes\tau)(\Delta(d_\ii))\,w[\ii]\\ 
 &=& \sum_{\ii}
                   \sum_{(\cc',\cc'')}
\left(\si(d_{\ii,\cc'})\otimes\tau(d_{\ii,\cc''})\right)
\,w[\ii],                   
\end{eqnarray*}
where the pairs $(\cc',\cc'')$ are as in Lemma~\ref{homo}.
Now it is clear that the coefficient of $w[\ii]$ in this last 
sum is the same as the coefficient of $w[\ii]$ in the shuffle
product
\[
\left(\sum_{\jj'}\si(d_{\jj'})\,w[\jj']\right)
\shuffle
\left(\sum_{\jj''}\tau(d_{\jj''})\,w[\jj'']\right)
\]
and the lemma is proved.
\end{proof}

%%%%%%%%%%%%%%%%%%%%%%%%%%%%%%%%%%%
\subsection{} % type A
%%%%%%%%%%%%%%%%%%%%%%%%%%%%%%%%%%%
Suppose that $\g$ is of type $\A_n$. 
Then $\M^*\simeq U(\n)^*\simeq\C[N]$,
where $N$ is the group of
unitriangular $(n+1)\times(n+1)$-matrices.

\subsubsection{}
Let us construct an explicit isomorphism
\[
\alpha: \C[N] \ra \M^*.
\] 
Let $t_{ij}$ denote the coordinate function assigning
to $n\in N$ its entry $n_{ij}$.
Then 
\[
\C[N] = \C[t_{ij}\mid 1 \le i<j \le n+1].
\]
It is known that in the isomorphism
$\C[N]\simeq U(\n)^*$, the natural basis 
of $\C[N]$ consisting of monomials in the $t_{ij}$
gets identified to the dual of the PBW-basis of
$U(\n)$ associated to the quiver
\[
Q_n:\qquad
\xymatrix{
1 & 2 \ar[l]_{\alpha_1} & \cdots \ar[l]_{\alpha_2} & n \ar[l]_{\alpha_{n-1}}}
\]
(see for example \cite[3.5]{LNT}).
The $G_\V$-orbits of $\rep(Q_n,\V)$ are naturally labelled by 
the multisegments of degree $|\V|$, and if we denote by
$\{\PBW^*_\msm\}$ the dual in $\M^*$ of the PBW-basis $\{\PBW_\msm\}$
in $\M$, then more precisely
the above isomorphism maps the monomial
$t_{i_1j_1}\cdots t_{i_rj_r}$ to the element $\PBW^*_\msm$
indexed by the multisegment
\[
\msm = [i_1,j_1-1]+\cdots+[i_r,j_r-1].
\]
For $i\le j$, let $x[i,j]$ denote an indecomposable representation of
$Q_n$ with socle $S_i$ and top $S_j$
(up to isomorphism there is exactly one such representation).
Then the orbit of $x[i,j]$ is open dense so 
$\delta_{x[i,j]}$ belongs to $\M^*$ by Lemma~\ref{lemmaopen}.
On the other hand, by Lemma~\ref{PBW}
\[
\delta_{x[i,j]}(\PBW_\msm) = 
\begin{cases}
1 &\text{if $\msm = [i,j]$},\\
0 &\text{otherwise.}
\end{cases}
\]
Hence $\PBW^*_{[i,j]} =\delta_{x[i,j]}$
and $\alpha$ is the algebra homomorphism determined by
$\alpha(t_{i,j+1})=\delta_{x[i,j]}$.

\subsubsection{}
If we regard the functions $t_{ij}$ 
as entries of a unitriangular $(n+1)\times(n+1)$ matrix~$T$
we may consider some special elements of $\C[N]$
given by the minors of this matrix.

\begin{Prop}\label{minor}
The images under $\alpha$ of 
all nonzero minors of the matrix $T$ belong to ${\mathcal S}^*$.
\end{Prop}

\begin{proof}
We shall use the embedding $\pi^*$. 
First note that since $x[i,j]$ has a unique composition series,
\[
\pi^*\alpha(t_{i,j+1}) = w[j,j-1,\ldots , i].
\]
Let $\varphi_{\iu\ju}$ be the $k\times k$-minor taken on the sequence of rows
$\iu=(i_1<\cdots <i_k)$ and the sequence of columns $\ju=(j_1<\cdots < j_k)$.
Since $T$ is a unitriangular matrix with algebraically
independent entries $t_{ij}$ above the diagonal, 
the function $\varphi_{\iu\ju}$ is nonzero if 
and only if $i_r\le j_r$ for every $r$.
We shall assume that this condition is satisfied.
Let 
\[
\lam=(j_k,j_{k-1}+1,\ldots ,j_1+k-1),\quad
\mu=(i_k,i_{k-1}+1,\ldots ,i_1+k-1).
\]
Then $\lam/\mu$ is a skew Young diagram.
We identify it with the following subset of $\ENp\times\ENp$
\[
\lam/\mu=\{(a,b) \mid 1\le b\le k,\ \mu_b < a\le \lam_b\}.
\]
Each pair $(a,b)$ is called a cell of $\lam/\mu$.
Let $y$ be a standard Young tableau of shape $\lam/\mu$,
that is, a total ordering $c_1<\cdots <c_t$
of the cells  of $\lam/\mu$ which is increasing
both on the rows and on the columns.
We associate to $y$ the element
\[
w[y] = w[a_t-b_t,\ldots ,a_1-b_1]
\]
of $F^*$, where $c_r=(a_r,b_r)\ (1\le r\le t)$.
Before we continue with the proof of Proposition~\ref{minor}
we need the following lemma.

\begin{Lem}\label{yshuf}
\[
\pi^*\alpha(\varphi_{\iu\ju}) = \sum_y w[y],
\]
where $y$ runs over the set of standard Young tableaux of
shape $\lam/\mu$.
\end{Lem}

\begin{proof}
Set 
\[
D_{\iu\ju}=\pi^*\alpha(\varphi_{\iu\ju}) \text{ and } S_{\lam\mu}= \sum_y w[y].
\]
We shall prove that $D_{\iu\ju}=S_{\lam\mu}$ by induction on 
the number $t$ of cells of $\lam/\mu$.
If $t=1$, then $D_{\iu\ju} = w[i] = S_{\lam\mu}$ for some $i$, and the
statement is clear. 
So suppose $t>1$. 

For $i=1,\ldots,n$ define $E_i \in \End_\C(F^*)$ by
\[
E_i(w[i_1,\ldots,i_s])= 
\begin{cases}
w[i_1,\ldots,i_{s-1}] & \text{if $i_1=i$},\\
0                     & \text{otherwise}.
\end{cases}
\]
It is immediate to check that $E_i$ is a derivation with
respect to the shuffle product, i.e.
\[
E_i(f\shuffle g) = E_i(f)\shuffle g + f\shuffle E_i(g),
\quad (f,g,\in F^*).
\]
It is also clear that $f=g$ if and only if $E_i(f)=E_i(g)$
for every $i$.
Note that $D_{\iu\ju}$ is the minor on rows $\iu$ and columns $\ju$ 
of the matrix
\[
W = \begin{pmatrix}
1 & w[1] & w[2,1] & \ldots & w[n,n-1,\ldots,1] \\
0 & 1    & w[2]   & \ldots & w[n,n-1,\ldots,2] \\
\vdots &\vdots  &\vdots & & \vdots                       \cr
0 & 0    & 0      & \ldots & 1
\end{pmatrix}
\]
where in the expansion of the determinant the shuffle product
of the entries is used.
It follows that, if $j+1\in \ju$ and $j\not\in \ju$ then 
$E_j(D_{\iu\ju}) = D_{\iu\ku}$ where $\ku$ is obtained from
$\ju$ by replacing $j+1$ by $j$,
and otherwise $E_j(D_{\iu\ju})= 0$.

On the other hand, 
\[
E_j(S_{\lam\mu})= \sum_z w[z],
\] 
where $z$ ranges over the Young tableaux whose 
shape is a skew Young diagram $\nu/\mu$ obtained from  
$\lam/\mu$ by removing an outer cell $c=(a,b)$ with $a-b=j$.
It is easy to check that there is one such
diagram only if $j+1\in \ju$ and $j\not\in \ju$, and that
this diagram then corresponds to the pair $(\iu,\ku)$
above. So, by induction
\[  
E_j(S_{\lam\mu}) = S_{\nu\mu} = D_{\iu\ku} = E_j(D_{\iu\ju})
\]  
in this case, and
$E_j(S_{\lam\mu})=0=E_j(D_{\iu\ju})$, 
otherwise. Therefore $S_{\lam\mu}=D_{\iu\ju}$.
This finishes the proof of Lemma~\ref{yshuf}.
\end{proof}

We continue with the proof of Proposition~\ref{minor}.
Let 
\[
\msm = [i_1,j_1-1] + \cdots + [i_k,j_k-1]
\]
be the multisegment corresponding to the pair $(\iu,\ju)$.
(Here we leave out $[i_l,j_l-1]$ in case $i_l=j_l$.) 
Following \cite{Ri2} 
this parametrizes a laminated $\Lam$-module $x[\msm]$,
that is, a direct sum of indecomposable subquotients of 
projective $\Lam$-modules.
Let $\V$ be the underlying $I$-graded vector space of $x[\msm]$.
It is known that the $G_\V$-orbit of $x[\msm]$ is open dense
in its irreducible component, hence the function
$\delta_{x[\msm]}$ belongs to the dual semicanonical basis.

Now it is easy to see that the types $\ii$ of composition 
series of $x[\msm]$ are
in one-to-one correspondence with the standard Young tableaux
of shape $\lam/\mu$, and that for each tableau, the corresponding
flag variety $\Phi_{\ii,x[\msm]}$ is reduced to a point.
Therefore, comparing with Lemma~\ref{yshuf} we see that 
\[
\pi^*(\delta_{x[\msm]})=\pi^*\alpha(\varphi_{\iu\ju}).
\]
Hence $\alpha(\varphi_{\iu\ju}) = \delta_{x[\msm]}$ belongs to the 
dual semicanonical basis.
This finishes the proof of Proposition~\ref{minor}.
\end{proof}

%%%%%%%%%%%%%%%%%%%%%%%%%%%%%%%%%%%%%%%%%%%%%%%%%%%%%%%%%%%%%%%%%%%%%%%%

\section{A Galois covering of $\Lam$ for type $\A_n$}
\label{Galois}

%%%%%%%%%%%%%%%%%%%%%%%%%%%%%%%%%%%%%%%%%%%%%%%%%%%%%%%%%%%%%%%%%%%%%%%%

In order to prove the ``if'' part of Theorem~\ref{ThmB} we need
to study the canonical decomposition of
$Z\in\irr(\Lam)$ for type $\A_n\ (n\le 4)$.
Our main tool for this will be the Auslander-Reiten quiver of $\Lam$,
which we will calculate by using a Galois covering $\tLam$ of $\Lam$.
This covering is in fact important for all $n$, and it
will also play an essential r\^ole in our
investigation of type $\A_5$ in the last sections of the paper.
So we shall work in type $\A_n$ for general $n$ in the next two
sections, and we shall specify which results are only valid
for $n\le 5$.
We will also exclude the trivial case $\A_1$ and assume
that $n\ge 2$.
 
%%%%%%%%%%%%%%%%%%%%%%%%%%%%%%%%%%%%%%%%%%%%%%
\subsection{} %$CQ_n$ and $\Lam_n$
%%%%%%%%%%%%%%%%%%%%%%%%%%%%%%%%%%%%%%%%%%%%%%
For $n\ge 2$, let again $Q_n$ be the quiver 
\[
\xymatrix{
1 & 2 \ar[l]_{\alpha_1} & \cdots \ar[l]_{\alpha_2} & n \ar[l]_{\alpha_{n-1}}}
\]
of Dynkin type $\A_n$.
Let 
$
\Lam_n = P(Q_n)
$
be the preprojective algebra corresponding to $Q_n$.
Thus
$
\Lam_n = \field \overline{Q_n}/J_n
$
where the double quiver $\overline{Q_n}$ of $Q_n$ is
\[
\xymatrix{
1\ar@<1ex>[r]^{\alpha^*_1}&\ar@<1ex>[l]^{\alpha_1} 2 
\ar@<1ex>[r]^{\alpha^*_2}&\ar@<1ex>[l]^{\alpha_2} %3\ar@<1ex>[r]^{\alpha^*_3}&
%\ar@<1ex>[l]^{\alpha_3} 4 
%\ar@<1ex>[r]^{\alpha^*_4}&
%\ar@<1ex>[l]^{\alpha_4} 
\ar@<1ex>[r]^{\alpha^*_{n-1}} 
\cdots & n \ar@<1ex>[l]^{\alpha_{n-1}},
}
\]
and the ideal $J_n$ is generated by  
\[
\alpha_1\alpha^*_1,\quad \alpha^*_{n-1}\alpha_{n-1},\quad
\alpha^*_i\alpha_i - \alpha_{i+1}\alpha^*_{i+1},\quad (1 \leq i \leq n-2).
\]

%%%%%%%%%%%%%%%%%%%%%%%%%%%%%%
\subsection{} % The covering
\begin{figure}[t]
\[
\def\objectstyle{\scriptstyle}
\def\lablestyle{\scriptstyle}
\xymatrix@-1.0pc{
&&&&\\
\tLam_5:
&{1_3}\ar[rd]^{\alpha^*_{13}}&{}\save[]+<0cm,2ex>*{\vdots}\restore
&\ar[ld]^{\alpha_{22}} 3_2
\ar[rd]^{\alpha^*_{32}}&{}\save[]+<0cm,2ex>*{\vdots}\restore&
\ar[ld]^{\alpha_{41}} 5_1\\
&&\ar[ld]^{\alpha_{12}} 2_2 \ar[rd]^{\alpha^*_{22}}&&
\ar[ld]^{\alpha_{31}} 4_1 \ar[rd]^{\alpha^*_{41}}\\
&{1_2}\ar[rd]^{\alpha^*_{12}}&
&\ar[ld]^{\alpha_{21}} 3_1
\ar[rd]^{\alpha^*_{31}}&&
\ar[ld]^{\alpha_{40}} 5_0\\
&&\ar[ld]^{\alpha_{11}} 2_1 \ar[rd]^{\alpha^*_{21}}&&
\ar[ld]^{\alpha_{30}} 4_0 \ar[rd]^{\alpha^*_{40}}\\
&1_1 \ar[rd]^{\alpha^*_{11}}&&\ar[ld]^{\alpha_{20}} 3_0
\ar[rd]^{\alpha^*_{30}}&&\ar[ld]^{\alpha_{4,-1}} 5_{-1}\\
&&\ar[ld]^{\alpha_{10}} 2_0 \ar[rd]^{\alpha^*_{20}}&&
\ar[ld]^{\alpha_{3,-1}} 4_{-1} \ar[rd]^{\alpha^*_{4,-1}}\\
&1_0 &{}\save[]+<0cm,-2ex>*{\vdots}\restore& 3_{-1} 
&{}\save[]+<0cm,-2ex>*{\vdots}\restore& 5_{-2}\\
&&&\ar[dd]^F&&\\
&&&&&\\
&&&&&\\
\Lam_5: & 1\ar@<1ex>[r]^{\alpha^*_1}&\ar@<1ex>[l]^{\alpha_1} 2 
\ar@<1ex>[r]^{\alpha^*_2}&\ar@<1ex>[l]^{\alpha_2} 3\ar@<1ex>[r]^{\alpha^*_3}&
\ar@<1ex>[l]^{\alpha_3} 4 
\ar@<1ex>[r]^{\alpha^*_4}&
\ar@<1ex>[l]^{\alpha_4} 5\\
}
\]
\caption{\label{Figu1} {\it The Galois covering}}
\end{figure}
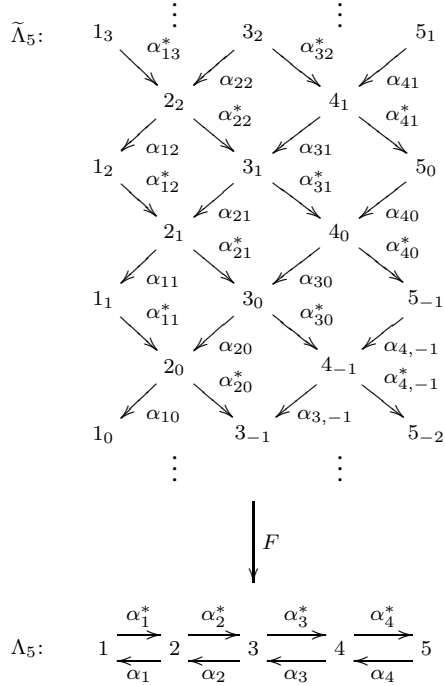

\label{covLam}
%%%%%%%%%%%%%%%%%%%%%%%%%%%%%%
Next, let 
$
\tLam_n = \field \tQ_n/\widetilde{J}_n
$
where $\tQ_n$ is the quiver with vertices
$\{ i_j \mid 1 \leq i \leq n, j \in \Z \}$
and arrows 
\[
\alpha_{ij}: (i+1)_j \to i_j,\quad \alpha^*_{ij}: i_j \to (i+1)_{j-1},\quad
(1 \leq i \leq n-1, j \in \Z),
\]
and the ideal $\widetilde{J}_n$ is generated by
\[
  \alpha_{1j}\alpha^*_{1,j+1},\quad \alpha^*_{n-1,j}\alpha_{n-1,j},\quad 
\alpha^*_{ij}\alpha_{ij} - \alpha_{i+1,j-1}\alpha^*_{i+1,j},\quad
(1 \leq i \leq n-2, j \in \Z).
\]
For $n=5$ we illustrate these definitions in Figure~\ref{Figu1}.

Denote by $\tQa{a,b}_n$ the full (and convex) subquiver of $\tQ_n$
which has as vertices the set
\[
\{i_j\in\tQ_n\mid a+3\leq i+2j\leq b+5\},
\]
and denote by $\Daa{a,b}_n$ the restriction of $\tLam_n$ to $\tQa{a,b}_n$,
see~\ref{reDatla} for an example.

%%%%%%%%%%%%%%%%%%%%%%%%%%%%%%%%%%%%%%%%%%%
\subsection{} % Action of $Z$ and realization of indec $\Lam$-modules
              % as push-downs               
%%%%%%%%%%%%%%%%%%%%%%%%%%%%%%%%%%%%%%%%%%%
The group $\Z$ acts on $\tLam_n$ by $\C$-linear automorphisms via
\[
z \cdot i_j = i_{j+z},\qquad
z \cdot \alpha_{ij} = \alpha_{i,j+z},\qquad
z \cdot \alpha^*_{ij} = \alpha^*_{i,j+z}.
\]
This induces an action 
\[
\Z \times \md(\tLam_n) \longrightarrow \md(\tLam_n)
\]
\[
(z,M) \mapsto {^{(z)}M},
\]
where ${^{(z)}M}$ denotes the $\tLam_n$-module obtained from
$M$ by twisting the action with $z$. Roughly speaking, $^{(1)}M$ is the
same $\tLam_n$-module as $M$, but ``translated one level upwards''.

If we consider $\tLam_n$ and $\Lam_n$ as locally bounded categories 
we have a functor $F: \tLam_n \longrightarrow \Lam_n$
defined by 
\[
i_j \mapsto i,\qquad 
\alpha_{ij} \mapsto \alpha_i,\qquad 
\alpha^*_{ij} \mapsto \alpha^*_i.
\]
This is a Galois covering of $\Lam_n$ in the sense of \cite[\S
3.1]{G},
with Galois group $\Z$.
It provides us with the push-down functor~\cite[\S 3.2]{G} 
\[
\md(\tLam_n) \longrightarrow \md(\Lam_n)
\]
which we also denote by $F$. It is defined as follows.
Let $x \in \md(\tLam_n,\V)$ be a $\tLam_n$-module with 
underlying graded vector space $\V=\bigoplus_{i,j}V_{i_j}$. 
Then $F(x)$ has the same underlying vector space with the
grading $\V = \bigoplus_i V_i$ where
$V_i=\bigoplus_j V_{i_j}$,
and $F(x)$ has maps
$F(x)_{\alpha_i} = \bigoplus_j\, x_{\alpha_{ij}}$ and
$F(x)_{\alpha^*_i}=\bigoplus_j\, x_{\alpha^*_{ij}}$.

The next lemma follows from \cite[\S 3.6]{G} for $n\le 4$, and 
from \cite{DoSk1} for $n=5$, as noted in~\cite[\S 6.3]{DlRi}.

\begin{Lem}\label{push}
Let $2\le n\le 5$.
Any finite-dimensional indecomposable
$\Lam_n$-module is isomorphic to $F(x)$ for some indecomposable
$\tLam_n$-module $x$, which is unique up to a translation
$x \mapsto {^{(z)}x}$ by the Galois group $\Z$.
\end{Lem}

%%%%%%%%%%%%%%%%%%%%%%%%%%%%%%%%%%%%%%%%%%%%%%%
\subsection{} % dimension vectors and pushdcown
%%%%%%%%%%%%%%%%%%%%%%%%%%%%%%%%%%%%%%%%%%%%%%%
For a dimension vector ${\bf d} = (d_{i_j})_{1 \leq i \leq n, j \in \Z}$ 
for $\tLam_n$ define
$d = (d_1, \ldots, d_n)$
where 
\[
d_i = \sum_{j \in \Z} d_{i_j}.
\]
We have a morphism of varieties
\[
\GL(d) \times \md(\tLam_n,{\bf d}) 
\stackrel{F_{{\bf d}}}{\longrightarrow} 
\md(\Lam_n,d) 
%\xrightarrow{\pi_d}  
%\md(\field Q_n,d),
\]
where the {\it push-down morphism} $F_{{\bf d}}$ is defined by
$F_{\bf d}(g,M)= g\cdot F(M)$.

For $z \in \Z$ let ${\mathbf e} = {^{(z)}{\bf d}}$ be the 
$z$th {\it shift} of ${\bf d}$, that is,
\[
e_{i_j} = d_{i_{j-z}},\qquad (1 \leq i \leq n,j \in \Z).
\]
Thus, if ${\bf d}$ is the dimension vector of a $\tLam_n$-module $M$,
then ${\mathbf e}$ is the dimension vector of the shifted module
$^{(z)}M$.
In this case we write ${\mathbf d} \sim {\mathbf e}$.
This defines an equivalence relation~$\sim$.
To simplify our notation, we shall not always distinguish between 
${\bf d}$ and its equivalence class.
However, it will always be clear from the context which one is meant.
To display a dimension vector (or its equivalence class) for $\tLam_n$
we just write down the relevant entries, all other entries are assumed to
be 0.

%%%%%%%%%%%%%%%%%%%%%%%%%%%%%%%%%%%%%%%%%%%%%%%%%%%%%
\subsection{} % AR-quivers
\label{ARQ}
%%%%%%%%%%%%%%%%%%%%%%%%%%%%%%%%%%%%%%%%%%%%%%%%%%%%%
Recall that to any finite-dimensional algebra $A$ over a field
$k$ (and more generally to any locally bounded $k$-algebra
or locally bounded category, see \eg \cite[p. 54]{Ri3}) 
is associated a translation quiver $\Gam_A$ called the 
{\em Auslander-Reiten quiver} of $A$.
It contains a lot of information about the category $\md(A)$. In particular if 
$A$ is representation-finite and standard, one can recover $\md(A)$
from $\Gam_A$.

The vertices of $\Gam_A$ consist of the isomorphism classes of indecomposable
finite-dimensional $A$-modules. 
If $M$ and $N$ are two such modules, then there are
$\dim_k \mathcal{R}(M,N)/\mathcal{R}^2(M,N)$ arrows from 
$M$ to $N$ in $\Gam_A$, 
where $\mathcal{R}$ denotes the
radical of the category $\md(A)$
(compare ~\cite[Chapter 2]{Ri}). 
This means that there is an arrow from
$M$ to $N$  if there exists a nonzero irreducible homomorphism
from $M$ to~$N$. 
Here, by  abuse of notation, we do not
distinguish between a module and its isomorphism class.
The quiver $\Gam_A$ is endowed with an injective map $\tau$,
{\em the translation},
defined on the subset of vertices corresponding to non-projective
modules. 
If $M$ is indecomposable and non-projective then $\tau M=N$
where $0 \ra N \ra E \ra M \ra 0$ is 
the {\em Auslander-Reiten sequence}
(or almost split sequence) ending in $M$.

The {\em stable Auslander-Reiten quiver} $\sGam_A$ of 
$A$ is obtained from $\Gam_A$ by removing all translates 
$\tau^{-n}P$ and $\tau^n I$ $(n \in \N)$ of the projective vertices
$P$ and the injective vertices $I$ 
as well as the arrows involving these vertices. 
Thus the translation $\tau$ induces a permutation on the vertices of 
$\sGam_A$.

We refer the reader to \cite[Chapter VII]{ARS} for more details on
Auslander-Reiten theory, or to~\cite{G0} and~\cite[Chapter 2]{Ri} 
for $k$ algebraically closed as we assume here. 
Note however that these papers  use slightly  different conventions.

\subsection{}  
In Section~\ref{AR} we display the Auslander-Reiten
quivers of $\Lam_2$, $\Lam_3$ and $\Lam_4$.
To calculate them one first determines 
the Auslander-Reiten quivers of their coverings $\tLam_2$,
$\tLam_3$ and $\tLam_4$.
Indeed the algebras $\tLam_n$ are {\em directed}, that is, there 
is no sequence of indecomposable $\tLam_n$-modules of the form
\[
M \ra M_1 \ra \cdots  \ra M_t \ra M
\]
with all homomorphisms being nonzero and non-invertible.
It follows that the Auslander-Reiten quiver of $\tLam_n$ can be calculated
by a combinatorial procedure known as the {\em knitting procedure},
see for example~\cite[\S 6.6]{G0}.
By applying the push-down functor~$F$ to this quiver, one
then obtains the Auslander-Reiten
quiver of $\Lam_n$ \cite[\S 3.6]{G}.

In our pictures, each indecomposable $\Lam_n$-module $M$
is represented by the dimension vector ${\bf d}$ of a $\tLam_n$-module
$M_{\bf d}$ such that $F(M_{\bf d})=M$.
One has to identify each vertex in the extreme left column with 
the vertex in the extreme right column represented by the same
dimension vector up to a shift by the Galois group $\Z$.
The Auslander-Reiten quiver of $\Lam_3$ is shaped like
a Moebius band, and for $\Lam_2$ and $\Lam_4$ one gets a cylinder.

In particular, we see that $\Lam_2$ has 4 isoclasses of indecomposable
modules, $\Lam_3$ has 12, and $\Lam_4$ has 40.

We should point out that for $n>5$ there are $\Lam_n$-modules which are not
obtained from a $\tLam_n$-module via the push-down functor.
%%%%%%%%%%%%%%%%%%%%%%%%%%%%%%%%%%%%%%%%%%%%%%%%%%%
\subsection{} %stable category 
%%%%%%%%%%%%%%%%%%%%%%%%%%%%%%%%%%%%%%%%%%%%%%
\label{Htheo}
For an algebra $A$ let
$\smd(A)$ be the {\em stable category}
of finite-dimensional $A$-modules (\cite[p.~55]{Ri3}).
By definition, the objects of $\smd(A)$ are the same as the
objects of $\md(A)$, and the morphism space $\underline{\Hom}_A(M,N)$
is defined as $\Hom_A(M,N)$ 
modulo the morphisms factoring through projective modules.
The isomorphism classes of indecomposable objects in
$\smd(A)$ correspond naturally to the isomorphism
classes of non-projective indecomposable $A$-modules.
%The Auslander-Reiten quiver of $\smd(A)$ is the
%{\em stable Auslander-Reiten quiver} as defined in \ref{stableAR}.

The stable category $\smd(A)$ is no longer abelian, but
if $A$ is a Frobenius algebra then $\smd(A)$ has the
structure of a triangulated category (\cite[\S I.2]{H}) with
translation functor $\Omega^{-1}$, the inverse of Heller's loop functor.

Moreover, in this situation we may identify the quiver of the 
stable category $\smd(A)$ (see for example~\cite[\S I.4.8]{H})
with the stable Auslander-Reiten quiver
$\sGam_A$ defined in~\ref{ARQ} - just remove the projective-injective
vertices from $\Gam_A$.

%%%%%%%%%%%%%%%%%%%%%%%%%%%%%%%%%%%%%%%%%%%%%%%%%%%%%%%%%%%%
\subsection{} % Happel's theorem
\label{Htheo2}
%%%%%%%%%%%%%%%%%%%%%%%%%%%%%%%%%%%%%%%%%%%%%%%%%%%%%%%%%%%
To a finite-dimensional algebra $B$ one associates
its {\em repetitive algebra} $\widehat{B}$ (\cite[p.~59]{H},
\cite[p.~57]{Ri3}).
This is an infinite-dimensional algebra without unit.
However it is locally bounded, and its
indecomposable projective modules coincide with
its indecomposable injective modules, that is, 
$\widehat{B}$ is a {\em Frobenius algebra}.
 
Hence for a finite-dimensional algebra $B$,
$\smd(\widehat{B})$ is a triangulated category, see~\ref{Htheo}.
Moreover if $B$ has finite global dimension, 
there is an equivalence of triangulated categories
between $\smd(\widehat{B})$
and the derived category $D^b(\md(B))$ of bounded
complexes of $B$-modules (Happel's Theorem \cite[\S II.4]{H}). Under this
equivalence the functor $\Omega^{-1}$ corresponds to the translation functor
$M \mapsto M[1]$.

%%%%%%%%%%%%%%%%%%%%%%%%%%%%%%%%%%%%%%%%%%%%%%%%%%%%%%%%%%%
\subsection{} % Grothendieck groups
\label{Gtgp}
%%%%%%%%%%%%%%%%%%%%%%%%%%%%%%%%%%%%%%%%%%%%%%%%%%%%%%%%%%%

If $B$ is a finite-dimensional algebra of finite global dimension
we have by~\cite[\S 3.2]{HR}
\begin{equation}
\label{eqn:decK}
K_0(\md(\widehat{B}))=K_0(\md(B))\oplus \mathcal{P}(\widehat{B}),
\end{equation}
where $\mathcal{P}(\widehat{B})$ is the subgroup of $K_0(\md(\widehat{B}))$
generated by the classes of projective $\widehat{B}$-modules. 
Thus we may
identify $K_0(\smd(\widehat{B}))$ with $K_0(\md(B))$.
This coincides with the Grothen\-dieck group of $\smd(\widehat{B})$
viewed as a triangulated category.

For a $\widehat{B}$-module $M$ we denote by $[M]$ its class in
$K_0(\md(B))=K_0(\smd(\widehat{B}))$. 
In particular, if $M$ is projective, then $[M] = 0$.
Notice that $[M]$ depends only on the dimension vector
$\dimv{M}\in K_0(\md(\widehat{B}))$ and~\eqref{eqn:decK} provides an
efficient method for calculating $[M]$.
In this context we have
\begin{equation}\label{eqn:Cox}
[\tau M]=\Phi [M]
\end{equation}
where $\Phi$ is the Coxeter transformation of $B$ and
$\tau$ is the self-equivalence of $\smd(\widehat{B})$ induced
by the Auslander-Reiten translation of $\md(\widehat{B})$. This follows
from Happel's Theorem~\ref{Htheo2}, the construction of the
Auslander-Reiten translate in~\cite[\S I.4.6]{H} and the definition of the
Coxeter transformation.

%%%%%%%%%%%%%%%%%%%%%%%%%%%%%%%%%%%%%%%%%%%%%%%%%%%%%%%%%%
\subsection{} % relation between \Da and \tilde\Lam
\label{reDatla}
%%%%%%%%%%%%%%%%%%%%%%%%%%%%%%%%%%%%%%%%%%%%%%%%%%%%%%%%%%

In case $2\leq n\leq 5$  we have
in $\smd(\tLam_n)$ a functorial isomorphism
\begin{equation} \label{eqn:shift}
^{(1)}M\cong \Omega^{-1}\tau^{-1} M 
\end{equation}
which is proved along the lines of~\cite[\S 6.4]{DlRi}. 
In particular this implies 
$[^{(1)}M] = -\Phi^{-1}[M]$.

Moreover, if we set $\Da_n=\Daa{0,n-4}_n$ it is easy to see that
$\tLam_n=\widehat{\Da}_n$. 
For example  $\Da_5$ is the algebra given by the following quiver 
with relations:
\begin{align*}
\Da_5= & 
\vcenter{\def\objectstyle{\scriptstyle}
\xymatrix@-1.2pc{
& \ar[ld] 2_2 \ar[rd] \ar@{--}[dd] & & \ar[ld] 
4_1 \ar[rd] \ar@{--}[dd]\\
1_2 \ar[rd] \ar@{--}@/^/[dd]& & \ar[ld] 3_1 \ar[rd] \ar@{--}[dd]
& & \ar[ld] 5_0 \ar@{--}@/_/[dd]\\
& \ar[ld] 2_1 \ar[rd] & & 4_0 \ar[rd] \ar[dl]\\
1_1 & & 3_0 & & 5_{-1}}}
\end{align*}
Here the dotted lines indicate zero relations and commutativity 
relations.

Thus by Happel's Theorem~\ref{Htheo2} we have
$\smd(\tLam_n)\cong D^b(\md(\Da_n))$. In particular we can apply~\ref{Gtgp}.
Note that $\Da_n$ is hereditary of type $\mathbb{A}_1$ (resp. $\mathbb{A}_3$)
for $n=2$ (resp. $n=3$), and tilted of type $\mathbb{D}_6$ in case
$n=4$. 
Happel's description of the Auslander-Reiten quiver of the
derived category of a (piecewise) hereditary algebra of Dynkin 
type~\cite[\S I.5.6, \S IV.4.3]{H} together with~\eqref{eqn:shift}
yields
the shape of the stable Auslander-Reiten quiver in these cases.
In case $n=5$ the algebra $\Da_5$ is tubular and we can derive the structure
of $\smd(\tLam_5)$ from the known structure of the derived category of a
tubular algebra~\cite{HR}. We will discuss this case in Section~\ref{tubweight}.
%
%%%%%%%%%%%%%%%%%%%%%%%%%%%%%%%%%%%%%%%%%%%%%%%%%%%%%%%%%%%%%%%%%%%%%%%%

\section{From Schur roots to indecomposable multisegments}
\label{indmult}

%%%%%%%%%%%%%%%%%%%%%%%%%%%%%%%%%%%%%%%%%%%%%%%%%%%%%%%%%%%%%%%%%%%%%%%%

%%%%%%%%%%%%%%%%%%%%%%%%%%%%%%%%%%
\subsection{} % Roots
%%%%%%%%%%%%%%%%%%%%%%%%%%%%%%%%%%
\label{IS}
Let $\Idv(\tLam_n)$ be the set of dimension vectors of the
indecomposable $\tLam_n$-modules. 
By $\Idv_S(\tLam_n)$ we denote the set of {\it Schur roots} of
$\tLam_n$, {\em i.e.} the set of dimension vectors of the 
(indecomposable) $\tLam_n$-modules $M$ with
$\End_{\tLam_n}(M) \cong \field$.
For ${\bf d} \in \Idv_S(\tLam_n)$ let 
\[
Z_{{\bf d}} = \overline{\{ M \in \md(\tLam_n,{{\bf d}}) \mid 
\End_{\tLam_n}(M) \cong \field \}}.  
\] 
For $2\le n\le 5$ this is an irreducible component of $\md(\tLam_n,{{\bf d}})$.
We write $Z_{\bf d} \sim Z_{\bf e}$ if ${\bf d} \sim {\bf e}$.

For an irreducible component $Z$ of $\md(\tLam_n,{\bf d})$ define
\[
\eta(Z) = \overline{F_{{\bf d}}(\GL(d) \times Z)}. 
\]
Clearly, $\eta(Z)$ is an irreducible subset of $\md(\Lam_n,d)$, 
but in general $\eta(Z)$ is not an irreducible component.
The following result was shown in \cite{GSc}.

\begin{Thm}\label{ThmGS}
Assume that $2\le n\le 5$. 
For ${\mathbf d}\in \Idv_S(\tLam_n)$, the variety 
$\eta(Z_{\bf d})$ is an indecomposable irreducible
component of $\md(\Lam_n,d)$. 
Moreover, $\eta$ induces a bijection from 
$\Idv_S(\tLam_n)/\!\!\sim$ to $\ind(\irr(\Lam_n))$.
\end{Thm}

%%%%%%%%%%%%%%%%%%%%%%%%%%%%%%%%%%%%%%%
\subsection{} %multisegments and orbits
%%%%%%%%%%%%%%%%%%%%%%%%%%%%%%%%%%%%%%%

Let $\MM(n)$ be the set of multisegments supported on $\{1,\ldots,n\}$,
that is, of the form
\[
\msm = \sum_{1 \leq i \leq j \leq n} m_{ij} [i,j]
\]
where $m_{ij} \in \ENn$.
The {\it degree} $d = (d_1, \ldots, d_n)$ of $\msm$ is given by
\[
d_k = \sum_{i\le k\le j} m_{ij}, \quad (1\le k\le n).
\]
There is a one-to-one correspondence $\msm \mapsto \orb_\msm$ between the set 
$\MM_d(n)$ of multisegments of degree $d$ in $\MM(n)$
and the set of $\GL(d)$-orbits in $\md(\field Q_n,d)$.

Let $\msm_{\rm max}(d)$ be the multisegment labelling
the unique dense orbit of $\md(\field Q_n,d)$.
It can be computed recursively by
\[
\msm_{\rm max}(d) = [a,b] + \msm_{\rm max}(d-e)
\]
where
\[
a = \min \{ i \mid d_i \not= 0 \} \text{ and }
b = \max \{ j \mid d_i \not=0 \text{ for all } a \leq i \leq j \},
\]
and $e_i = 1$ if $i \in [a,b]$ and $e_i = 0$, otherwise.
For example, 
\[
\msm_{\rm max}(1,2,3,1,2) = [1,5]+[2,3]+[3,3]+[5,5].
\]

%%%%%%%%%%%%%%%%%%%%%%%%%%%%%%%%%%%%%%%%%%%%%%%%%%
\subsection{} % irreducible components and orbits
%%%%%%%%%%%%%%%%%%%%%%%%%%%%%%%%%%%%%%%%%%%%%%%%%%
Let $\pi_d$ be the projection morphism from 
$\md(\Lam_n,d)$ to $\md(\field Q_n,d)$ obtained by forgetting
the maps labelled by the arrows $\alpha^*_i$.
Lusztig \cite{Lu91} proved that the irreducible components of $\md(\Lam_n,d)$
are the closures of the sets $\pi_d^{-1}(\orb_\msm)$, $\msm \in \MM_d(n)$.
Define
\[
Z_\msm = \overline{\pi_d^{-1}(\orb_\msm)}.
\]
A multisegment $\msm$ is said to be {\it indecomposable} if $Z_\msm$ is 
indecomposable.
By $\ind(\MM(n))$ we denote the set of indecomposable multisegments
in $\MM(n)$.

%%%%%%%%%%%%%%%%%%%%%%%%%%%%%%%
\subsection{} % Z_max and m_max
%%%%%%%%%%%%%%%%%%%%%%%%%%%%%%%

For a dimension vector 
${\bf d} = (d_{i_j})_{1 \leq i \leq n, j \in \Z}$ for $\tLam_n$ set
\[
\msm_{\rm max}({\bf d}) = \sum_{j \in \Z} \msm_{\rm max}(d(j))
\]
where
$d(j) = (d_{1_j}, d_{2_j},\ldots,d_{n_j})$.
The variety $\md(\tLam_n,d(j))$ is irreducible and, in fact, isomorphic
to $\md(\field Q_n,d(j))$.
There is a $\tLam_n$-module $M_{d(j)}$ such that
the orbit $\orb(M_{d(j)})$ is dense in $\md(\tLam_n,d(j))$,
and $M_{d(j)}$ is uniquely determined up to isomorphism.

Let ${\mathcal F}({\bf d})$ be the set of modules 
$M \in \md(\tLam_n,{\bf d})$ which have submodules $(M_j)_{j \in \Z}$ 
such that
$M_{j-1} \subseteq M_j$ and $M_j/M_{j-1} \cong M_{d(j)}$ for all $j \in \Z$.
Thus $M_j = 0$ if and only if $d(i) = 0$ for all $i \leq j$,
and $M_{j-1} = M_j$ if and only if $d(j) = 0$.
Set
\[
Z_{\rm max}({\bf d}) = \overline{{\mathcal F}({\bf d})}.
\]

For $Z\in\irr(\tLam_n)$ 
let $\mu(Z)$ be the unique multisegment in $\MM(n)$ 
such that 
\[
\overline{\orb_{\mu(Z)}} = \pi_d(\eta(Z)).
\]
This defines a map
\[
\mu: \irr(\tLam_n) \longrightarrow \MM(n).
\]
If $\eta(Z)$ is an irreducible component, say $Z_\msm$, of $\md(\Lam_n,d)$, then
$\mu(Z) = \msm$.

\begin{Lem}
For all dimension vectors ${\bf d}$ for $\tLam_n$ the set
$Z_{\rm max}({\bf d})$ is an irreducible
component of $\md(\tLam_n,{\bf d})$, and we have
\[
\mu(Z_{\rm max}({\bf d})) = \msm_{\rm max}({\bf d}).
\]
\end{Lem}

\begin{proof}
One can easily see that 
\[
\Ext_{\tLam_n}^1(M_{d(i)},M_{d(j)}) = 0 \text{ for all } i \leq j.
\]
Then \cite[Theorem 1.3]{CBSc} yields that $Z_{\rm max}({\bf d})$ is an irreducible
component.
The second part of the lemma follows directly from the definitions.
\end{proof}

%%%%%%%%%%%%%%%%%%%%%%%%%%%%%%%%%%%%
\subsection{} % The Theorem
%%%%%%%%%%%%%%%%%%%%%%%%%%%%%%%%%%%%

Assume that $2\le n\le 5$. Define $\psi : \Idv_S(\tLam_n) \rightarrow \MM(n)$ by
\[
\psi({{\bf d}}) = \begin{cases}
2\,[1,1]+[2,2]+[2,4]+[3,3]+[4,5] & 
\text{if $n=5$ and 
${{\bf d}} = \left( \begin{smallmatrix}
1&&0&&0\\
&1&&1&\\
1&&2&&1\\
&1&&1&
\end{smallmatrix} \right)$},
\vspace{0.5cm}\\
[1,2]+[2,4]+[3,3]+[4,4]+2\,[5,5] & 
\text{if $n=5$ and 
${{\bf d}} = \left( \begin{smallmatrix}
&1&&1&\\
1&&2&&1\\
&1&&1&\\
0&&0&&1
\end{smallmatrix} \right)$},\\
\msm_{\rm max}({{\bf d}}) & \text{otherwise}.
\end{cases}
\] 

\begin{Thm}\label{Thm1}
Let $2\le n \leq 5$. 
The map $\psi$ establishes a bijection
from $\Idv_S(\tLam_n)/\!\!\sim$ to $\ind(\MM(n))$.
Moreover the following diagram commutes:
\[
\xymatrix{
\Idv_S(\tLam_n)/\!\!\sim
\ar[rr]^{\theta: {{\bf d}} \mapsto Z_{{\bf d}}} \ar[dd]^{\psi} && 
\ind(\irr(\tLam_n))/\!\!\sim \ar[dd]^{\eta: Z_{{\bf d}} \mapsto 
\overline{F_{{\bf d}}(\GL(d) \times Z_{{\bf d}})}.}\\
&&\\
\ind(\MM(n)) \ar[rr]^{\phi: \msm \mapsto Z_\msm} & &
\ind(\irr(\Lam_n))
}
\]
\end{Thm}
The proof of Theorem~\ref{Thm1} will be given in Section~\ref{proof1}.

%%%%%%%%%%%%%%%%%%%%%%%%%%%%%%%%%%%%%%%%%%%%%%%%%%%%%%%%%%%%%%%%%%%%%%

\section{Cases $\A_2$, $\A_3$, $\A_4$: the component graph}
\label{comp234}

%%%%%%%%%%%%%%%%%%%%%%%%%%%%%%%%%%%%%%%%%%%%%%%%%%%%%%%%%%%%%%%%%%%%%%

%%%%%%%%%%%%%%%%%%%%%%%%%%%%%%%%%%%%%%%%%%%%%%%%%%%%%%%%%%%%
\subsection{} % results of previous sections for finite type
%%%%%%%%%%%%%%%%%%%%%%%%%%%%%%%%%%%%%%%%%%%%%%%%%%%%%%%%%%%%
In the case that $\Lam_n$ is of finite representation
type, the results of the previous sections simplify greatly.

\begin{Thm}\label{cases234}
Assume that $n =2,3,4$. Then
for each ${\bf d} \in \Idv(\tLam_n)$ there exists (up to isomorphism)
exactly one indecomposable $\tLam_n$-module $M_{{\bf d}}$ with dimension 
vector ${\bf d}$.
Furthermore, $\End_{\tLam_n}(M_{{\bf d}}) \cong \field$ for all ${\bf d}$, 
i.e. 
$
\Idv_S(\tLam_n) = \Idv(\tLam_n).
$
Therefore,
\[
\ind(\irr(\Lam_n)) = \{ \eta(Z_{\bf d}) \mid {\bf d} \in 
\Idv(\tLam_n)/\!\!\sim \}
\]
and
\[
\ind(\MM(n)) = \{ \msm_{\rm max}({\bf d}) \mid {\bf d} \in 
\Idv(\tLam_n)/\!\!\sim \}.
\]
\end{Thm}
\begin{proof}
The first two statements follow from the general theory of directed
(simply connected representation-finite) algebras~\cite[\S 2.4.(8)]{Ri}.
Then we apply Theorem~\ref{ThmGS} and Theorem~\ref{Thm1}.
\end{proof}

%%%%%%%%%%%%%%%%%%%%%%%%%%%%%%%%%%%%%%%%%%%%%%%%%%%%%%%%%%%%%%%
\subsection{} % extensions
%%%%%%%%%%%%%%%%%%%%%%%%%%%%%%%%%%%%%%%%%%%%%%%%%%%%%%%%%%%%%%%

In Section~\ref{A4multi} we list the 40 indecomposable multisegments in
$\MM(4)$ labelling the 40 indecomposable $\Lam_4$-modules,
and we redisplay the Auslander-Reiten quiver of $\Lam_4$ 
with vertices these multisegments.
The translation $\tau$ can be read off by going horizontally two units
to the left, for example 
\[
\tau(\msm_1)=\msm_{16}, \quad \tau(\msm_2) = \msm_{29},
\quad \tau(\msm_3)=\msm_{27},
\]
and so on. 
Note that $\msm_{37}, \msm_{38}, \msm_{39}, \msm_{40}$ have no
$\tau$-translate because they are the projective vertices.

%%%%%%%%%%%%%%%%%%%%%%%%%%%%%%%%%%%%%%%%%%%%%%%%%%%%%%%%%%%
\subsection{} % results
\label{results_extensions}
%%%%%%%%%%%%%%%%%%%%%%%%%%%%%%%%%%%%%%%%%%%%%%%%%%%%%%%%%%%
We shall use $\Gam_{\Lam_n}$ to describe the pairs $(X,Y)$ of indecomposable
$\Lam_n$-modules such that $\Ext_{\Lam_n}^1(X,Y)=0$.
By Proposition~\ref{PropCB}, $\Ext_{\Lam_n}^1(X,Y)=0$ if and 
only if $\Ext_{\Lam_n}^1(Y,X)=0$, so this is a symmetric relation.

Recall, that for any finite-dimensional algebra $\Lam$ the Auslander-Reiten 
formula~\cite[Proposition 4.5]{ARS} gives us 
$D\Ext_\Lam^1(X,Y)\cong\sHom_\Lam(\tau^{-1}Y,X)$.
Now, if $\Lam$ is selfinjective, representation-finite and admits a simply
connected Galois covering $\tLam$ with Galois group $G$, we get
\[
D\Ext_\Lam^1(X,Y)\cong\sHom_\Lam(\tau^{-1}Y,X)\cong
\bigoplus_{g\in G} \sHom_\tLam ( ^{(g)}(\tau\tilde{Y}),\tilde{X}).
\]
Here
$X$ and $Y$ are indecomposable $\Lam$-modules, and $\tilde{X}$ and $\tilde{Y}$
are indecomposable $\tLam$-modules which under push-down give $X$ and $Y$,
respectively. 
In this situation it is easy to determine the dimensions of
the summands in the last term  using additive functions on the stable 
Auslander-Reiten quiver $\sGam_{\tLam}=\Z\Theta$ for some Dynkin quiver
$\Theta$, see~\cite[\S 6.5]{G0}. 

As we have seen, this is exactly the situation for $\Lam=\Lam_n\ (n \le 4)$,
and by~\eqref{eqn:shift} there will be at most one $i\in\Z$ with
$\sHom_{\tLam}(^{(i)}(\tau\tilde{Y}),\tilde{X})\neq 0$.

Since $\tau$ induces a self-equivalence of $\smd(\Lam)$ we have to
do this calculation only for one representative $X$ of each $\tau$-orbit.

The stable Auslander-Reiten quivers of $\Lam_2$, $\Lam_3$ and $\Lam_4$ have
1, 2 and 6 non-trivial $\tau$-orbits, respectively.

For $n=2$ there are only two indecomposable non-projective
$\Lam_n$-modules, say $X$ and $Y$, and
$\Ext_{\Lam_2}^1(X,Y) \not= 0$.

In Sections~\ref{extA3} and \ref{extA4} we display
several copies of the stable Auslander-Reiten quivers of $\Lam_3$ and $\Lam_4$.
In each copy we pointed a representative $X_i$ of a $\tau$-orbit,
and we marked with the sign $\bullet$ all indecomposable
$\Lam_n$-modules $M$ such that
\[
\Ext_{\Lam_n}^1(M,X_i) \not= 0. 
\]
For example, there are 21 indecomposable $\Lam_4$-modules
$M$ such that
\[
\Ext_{\Lam_4}^1(M,X_2) \not= 0.
\]

%%%%%%%%%%%%%%%%%%%%%%%%%%%%%%%%%%%%%%%%%%%%%%%%%%
\subsection{} % components
%%%%%%%%%%%%%%%%%%%%%%%%%%%%%%%%%%%%%%%%%%%%%%%%%%
We note that the previous description shows in particular
that for $n \le 4$ every indecomposable $\Lam_n$-module $X$ satisfies
\[
\Ext_{\Lam_n}^1(X,X) = 0.
\]
This was first observed by Marsh and Reineke.
It follows that the orbit closures of
the indecomposable $\Lam_n$-modules are 
the indecomposable irreducible components of the varieties of
$\Lam_n$-modules.
Therefore the results of \ref{results_extensions} give a complete 
description of the component graph
${\mathcal C}(\Lam_n)$ for $n \leq 4$.

%%%%%%%%%%%%%%%%%%%%%%%%%%%%%%%%%%%%%%%%%%%%%%%%%%%%%%%%%%%%%%%%%%%%%%%%%%

\section{Cases $\A_2$, $\A_3$, $\A_4$: the graph of prime elements of $\B^*$}
\label{prime234}

%%%%%%%%%%%%%%%%%%%%%%%%%%%%%%%%%%%%%%%%%%%%%%%%%%%%%%%%%%%%%%%%%%%%%%%%%%%%%

In this section we consider the dual canonical basis $\B^*$ of $\C[N]$ for
type $\A_n\ (n\le 4)$.
Using the bijection $\msm \mapsto Z_\msm$ we may label its elements
by multisegments $\msm$ or irreducible components $Z$.
We will write 
$\B^*=\{b^*_\msm \mid \msm \in \MM(n)\}$
or
$\B^*=\{b^*_Z \mid Z \in \irr(\Lam_n)\}$
depending on the context.

%%%%%%%%%%%%%%%%%%%%%%%%%%%%%%%%%%%%%%%%%%%%%%%%%%%
\subsection{} % Definition of the graphs G_n
%%%%%%%%%%%%%%%%%%%%%%%%%%%%%%%%%%%%%%%%%%%%%%%%%%%
An element $b^* \in \B^* - \{ 1 \}$ is called 
{\it prime} if it does not have a
non-trivial factorization
$b^*=b^*_1b^*_2$ with $b^*_1,b^*_2 \in \B^*$. 
Let ${\mathcal P_n}$ be the set of prime elements in $\B^*$.

Let ${\mathcal G}_n$ be the {\it graph of prime elements}.
The set of vertices of ${\mathcal G}_n$ is ${\mathcal P}_n$, 
and there is an edge between $b^*_1$ and $b^*_2$ if and only if the 
product $b^*_1b^*_2$ is in $\B^*$.
These graphs give a complete description of the basis $\B^*$
\cite{BZ}.
Indeed, $\B^*$ is the collection of all monomials of the form
\[
\prod_{\msm\in\ind(\MM(n))} (b^*_{\msm})^{k(\msm)}
\]
where the $k(\msm)\in\ENn$ satisfy for $\msm\neq\msm'$
\[
(k(\msm)k(\msm') \neq 0)\quad 
\Longrightarrow \quad ((\msm,\msm') \in {\mathcal G}_n).
\]
Note that for all $b^*$ in $\B^*$, the square of $b^*$ also
belongs to $\B^*$. This is a particular case of the BZ-conjecture
which holds for $\A_n\ (n\le 4)$.
Therefore there is a loop at each vertex of~${\mathcal G}_n$.
Similarly, the $n$
vertices labelled by an irreducible component containing 
an indecomposable projective module are connected to
every other vertex (because the corresponding elements of 
$\B^*_q$ belong to the $q$-center).
Let ${\mathcal G}_n^\circ$ be the graph obtained from
${\mathcal G}_n$ by deleting these $n$
vertices and all the loops.
Clearly ${\mathcal G}_n^\circ$ contains all the information.

%%%%%%%%%%%%%%%%%%%%%%%%%%%%%%%%%%%%%%%%%%%%%%%%%%%%%%%%%%%%%%%%
\subsection{} % Results on the graphs G_n and on the basis \B^*
\label{B^*A4}
%%%%%%%%%%%%%%%%%%%%%%%%%%%%%%%%%%%%%%%%%%%%%%%%%%%%%%%%%%%%%%%%
The graphs ${\mathcal G}_2^\circ$ and ${\mathcal G}_3^\circ$
have been determined by Berenstein
and Zelevinsky \cite{BZ}.
They are respectively dual to an associahedron of type $\A_1$ 
and $\A_3$ in the terminology of \cite{CFZ,FZ2,FZ3}.
In \ref{graphA4} we display the graph 
${\mathcal G}_4^\circ$, which has $36$ vertices (corresponding to
the multisegments $\msm_1, \ldots, \msm_{36}$ in the list of \ref{A4multi})
and 330 edges.
We have calculated it by computer using the BZ-conjecture.
As suggested by Zelevinsky, ${\mathcal G}_4^\circ$ 
is dual to an associahedron of type $\D_6$. 
The maximal complete subgraphs of ${\mathcal G}_n^\circ$ all have 
the same cardinality, namely $1$, $3$ and $6$ for $n = 2$, $3$
and $4$, and are called {\it clusters}.
There are respectively $2$, $14$ and $672$ clusters.
%Each cluster $S$ gives rise to an infinite family of 
%elements of $\B^*$ consisting of all monomials in
%the primes belonging to $S$ (and in the $n$ primes labelled
%by an indecomposable projective $\Lam$-module). 

%%%%%%%%%%%%%%%%%%%%%%%%%%%%%%%%%%%%%%%%%%%%%%%%%%%%%%%%%
\subsection{} % The theorem of Marsh and Reineke 
%%%%%%%%%%%%%%%%%%%%%%%%%%%%%%%%%%%%%%%%%%%%%%%%%%%%%%%%%
The following theorem was proved by Marsh and Reineke 
for $n\le 3$ and conjectured for $n=4$ \cite{MR}.

\begin{Thm}\label{primethm}
For $n\le 4$ the graph ${\mathcal G}_n$ is isomorphic 
to ${\mathcal C}(\Lam_n)$ via the map $b^*_\msm \mapsto Z_\msm$.
\end{Thm}

\begin{proof}
This is checked by using the explicit descriptions of both graphs.
For example, from the first quiver of Section~\ref{extA4} we get
that the vertex labelled by $\msm_{4}$ in  ${\mathcal C}(\Lam_4)$
is connected to all other vertices except 
\[ \msm_3,
\msm_7, \msm_8, \msm_{11}, \msm_{12}, \msm_{13}, \msm_{14}, \msm_{15},
\msm_{31}, \msm_{32}.
\]
The same happens in the graph ${\mathcal G}_4$, as can be seen
from Section~\ref{graphA4}.
\end{proof}

%%%%%%%%%%%%%%%%%%%%%%%%%%%%%%%%%%%%%%%%%%%%%%%%%%%%%%%%%%%%%%%%%%%%%%

\section{End of the proof of Theorem~\ref{ThmB}}

\label{compare}

%%%%%%%%%%%%%%%%%%%%%%%%%%%%%%%%%%%%%%%%%%%%%%%%%%%%%%%%%%%%%%%%%%%%%%

In this section we prove the ``if'' part of Theorem~\ref{ThmB}.

\subsection{}
Let $\g$ be of type $\A_4$.
For brevity set $\rho_\msm = \rho_{Z_\msm}$ for $\msm \in \MM(4)$.

\begin{Prop}\label{egalite}
For every indecomposable multisegment $\msm$ of $\MM(4)$ we have 
\[
b^*_\msm=\rho_\msm.
\]
\end{Prop}

\begin{proof}
For 34 multisegments $\msm_i$ out of the 40 elements of $\ind(\MM(4))$
the vector $b^*_{\msm_i}$ is a minor, and the result follows from 
Proposition~\ref{minor}.

The six elements which are not minors are
$b^*_{\msm_i}$ where $ 31 \le i \le 36$. 
Denoting by $\zeta$ the multisegment duality of
Zelevinsky \cite{Z2}, we have
\[
\zeta(\msm_{31})=\msm_{32},\quad
\zeta(\msm_{33})=\msm_{34},\quad
\zeta(\msm_{35})=\msm_{36}.
\]
Denote also by $\zeta$ the linear involution on $F^*$
given by 
\[
\zeta(w[i_1,\ldots ,i_k])=w[i_k,\ldots,i_1].
\]
It follows from \cite{Lu91,Lu00} that 
\[
\zeta(\pi^*(\B^*)) = \pi^*(\B^*) 
\text{ and }
\zeta(\pi^*({\mathcal S}^*)) = \pi^*({\mathcal S}^*).
\]
Moreover it is known \cite{Z2} that for any multisegment $\msm$
we have
\[
\zeta(\pi^*(b^*_\msm)) = \pi^*(b^*_{\zeta(\msm)}).
\]
Hence it is enough to prove the lemma for $i=31, 33, 35$.
This can be checked by an explicit computation in $F^*$.
The calculation of $\pi^*(b^*_\msm)$ is easy to perform
using the algorithm of \cite{Le}. 
On the other hand, for $x_i$ a point in the dense orbit of
the irreducible component $Z_{\msm_i}$, we
have $\rho_{\msm_i} = \delta_{x_i}$ and 
$\pi^*(\delta_{x_i})$ can be computed via (\ref{eqdel}) in Section~\ref{emb}.
Thus, for $i=31$ we obtain the following expression
for both $\pi^*(b^*_{\msm_i})$ and $\pi^*(\rho_{\msm_i})$:
\begin{eqnarray*}
&&2\,w[4, 2, 3, 3, 1, 2] + 2\,w[2, 4, 3, 1, 3, 2] + 2\,w[2, 4, 3, 3, 1, 2]
      + w[2, 1, 4, 3, 2, 3]
\\&&
 +\ w[2, 1, 3, 4, 3, 2] + 2\, w[2, 1, 4, 3, 3, 2]
      + w[2, 3, 4, 1, 3, 2] + w[4, 2, 1, 3, 2, 3] 
\\&&
+\ w[2, 3, 4, 3, 1, 2] 
+ w[4, 3, 2, 1, 3, 2] + w[4, 3, 2, 3, 1, 2] + 2\,w[4, 2, 1, 3, 3, 2]
\\&&
+\ w[2, 4, 3, 1, 2, 3] + w[2, 4, 1, 3, 2, 3] + 2\,w[2, 4, 1, 3, 3, 2]
+ w[2, 3, 1, 4, 3, 2] 
\\&&
+ \ 2\,w[4, 2, 3, 1, 3, 2] + w[4, 2, 3, 1, 2, 3]    
\end{eqnarray*}
The calculations for $i=33$ and $i=35$ are similar and we omit them.
\end{proof}

\subsection{}

We can now finish the proof of Theorem~\ref{ThmB}.
Assume that $\g$ is of type $\A_4$.
Let
\[
Z = \overline{Z_1 \oplus \cdots \oplus Z_r}
\]
be the canonical decomposition of an irreducible component 
$Z \in \irr(\Lam_4)$.
By Corollary~\ref{decomp}
\[
\rho_Z = \rho_{Z_1} \cdots \rho_{Z_r}.
\] 
All components $Z_k$, $1 \le k \le r$
are of the form 
\[
Z_k = Z_{\msm_{i_k}}
\]
for some indecomposable multisegment $\msm_{i_k}$.
By Proposition~\ref{egalite} we thus have $\rho_{Z_k} = b^*_{Z_k}$.
Moreover, using Theorem~\ref{primethm} we get
\[
b^*_{Z_1} \cdots b^*_{Z_r} \in \B^*.
\]
Hence $\rho_Z$ belongs to $\B^*$.
Thus for Dynkin types $\A_n$ with $n \leq 4$ the dual canonical basis
and the dual semicanonical basis coincide.
\qed

%%%%%%%%%%%%%%%%%%%%%%%%%%%%%%%%%%%%%%%%%%%%%%%%%%%%%%%%%%%%%%%%%%%%%
%   PART A5
%%%%%%%%%%%%%%%%%%%%%%%%%%%%%%%%%%%%%%%%%%%%%%%%%%%%%%%%%%%%%%%%%%%%%

%%%%%%%%%%%%%%%%%%%%%%%%%%%%%%%%%%%%%%%%%%%%%%%%

\section{Case $\A_5$: the tubular algebra $\Da$ and the weighted
projective line $\X$} \label{tubweight}

%%%%%%%%%%%%%%%%%%%%%%%%%%%%%%%%%%%%%%%%%%%%%%%%

For the rest of this article we set 
$\Lam = \Lam_5$, $\tLam = \tLam_5$, $Q=Q_5$, $\Daa{a,b}=\Daa{a,b}_5$, 
$\Da=\Da_5=\Daa{0,1}_5$ and $\MM=\MM(5)$. 
For our convenience we define moreover 
$\Daa{a}=\Daa{a,a}$, $\Da_0=\Daa{0}$, $\Da_\infty=\Daa{1}$
and $\Da^*=\Daa{-1,0}$. 
Note that $\Da^*\cong\Da^{\text{op}}$.

%%%%%%%%%%%%%%%%%%%%%%%%%%%%%%%%%%%%%%%%%%%%%
\subsection{} % tube
\label{tub-fam}
%%%%%%%%%%%%%%%%%%%%%%%%%%%%%%%%%%%%%%%%%%%%%

Almost all components of the Auslander-Reiten quiver
of $\Da$ are {\em tubes}.
This plays a crucial r\^ole in our
results, so we shall recall the definition of a tube
(see \cite[p.~287]{ARS}, \cite[p.~113]{Ri}).

Let $\Z \A_\infty$ be the quiver with vertices
$\{ i_j \mid i \in \Z, j \in \ENp \}$
and arrows  
\[
\{ i_j \to i_{j+1} \mid i \in \Z,\ j \geq 1 \}
\cup 
\{ i_j \to (i+1)_{j-1} \mid i \in \Z,\ j \geq 2 \}.
\]
Define a map $\tau$ on the set of vertices by
$\tau((i+1)_j) = i_j$.
For a vertex $x$ in $\Z \A_\infty$ and $r \geq 1$ let
\[
[x]_r = \{ \tau^{ri}(x) \mid i \in \Z \}.
\]
Then $T_r=\Z \A_\infty/(\tau^r)$ is the quiver with
vertices $[x]_r$, 
and having an arrow $[x]_r \to [y]_r$ if and only if there is an arrow
$x' \to y'$ for some $x' \in [x]_r$ and some $y' \in [y]_r$.
The vertex $[i_j]_r$ in $T_r$ is said to have
{\it quasi-length} $j$.
The map $\tau$ induces a map on $T_r$ again denoted by $\tau$, 
given by 
\[
\tau([x]_r) = [\tau(x)]_r
\]
and called the {\em translation}.
In this way, $T_r$ becomes a translation quiver in the
sense of \cite[p.~47]{Ri}. 
One calls $T_r$ a {\it tube of rank} $r$. 
A tube of rank 1 is called a {\it homogeneous tube}. 
The {\em mouth} of a tube is the subset of vertices of
quasi-length $1$. 
Sometimes we also consider as a tube a translation quiver $\Gam$ 
whose stable part $\underline{\Gam}$, obtained by deleting the translates
of the projective and injective vertices, is a tube.

%%%%%%%%%%%%%%%%%%%%%%%%%%%%%%%%%%%%%%%%%%%%%%%%%%%%%%
\subsection{}  %tubular families
\label{tubfam}
%%%%%%%%%%%%%%%%%%%%%%%%%%%%%%%%%%%%%%%%%%%%%%%%%%%%%%
For an algebra $A$, we call a class of 
indecomposable $A$-modules a tube if the vertices of $\Gam_A$
that belong to that class form a tube~\cite[\S 3.1]{R}. 
A family of tubes $\Tub=(\Tub_x)_{x\in X}$ is called a {\em tubular
$X$-family}.

In our situation the index set $X$ will always be the projective line 
$\mathbb{P}_1(\field)$. 
Such a tubular family is said to be of {\em type} $(m_1,\ldots,m_n)$ if for 
certain points $x_1,\ldots,x_n$ the corresponding tubes have rank 
$m_1,\ldots,m_n$ respectively, and for all remaining points 
the tubes are homogeneous. 
Because of these exceptional points, a better index set is 
provided by the weighted projective line $\X$ 
in the sense of Geigle-Lenzing \cite{LM}, with exceptional
points $x_1,\ldots,x_n$ having respective weights $m_1,\ldots,m_n$.
Below we will point out some strong relations between
the representation theory of $\Lam$ and the 
weighted projective line $\X$ of weight type 
$(m_1,m_2,m_3)=(6,3,2)$ (see \ref{ARtLam}, \ref{+roots}). 

%%%%%%%%%%%%%%%%%%%%%%%%%%%%%%%%%%%%%%%%%%%%%%%%%%%%%%%%%%
\subsection{} % \Da is tubular
\label{datub}
%%%%%%%%%%%%%%%%%%%%%%%%%%%%%%%%%%%%%%%%%%%%%%%%%%%%%%%%%%
The algebra $\Da_\infty$ is a {\em tame concealed algebra} of type
$\widetilde{\D}_6$. This means that $\Da_\infty$ is obtained from
a hereditary algebra of type $\widetilde{\D}_6$ by tilting
with respect to a preprojective tilting module 
\cite[\S~4.3]{Ri}.
The tame concealed algebras have been classified by
Happel and Vossieck \cite{HV}, and one can check
that $\Da_\infty$ belongs to one of the frames in this list.

Similarly $\Da_0$ is a tame concealed algebra of type 
$\widetilde{\E}_7$. 
Its tubular type is $(4, 3, 2)$ (see \cite[p.~158]{Ri}). 
There are two indecomposable $\Da_0$-modules $M$ and $M'$ 
completely determined (up to isomorphism) by their 
respective dimension vectors
\[
\left( \bsm 1&&1&&0\\&1&&1\\0&&1&&1 \esm \right),
\qquad
\left( \bsm 0&&1&&1\\&1&&1\\1&&1&&0 \esm \right).
\]
It can be shown that $M$ and $M'$ lie at the mouth of a tube
of the Auslander-Reiten quiver of $\Da_0$, and that this tube has rank $4$.
Moreover, the tubular extension $\Da_0[M][M']$ 
(see \cite[\S 4.7]{Ri})
is isomorphic to $\Da$. 
It follows that $\Da$ is a {\em tubular algebra} of tubular
type $(6,3,2)$ (see \cite[\S 5]{Ri}).

Similarly, $\Da$ can be regarded as a tubular coextension of 
$\Da_\infty$.

Note that the algebras $\Daa{2i}$ and $\Daa{2i+1}$ are isomorphic to
$\Da_0$ and $\Da_\infty$, respectively, so they are
tame concealed.
Similarly the algebras $\Daa{2i,2i+1}$ and $\Daa{2i-1,2i}$
are isomorphic to $\Da$ and $\Da^*$, respectively, 
hence they are tubular algebras.
Thus we may speak of preprojective, regular and preinjective 
$\Daa{i}$-modules
or $\Daa{i,i+1}$-modules.

%%%%%%%%%%%%%%%%%%%%%%%%%%%%%%%%%%%%%%%%%%%%%%%%%%%%%%%%%%%%%%%%%%%%
\subsection{} % description of \md{\tLam}
\label{destLam}
%%%%%%%%%%%%%%%%%%%%%%%%%%%%%%%%%%%%%%%%%%%%%%%%%%%%%%%%%%%%%%%%%%% 

We are going to define some tubular families of $\tLam$-modules.
Following~\cite[\S 2]{HR} we introduce the following classes of modules.
Let $\cTa{i}$ be the class of indecomposable $\Daa{i-1,i+1}$-modules
$M$ such that the restriction of $M$ to $\Daa{i}$ is regular and nonzero.
Similarly, let $\cMa{i,i+1}$ be the class of indecomposable
$\Daa{i,i+1}$-modules $M$ such that the restriction of $M$ to
$\Daa{i}$ is preinjective and the restriction to $\Daa{i+1}$ is preprojective.

Clearly, we may also interpret the classes $\cTa{i}$ and $\cMa{i,i+1}$ 
as classes of $\tLam$-modules, on which the Galois group $\Z$ acts by
\begin{equation}\label{eq:trans}
^{(j)}\cTa{i}=\cTa{i+2j}, \qquad
^{(j)}\cMa{i,i+1}=\cMa{i+2j,i+2j+1}, \qquad (j\in\Z).
\end{equation}

The classes $\mathcal{T}_0=\cTa{0}$ and $\mathcal{T}_\infty=\cTa{1}$ are
tubular families of type $(6,3,2)$. 
In Section~\ref{Tnull} and~\ref{Tinfty} we display the dimension vectors 
of the
non-homogeneous tubes in these families. The remaining tubes consist only of
the homogeneous $\Daa{i}$-modules.

The classes $\cMa{i,i+1}$ are the regular modules of the tubular algebra
$\Daa{i,i+1}$ which do not belong to the first or last tubular family.
We may decompose $\cMa{i,i+1}$ into a disjoint collection of tubular families 
$(\cTa{i,i+1}_\lam)_{\lam\in\Q^+}$ of type $(6,3,2)$ (see \cite[\S 5.2]{Ri}
and also~\ref{expldec} for more details).

%%%%%%%%%%%%%%%%%%%%%%%%%%%%%%%%%%%%%%%%%%%%%%%%%%%%%%%%%%%%%%%%%%%%%%%%%
\subsection{}   % The class H
%%%%%%%%%%%%%%%%%%%%%%%%%%%%%%%%%%%%%%%%%%%%%%%%%%%%%%%%%%%%%%%%%%%%%%%%%
Define 
\[
\cHe=\cMa{-1,0}\cup\cTa{0}\cup\cMa{0,1}\cup\cTa{1}
\]
and
\[
H=\{\dimv{M}\in K_0(\md(\tLam))\mid M\in\cHe\}\subset\Idv(\tLam).
\]
The next proposition shows that $\cHe$ is a 
transversal of the action of the Galois group on the set of
indecomposable $\tLam$-modules.

\begin{Prop} \label{prp:transv}
For each indecomposable $\tLam$-module $M$ there is a unique $j\in\Z$ 
such that $^{(j)}M\in\cHe$. Moreover for each $\dva{d}\in\Idv(\tLam)$
there is a unique $i\in\Z$ with $^{(i)}\dva{d}\in H$.
\end{Prop}

\begin{proof} By repeating the argument in~\cite[\S 1]{HR} we find that each
indecomposable module is contained in
$\bigcup_{i\in\Z}(\cMa{i-1,i}\cup\cTa{i})$. 
By~\eqref{eq:trans} it only remains to
show the unicity of $j$. This follows from the fact that the classes 
$\cMa{-1,0},\cTa{0},\cMa{0,1},\cTa{1}$ are pairwise disjoint.
Indeed $\cTa{0}$ and $\cTa{1}$ are disjoint as we can see from the dimension
vectors of their objects, see~\ref{Tnull} and~\ref{Tinfty}. 
On the other hand $\cMa{-1,0}$ and $\cMa{0,1}$ are disjoint
because a $\Daa{i}$-module can not be 
preprojective and preinjective at the same time.

The second claim follows from the first.
Indeed, we can determine the class $\cTa{i}$ or $\cMa{i,i+1}$ which
an indecomposable module belongs to by means of $\dimv{M}$
only (see \cite[5.2]{Ri} and \ref{expldec} below).
\end{proof}

%%%%%%%%%%%%%%%%%%%%%%%%%%%%%%%%%%%%%%%%%%%%%%%%%%%%%%%%%%%%%%%%%%%%%%%%%
\subsection{}   % AR-quiver of md(\tLam)
\label{ARtLam}
%%%%%%%%%%%%%%%%%%%%%%%%%%%%%%%%%%%%%%%%%%%%%%%%%%%%%%%%%%%%%%%%%%%%%%%%%
Proposition~\ref{prp:transv} shows in particular the (known) fact that
$\tLam$ is locally support finite. 
Moreover 
the additive closure $\operatorname{add}(\cHe)$ 
({\em i.e.} the full subcategory of $\md(\tLam)$ consisting of modules which
are direct sums of modules of $\cHe$) is closed under
Auslander-Reiten sequences. Thus by~\cite{DoSk} and~\cite[\S 3.6]{G} we may
identify the Auslander-Reiten quivers of $\operatorname{add}(\cHe)$ and
of $\md(\Lam)$.

These considerations together with~\eqref{eqn:shift} imply that
the image $\underline{\cHe}$  of $\cHe$ in $\smd(\tLam)$ is a transversal
for the indecomposable objects in $\smd(\tLam)$ under the action of
$\Z$ by powers of the translation functor $\Omega^{-1}$.

On the other hand, since $\Da$ is a tubular algebra of tubular type $(6,3,2)$,
we know from~\cite{GL} that we have equivalences of triangulated categories
\[
\smd(\tLam)\cong D^b(\md(\Da))\cong D^b(\coh(\X)),
\]
where $\X$ is a weighted projective line of weight type $(6,3,2)$. 
Since $\coh(\X)$ is a hereditary category
we have the useful decomposition
\[
D^b(\coh(\X))=\bigvee_{i\in\Z}\coh(\X)[i].
\]
This is similar to the decomposition
\[
\smd(\tLam)=\bigvee_{i\in\Z}(\operatorname{add}(\underline{\cHe}))[i]
\]
that we have just explained.
Although a lot of our intuition on $\smd(\tLam)$ comes from the
comparison with $\coh(\X)$, we prefer to avoid using this machinery.

%%%%%%%%%%%%%%%%%%%%%%%%%%%%%%%%%%%%%%%%%%%%%%%%%%%%%%%%%%%

\section{Case $\A_5$: the root system}

%%%%%%%%%%%%%%%%%%%%%%%%%%%%%%%%%%%%%%%%%%%%%%%%%%%%%%%%%%
\label{roots}

%%%%%%%%%%%%%%%%%%%%%%%%%%%%%%%%%%%%%%%%%%%%%%%%%%%%%%%%%%
\subsection{} % Ringel form
%%%%%%%%%%%%%%%%%%%%%%%%%%%%%%%%%%%%%%%%%%%%%%%%%%%%%%%%%%

Write $K_0(\Da)$ for the Grothendieck group of $\md(\Da)$.
We have $K_0(\Da)=\Z^{10}$ by identifying the class  of
a $\Da$-module $M$ with its dimension vector $\dimv{M} \in \ENn^{10}$.

Let 
$
\bil{-,-}: \Z^{10} \times \Z^{10} \longrightarrow \Z
$
be the {\it Ringel bilinear form} of $\Da$, given by
\[
\bil{\dimv{M},\dimv{N}}=\sum_{l\ge 0}(-1)^l \dim \Ext^l_\Da(M,N),
\qquad (M, N \in \md(\Da)),
\]
see \cite[p.70, 71]{Ri}.
The algebra $\Da$ being tubular has global dimension 2,
so $\Ext^j_\Da(M,N) = 0$ for $j>2$, and this infinite sum is
in fact finite.
Explicitly for
\[
\cla{d} = \left( \begin{smallmatrix}
&d_{2_2}&&d_{4_1}\\
d_{1_2}&&d_{3_1}&&d_{5_0}\\
&d_{2_1}&&d_{4_0}\\
d_{1_1}&&d_{3_0}&&d_{5_{-1}} 
\end{smallmatrix} \right),
\qquad \cla{e} = \left( \begin{smallmatrix} 
&e_{2_2}&&e_{4_1}\\
e_{1_2}&&e_{3_1}&&e_{5_0}\\
&e_{2_1}&&e_{4_0}\\
e_{1_1}&&e_{3_0}&&e_{5_{-1}} 
\end{smallmatrix} \right)
\]
we have
\[
\bil{\cla{d},\cla{e}} = 
\sum d_{i_j}e_{i_j}  
- \sum_{i_j \to k_l} d_{i_j}e_{k_l} 
+ \sum_{i_j,\,k_l} r_{i_j,k_l} d_{i_j}e_{k_l}  
\]
where the first sum runs over all indices $i_j$,
the second one over all pairs $(i_j,k_l)$ such that there 
is an arrow from $i_j$ to $k_l$, and in the third sum
$r_{i_j,k_l}$ is the number of relations from $i_j$ to $k_l$, i.e.
\[
r_{i_j,k_l} = \begin{cases}
1 & \text{if $(i_j,k_l) \in 
\{ (2_2,2_1),(4_1,4_0),(1_2,1_1),(3_1,3_0),(5_0,5_{-1}) \}$},\\
0 & \text{otherwise}.
\end{cases}
\]

%%%%%%%%%%%%%%%%%%%%%%%%%%%%%%%%%%%%%%%%%%%%%%
\subsection{} % Coxeter transformation
%%%%%%%%%%%%%%%%%%%%%%%%%%%%%%%%%%%%%%%%%%%%%%

Let $E$ be the matrix giving the Ringel form:
\[
\bil{\cla{d},\cla{e}} = \cla{d}^{\rm t} E \cla{e},
\]
where $\cla{d}$ and $\cla{e}$ are interpreted as column vectors in $\Z^{10}$
and $\cla{d}^{\rm t}$ means transposition.
Define the {\it Coxeter matrix} of $\Da$ by 
$\Phi = - E^{-1} E^{\rm t}$ (see \cite[p.71]{Ri}).
It is easy to check that 
\begin{equation}\label{Cox}
\bil{\cla{d},\cla{e}} = - \bil{\cla{e},\Phi(\cla{d})} = 
\bil{\Phi(\cla{d}),\Phi(\cla{e})}, \qquad (\cla{d},\cla{e}\in \Z^{10}).
\end{equation}
The data ${\mathcal L}=(\Z^{10},\bil{-,-},\Phi)$ is called a 
{\em bilinear lattice}~\cite{Len}.

%%%%%%%%%%%%%%%%%%%%%%%%%%%%%%%%%%%%%%%%%%%%%%%%%%%%%%%%%%%%%%%%%%%%%%%%%%%%
\subsection{} % Isometries
\label{sus:isomt}
%%%%%%%%%%%%%%%%%%%%%%%%%%%%%%%%%%%%%%%%%%%%%%%%%%%%%%%%%%%%%%%%%%%%%%%%%%%
Let $\cla{d}=(d_{i_j})_{1\leq i\leq 5,\,j\in\Z}$ be an element of
$K_0(\tLam)=K_0(\md(\tLam))$. The {\em support} $\supp(M)$ of 
$M$ is defined as the set of vertices $i_j$ of $\tQ$  
such that the $i_j$th component of $\dimv{M}$ is nonzero.
For $i \le j$ 
it will be convenient to identify the Grothendieck group 
$K_0(\Daa{i,j})$ with the subgroup of $K_0(\tLam)$ of elements with support
in $\tQa{i,j}$, see~\ref{covLam}.
We shall denote by $\bia{i,j}{-,-}$ the Ringel form on
$K_0(\Daa{i,j})$, so that $\bil{-,-}=\bia{0,1}{-,-}$.

Let $\dva{p}_1,\dva{p}_3,\dva{p}_5$ (resp. $\dva{p}_2,\dva{p}_4$)
be the dimension vectors of the indecomposable projective 
$\tLam$-modules with support in $\tQa{0,2}$ (resp. in
$\tQa{-1,1}$),
that is, 
\[
\dvp{1} = \left( \bsm 1&&0&&0\\&1&&0\\0&&1&&0\\&0&&1\\0&&0&&1 \esm \right),\quad
\dvp{3} = \left( \bsm 0&&1&&0\\&1&&1\\1&&1&&1\\&1&&1\\0&&1&&0 \esm \right),\quad 
\dvp{5} = \left( \bsm 0&&0&&1\\&0&&1\\0&&1&&0\\&1&&0\\1&&0&&0 \esm \right),
\]
\[
\dvp{2} = \left( \bsm &1&&0\\1&&1&&0\\&1&&1\\0&&1&&1\\&0&&1 \esm \right),\qquad
\dvp{4} = \left( \bsm &0&&1\\0&&1&&1\\&1&&1\\1&&1&&0\\&1&&0 \esm \right).
\]
Up to shift these are the dimension vectors of all indecomposable projective 
$\tLam$-modules. Note that these are precisely the dimension vectors of
the projective $\tLam$-modules that belong to $\cHe$.
The following result is verified by a direct calculation.

\begin{Lem} \label{lem:isom}
The map $K_0(\tLam)\lra K_0(\smd(\tLam))=K_0(\Da),
\dimv{M}\mapsto [M]$
(see~\ref{Gtgp}) induces isometries
\begin{align*}
\iota_+\colon K_0(\Daa{-1,0})\lra K_0(\Da), \quad &
\dva{x}\mapsto \dva{x}-x^\phg_{4_{-1}}\dvp{2}-x^\phg_{2_0}\dvp{4} \\
\iota_-\colon K_0(\Daa{1,2})\lra K_0(\Da),\quad &
\dva{x}\mapsto 
\dva{x}-x^\phg_{1_3}\dvp{1}-x^\phg_{3_2}\dvp{3}-x^\phg_{5_1}\dvp{5}
\end{align*}
of bilinear lattices, that is,
\[
\bia{-1,0}{\cla{x},\cla{y}}=\bil{\iota_+(\cla{x}),\iota_+(\cla{y})},\qquad 
\bia{1,2}{\cla{x},\cla{y}}=\bil{\iota_-(\cla{x}),\iota_-(\cla{y})}.
\]
\end{Lem}

%%%%%%%%%%%%%%%%%%%%%%%%%%%%%%%%%%%%%%%%%%%%%%%%
\subsection{} % The quadratic form and the radical
%%%%%%%%%%%%%%%%%%%%%%%%%%%%%%%%%%%%%%%%%%%%%%%%

Let $\qda{i,j}: K_0(\Daa{i,j}) \lra \Z$
be the quadratic form 
$\qda{i,j}(\dva{d}) = \bia{i,j}{\dva{d},\dva{d}}$. 
For simplicity we write $\qda{i}=\qda{i,i}$ and $q=\qda{0,1}$.
Define
$\bha{i}$ to be the positive generator of the one-dimensional radical
of $\qda{i}$ (recall that $\Daa{i}$ is tame concealed). 
This means for example that the support of $\bha{0}$ and of $\bha{1}$ 
is contained in $\tQa{0,1}$. 
Explicitly
\[ 
\bha{0} = \left( 
\begin{smallmatrix} &0&&0\\1&&2&&1\\&3&&3\\1&&2&&1 \end{smallmatrix}
\right),
\qquad
\bha{1} = \left( \begin{smallmatrix} &1&&1\\1&&2&&1\\&1&&1\\0&&0&&0 
\end{smallmatrix} \right).
\]
Moreover  $\bha{2i}={^{(i)}\bha{0}}$ and $\bha{2i+1}={^{(i)}\bha{1}}$.
Since $\qda{i,i+1}$ is positive semidefinite of corank 2 it is easy to see that
$\rad(\qda{i,i+1}) = \Z\,\bha{i}\oplus\Z\,\bha{i+1}$.
Notice that (Lemma~\ref{lem:isom})
\begin{equation}
\label{eq:lah}
\iota_+(\bha{-1})=-\bha{1}, \qquad \iota_+({\bha{0}}) = \bha{0}.
\end{equation}
It will be convenient to set $\dva{h}_0=\bha{0}$
and $\dva{h}_\infty=\bha{1}$.
These are the two generators of $\rad(q)$.

%%%%%%%%%%%%%%%%%%%%%%%%%%%%%%%%%%%%%%%%%%%%%%%%%%%%%
\subsection{} %  root system
%%%%%%%%%%%%%%%%%%%%%%%%%%%%%%%%%%%%%%%%%%%%%%%%%%%%%
Let
\[
R = \{ \cla{d} \in \Z^{10} \mid q(\dva{d}) = 0,1,\ \cla{d} \not = 0\}
\]
be the set of {\it roots} of $q$.
A root $\dva{d} \in R$ is called {\it imaginary} if $q(\dva{d}) = 0$, and 
{\it real} otherwise.
Let $R^{\rm im}$ (resp. $R^{\rm re}$) be the set of imaginary (resp. real) 
roots of $q$.

The form $q$ being positive semidefinite, $R^{\rm im}$ consists of
the nonzero elements of $\rad(q)$.
Note that, by (\ref{Cox}), $\dva{d}$ is a radical vector if and only 
if $\Phi(\dva{d})=\dva{d}$.
Indeed, 
\[
\bil{\dva{d},\dva{e}}+\bil{\dva{e},\dva{d}}=\bil{\dva{d}-\Phi(\dva{d}),\dva{e}}
\]
and the Ringel form $\bil{-,-}$ is nondegenerate.

Since $\Da$ has tubular type $(6,3,2)$,
the form~$q$ induces a positive definite form of type $\E_8$
on $\Z^{10}/\rad(q)$.
(Note that the Dynkin diagram of type $\widetilde{\E}_8$ is a star with
three branches of lengths $6, 3, 2$.)
Thus, $R$ is an elliptic root system of type $\E_8^{(1,1)}$
in the sense of Saito \cite{Sa}.

%%%%%%%%%%%%%%%%%%%%%%%%%%%%%%%%%%%%%%%%%%%%%%%%%%%%%%%%%%%%%%%%%%%%%%%%%%%%%
\subsection{}  % Positive roots
%%%%%%%%%%%%%%%%%%%%%%%%%%%%%%%%%%%%%%%%%%%%%%%%%%%%%%%%%%%%%%%%%%%%%%%%%%%%
\label{+roots}
By \cite[p.278]{Ri}, the elements of $R \cap \ENn^{10}$
are precisely the dimension vectors of the indecomposable
$\Da$-modules.
Note however that there exist elements in $R$ whose
coordinates are not all nonnegative or all nonpositive.
For example
\[
\left( 
\begin{smallmatrix} &-1&&0\\0&&0&&0\\&1&&0\\0&&0&&0 \end{smallmatrix}
\right)
\in R.
\]
Following \cite{LM} and \cite{GSc} we
will introduce another notion of positive root in~$R$.
Define the subset of {\it positive roots} of $q$ as
\[
R^+ = \{ \dva{r} \in R \mid \bil{\dva{r},\dva{h}_\infty} > 0 \text{ or } 
(\bil{\dva{r},\dva{h}_\infty} = 0 \text{ and } \bil{\dva{h}_0,\dva{r}} > 0) \},
\]
and set $R^- = \{-\dva{d} \mid \dva{d}\in R^+\}$.

In our setting $R^+$ is the set of classes in $K_0(\Da)=K_0(\smd(\tLam))$
of the non-projective (indecomposable) $\tLam$-modules which belong to $\cHe$
as we will see in Corollary~\ref{cor:eta}. 
Then it follows from~\ref{ARtLam} that
$R = R^+ \cup R^-$, see also~\cite{LM}. One could also derive this
fact from a careful analysis of $R$.

Our definition  of $R^+$ is partially inspired by the
interpretation of the bilinear lattice ${\mathcal L}$ as the 
Grothendieck group of the category ${\rm coh}(\X)$, see~\ref{ARtLam}.
Indeed, since ${\rm coh}(\X)$ and $\md(\Da)$ are derived equivalent~\cite{LM},
they have isomorphic Grothendieck groups \cite[p.101]{H}.
The linear forms $\bil{-,\dva{h}_\infty}$ and 
$\bil{\dva{h}_0,-}$ correspond to the
functions ``rank'' and ``degree'' on $K_0(\coh(\X))$,
and $R^+$ corresponds to the set of classes of indecomposable
sheaves.

%%%%%%%%%%%%%%%%%%%%%%%%%%%%%%%%%%%%%%%%%%%%%%%%%%%%%%%
\subsection{} % slope
%%%%%%%%%%%%%%%%%%%%%%%%%%%%%%%%%%%%%%%%%%%%%%%%%%%%%%%
Since for $\dva{d}\in R^+$ we have $\bil{\dva{d},\dva{h}_\infty} \neq 0$
or $\bil{\dva{h}_0,\dva{d}} \neq 0$, the quotient 
$\bil{\dva{h}_0,\dva{d}}/\bil{\dva{d},\dva{h}_\infty}$
is a well-defined element of $\Q_\infty$.
For $\lam \in \Q_\infty$ let
$
R^\lam = \{ \dva{d} \in R^+ \mid 
\bil{\dva{h}_0,\dva{d}}/\bil{\dva{d},\dva{h}_\infty} = \lam \}
$
be the set of {\it positive roots of slope} $\lam$.
Thus
\[
R^+ = \bigcup_{\lam \in \Q_\infty} R^\lam.
\]
For $a \in \ENn$ and $b \in \Z$ such that $a > 0$ or ($a=0$ and $b>0$) 
set
\[
\dva{h}_{a,b} = a\dva{h}_0 + b\dva{h}_\infty.
\]
We have
$
R^\lam \cap R^{\rm im} = \{ \dva{h}_{a,b} \mid b/a = \lam \}.
$
This follows immediately from the equalities
\[
\bil{\dva{h}_0,\dva{h}_\infty}=-\bil{\dva{h}_\infty,\dva{h}_0}=6,\quad
\bil{\dva{h}_0,\dva{h}_0}=\bil{\dva{h}_\infty,\dva{h}_\infty}=0.
\]

Set
\[
R^{++} = \bigcup_{\lam \in \Q^+} R^\lam, \qquad 
R^{+-} = \bigcup_{\lam \in \Q^-} R^\lam.
\]
Thus 
$
R^+ = R^{+-} \cup R^0 \cup R^{++} \cup R^\infty.
$
%%%%%%%%%%%%%%%%%%%%%%%%%%%%%%%%%%%%%%%%%%%%%%%%%%%%%%%%%%%%%%%%%%%
\subsection{} % explain decomposition
\label{expldec}
%%%%%%%%%%%%%%%%%%%%%%%%%%%%%%%%%%%%%%%%%%%%%%%%%%%%%%%%%%%%%%%%%%%
The above decomposition of $R^+$ is motivated by the following fact.
For $\lam\in\Q^+$ let $\cTa{i,i+1}_\lam$ be the class of indecomposable
$\Daa{i,i+1}$-modules $M$ with 
\[
\frac{\bia{i,i+1}{\bha{i},\dimv{M}}}{\bia{i,i+1}{\dimv{M},\bha{i+1}}}=\lam.
\]
Then $\cTa{i,i+1}_\lam$ is a tubular family of type $(6,3,2)$ 
and $\bigcup_{\lam\in\Q^+}\cTa{i,i+1}=\cMa{i,i+1}$.
Moreover
the set of dimension vectors $\{\dimv{M}\mid M\in\cTa{i,i+1}_\lam\}$ coincides
with the set of roots $\dva{v}$ of $K_0(\Daa{i,i+1})$ such that
\[
\frac{\bia{i,i+1}{\bha{i},\dva{v}}}{\bia{i,i+1}{\dva{v},\bha{i+1}}}=\lam,
\]
see~\cite[5.2]{Ri}.
For $\lam\in\Q_\infty$ we may now define the tubular family
\[
\Tub_\lam =\begin{cases}
\cTa{1}              &\text{ for } \lam=\infty,\\
\cTa{0,1}_\lam       &\text{ for } 0<\lam<\infty,\\
\cTa{0}              &\text{ for } \lam=0,\\
\cTa{-1,0}_{-1/\lam} &\text{ for } \lam < 0.
\end{cases}
\]
Clearly we have $\bigcup_{\lam\in\Q_\infty} \Tub_\lam =\cHe$.
We shall denote by  $T_{\lam,x}$ the tube of the family $\Tub_\lam$
indexed by $x\in\X$.
Thus, for all points $x$ except the three exceptional points 
$x_1, x_2, x_3$, $T_{\lam,x}$ is a homogeneous tube,
and $T_{\lam,x_1}, T_{\lam,x_2},T_{\lam,x_3}$ have respective ranks 
$6, 3$ and $2$.

\begin{Lem} \label{lem:tubes}
For $\lam\in\Q_\infty$ we have 
$R^{\lam}=\{[M]\mid M\in\Tub_\lam\text{\ and $M$ is non-projective}\}$.
\end{Lem}

\begin{proof}
For $\lam\in\{0,\infty\}$ this may be verified directly by comparing the
dimension vectors of the objects in $\Tub_0$ (resp. $\Tub_\infty$) which we 
describe in~\ref{Tnull} (resp.~\ref{Tinfty}) with
the classes in $R^0$ (resp. $R^\infty$) which we can construct explicitly,
see Theorem~\ref{Thm3}.
For $0<\lam<\infty$ the result is clear from the case $i=0$ of
Ringel's result which we explained above. For $\lam<0$
it follows from the same result for $i=-1$ after applying the 
isometry $\iota_+$ from
Lemma~\ref{lem:isom} and~\eqref{eq:lah}.
\end{proof}

\begin{Cor} \label{cor:eta}
The restriction $\xi\colon H\lra R^+\cup\{0\}$ of the map
$K_0(\tLam)\lra K_0(\Da),\dimv{M}\mapsto [M]$, is well-defined and surjective.
Moreover, we have
\begin{align*}
\xi^{-1}(0) &=\{\dvp{1},\dvp{2},\dvp{3},\dvp{4},\dvp{5}\},\\
\xi^{-1}(m\dva{h}_0)
&=\{m\dva{h}_0,m\dva{h}_0+\dvp{2},m\dva{h}_0+\dvp{4}\},
\quad (m\in\ENp),\\
\xi^{-1}(m\dva{h}_\infty) &=
\{m\dva{h}_\infty,m\dva{h}_\infty+\dvp{1},m\dva{h}_\infty+\dvp{3},
m\dva{h}_\infty+\dvp{5}\},\quad (m\in\ENp).
\end{align*}
If $\dva{r}\in R^+ \setminus \{\ENp\dva{h}_0\cup\ENp\dva{h}_\infty\}$ the set
$\xi^{-1}(\dva{r})$ contains precisely one element.
\end{Cor}
We leave the proof to the reader.

%%%%%%%%%%%%%%%%%%%%%%%%%%%%%%%%%%%%%%%%%%%%%%%%%%
\subsection{} % rank and quasi-length
%%%%%%%%%%%%%%%%%%%%%%%%%%%%%%%%%%%%%%%%%%%%%%%%%%
Since $\Da$ is tubular of tubular type $(6,3,2)$,
we have that $\Phi^6=1$ (see \cite{Len,Len2}). 
Of course this can also be verified 
directly by a simple calculation.

Define the {\it rank} $\rk(\cla{d})$ of a vector $\cla{d}\in\Z^{10}$ 
as the minimal $n \geq 1$ such that $\Phi^n(\cla{d}) = \cla{d}$. 
Hence, if $\cla{d}$ is an imaginary root $\rk(\cla{d})=1$, and
if $\cla{d}$ is a real root $\rk(\cla{d}) = 2, 3$ or $6$.

The {\it quasi-length} $\ql(\cla{d})$ of $\cla{d}$ is defined as the greatest
common divisor of the entries in the imaginary root
\[
h(\cla{d}) = \sum_{i=1}^{\rk(\cla{d})} \Phi^i(\cla{d}).
\]
In other words, writing $h(\cla{d}) = \dva{h}_{a,b}$ we have 
$\ql(\cla{d}) = \gcd(a,b)$.
Set also
\[
\iso(\cla{d}) = \frac{1}{\ql(\cla{d})} h(\cla{d}).
\]
Then $\cla{d} \in R^\lam$ if and only if $\iso(\cla{d})=\cla{h}_{a,b}$
where $\lam=b/a$ and $\gcd(a,b)=1$.

These definitions are motivated by the fact that for a non-projective 
(indecomposable)
$\tLam$-module $M$ which has quasi-length $\ell$ in a tube of rank $r$ we have
\begin{equation}\label{eqn:ql}
\ql([M])=\begin{cases}\ell    &\text{ if $r \nmid \ell$},\\
                      \ell/r  &\text{ if $r\mid  \ell$},
\end{cases}\qquad
\rk([M])=\begin{cases} r   &\text{ if $r\nmid \ell$},\\
                       1   &\text{ if $r\mid  \ell$}.
\end{cases}
\end{equation}
(Of course our $\rk$ should not be confused with the rank function on
$K_0(\coh(\X))$).

%%%%%%%%%%%%%%%%%%%%%%%%%%%%%%%%%%%%%%%%%%%%%%%%%%%%%
\subsection{} \label{R-lambda-l-i}% The subsets R^\lam_l(i)
%%%%%%%%%%%%%%%%%%%%%%%%%%%%%%%%%%%%%%%%%%%%%%%%%%%%%
For $\lam \in \Q_\infty$, 
$i \in \{ 1,2,3,6 \}$ and $\ell \in \ENp$ let 
\[
R^\lam_\ell(i) = \{ \cla{d} \in R^\lam \mid \rk(\cla{d}) = i,
\ \ql(\cla{d}) = \ell \}
\]
so that
\[
R^\lam = \bigcup_{\ell \in \ENp,\ i \in  \{ 1,2,3,6 \}} R^\lam_\ell(i).
\]

\begin{Lem} \label{lemRli}
\begin{itemize}
\item[(a)]
If $\cla{d}\in R^\lam_\ell(i)$ then 
$R^\lam_\ell(i) = \{ \Phi^j(\cla{d}) \mid j=1,\ldots, i\}$;
\item[(b)]
If $i > 1$ and $i$ divides $\ell$, then $R^\lam_\ell(i) = \emptyset$;
\item[(c)]
If $i > 1$ and $i$ does not divide $\ell$, then $R^\lam_\ell(i)$ consists
of $i$ real roots;
\item[(d)]
If $i = 1$, then
$R_\ell^\lam(1) = \{ h_{a,b} \}$
where $b/a = \lam$ and $\gcd(a,b) = \ell$.
\end{itemize}
\end{Lem}

\begin{proof}
Let $\dva{d}\in R^\lam_\ell(i)$.
Clearly, $\rk(\Phi(\cla{d})) = \rk(\cla{d})$, $h(\Phi(\cla{d}))=h(\cla{d})$,
and $\ql(\Phi(\cla{d})) = \ql(\cla{d})$.
Moreover, by Equation~(\ref{Cox}) we have $\Phi(\cla{d})\in R^\lam$.
Hence the $\Phi$-orbit of $\cla{d}$ is contained in $R^\lam_\ell(i)$.

Conversely, if $R^\lam_\ell(i)$ is nonempty, then it
consists of a single $\Phi$-orbit.
This  follows from~\eqref{eqn:Cox}
and Lemma~\ref{lem:tubes} (remember that $\Tub_\lam$ is a tubular family of
type $(6,3,2)$).
Note that all the indecomposable modules lying on some
homogeneous tube of $\Tub_\lam$ and having the same quasi-length~$\ell$
have the same class $\cla{h}_{a,b}$ in $K_0(\Da)$, where $b/a=\lam$
and $\gcd(a,b)=\ell$.
Moreover for $i=2, 3, 6$ the element $\cla{h}_{a,b}$ is also the class of  
$i$ modules with quasi-length $\ell$ lying on the tube of 
rank $i$ of $\Tub_\lam$.
It follows that there is no positive root $\cla{d}$ with $\rk(\cla{d})=i>1$
and $\ql(\cla{d})$ divisible by $i$.
\end{proof}

%%%%%%%%%%%%%%%%%%%%%%%%%%%%%%%%%%%%%%%%%
\subsection{} %Description of $R^+$
%%%%%%%%%%%%%%%%%%%%%%%%%%%%%%%%%%%%%%%%%
Let $i\in\{2, 3, 6\}$.
For $0\le m, n\le i-1$ and $(m,n) \not = (0,0)$ set 
\[
R_{[m,n]}(i) = R^\lam_\ell(i)
\]
where $\lam = n/m$ and $\ell = \gcd(m,n)$.
By Lemma~\ref{lemRli}, $R_{[m,n]}(i)$ has cardinality $i$.
In Section~\ref{rootlist}
we list all elements of the sets $R_{[m,n]}(i)$,
that is, a distinguished set of 
\[
2\cdot 3 + 3\cdot 8 + 6\cdot 35 =  240
\]
real elements, say $\cla{r}_1, \ldots, \cla{r}_{240}$, of $R^+$.

We constructed this set in the following way. 
Note first, that it coincides with the set
\[
\{\cla{r}\in R^+\mid 
0\leq \bil{\cla{h}_0,\cla{r}},\bil{\cla{r},\cla{h}_\infty}\leq 5\}.
\]
Now consider the algebra $\Da'$ which is obtained by restricting $\Da$ to the
full subquiver of $\tQa{0,1}$ which is obtained by removing the vertices $1_1$
and $4_1$. Notice that if $\cla{r}$ is a real root, then 
$\cla{r}-r_{1_1}\cla{h}_0- r_{4_1}\cla{h}_\infty$ is a real root with support
in $\Da'$. This is tilted of type $\E_8$, thus the associated quadratic form
$q'$ has 240 roots $\dva{r}'_1,\ldots,\dva{r}'_{240}$. Find them for example as
the orbits of the dimension vectors of the $8$ indecomposable projective
$\Da'$-modules under the corresponding Coxeter transformation $\Phi'$.
Since $\Da'$ is tilted of type $\E_8$, the transformation $\Phi'$
has order $30$ (the Coxeter number of $\E_8$). Next, define integers
$a'_i, a''_i, b'_i, b''_i$ by
\begin{align*}
\bil{\cla{r}'_i,\cla{h}_\infty} =& 6a'_i + a''_i \text{ with } 0\leq a''_i\leq 5,\\
\bil{\cla{h}_0,\cla{r}'_i}      =& 6b'_i + b''_i \text{ with } 0\leq b''_i\leq 5, 
\end{align*}
for $1\leq i\leq 240$. Finally set
$
\cla{r}_i=\cla{r}'_i- a'_i\cla{h}_0-b'_i\cla{h}_\infty.
$

We can now give the following explicit construction of all
real roots in $R^+$.

\begin{Thm}[Construction of $R^\lam_\ell(i)$]
\label{Thm3}
Let $\lam \in \Q_\infty$, $i\in\{ 2,3,6\}$
and $\ell \in \ENp$ not divisible by $i$.
Write $\lam = b/a$ with $a \in \ENn$, $b \in \Z$ and $\gcd(a,b) = \ell$
(if $\lam = \infty$ set $a=0$ and $b=\ell$).
Write
\[
a = ia' + a'', \qquad b = ib' + b''
\]
with
$a',a'',b',b'' \in \Z$ such that $0 \leq a'',b'' \leq i-1$.
Then 
\[
R^\lam_\ell(i) = 
\{ \dva{h}_{a',b'} + \dva{r} \mid \dva{r} \in R_{[a'',b'']}(i) \}.
\]
\end{Thm}
\begin{proof}
Let $\cla{r}$ be a real root of rank $i$, $\cla{h}$ an imaginary root
and set $\cla{d}=\cla{r}+\cla{h}$. Then $\cla{d}$ is a real root and 
$\Phi(\cla{d})=\Phi(\cla{r})+\cla{h}$,
hence $\cla{d}$ has also rank $i$.
Moreover, if $\cla{h}=\cla{h}_{a',b'}$ and 
$\cla{r}\in R_{[a'',b'']}(i)$ then
$h(\cla{d})=i\cla{h}_{a',b'}+\cla{h}_{a'',b''}=\cla{h}_{a,b}$, therefore 
$\cla{d}$ belongs to $R^\lam_\ell(i)$.
Thus, since $R^\lam_\ell(i)$ and $R_{[a'',b_i]}(i)$ both
have cardinality $i$ we see that the first subset is the translate
of the second one by $\cla{h}_{a',b'}$.
\end{proof}
Thus, the 240 positive real roots listed in Section~\ref{rootlist} yield
a complete description of the infinite set of all positive real roots of $R$.
Note that the classes of these 240 roots in $\Z^{10}/\rad(q)$ form a finite
root system of type $\E_8$, and we recover that
$R$ is an elliptic root system of type $\E_8^{(1,1)}$.

%%%%%%%%%%%%%%%%%%%%%%%%%%%%%%%%%%%%%%%%%%%%%%%%%%%%%%%
\subsection{}  % Schur roots
%%%%%%%%%%%%%%%%%%%%%%%%%%%%%%%%%%%%%%%%%%%%%%%%%%%%%%% 
The set $R_S^+$ of {\it Schur roots} is defined as 
\[
R_S^+ = \{ \cla{d}\in R^+ \mid 
\gcd(\bil{\cla{h}_0,\cla{d}},\bil{\cla{d},\cla{h}_\infty}) \leq 6 \}.
\]
It is easy to see that the set of imaginary Schur roots consists
of the $h_{a,b}$ with $\gcd(a,b)=1$.

The Schur roots can also be characterized in terms of rank and
quasi-length, namely
\[
R_S^+ = \{ \cla{d} \in R^+ \mid \ql(\cla{d}) \leq \rk(\cla{d}) \}.
\]
This comes from the identity
\[
\gcd(\bil{\cla{h}_0,\cla{d}},\bil{\cla{d},\cla{h}_\infty})=
{6\,\frac{\ql(\cla{d})}{\rk(\cla{d})}}.
\]
Thus the set of real Schur roots is equal to the union of
all subsets $R^\lam_\ell(i)$ for $\lam\in \Q_\infty$,
$i=2, 3, 6$ and $1\le \ell<i$.
Using \ref{R-lambda-l-i}, this implies that $R_S^+ \cap R^\lam$ 
contains exactly 
\[
5 \cdot 6 + 2 \cdot 3 + 1 \cdot 2 + 1 = 39
\]
roots, one imaginary and the others real.
Note that all the 240 roots listed in Section~\ref{rootlist}
are Schur roots.

The Schur roots of $R^+$ are related to the Schur roots of
$\tLam$ (see \ref{IS}) in the following way.
Suppose  that $M\in\cHe$ is non-projective and
has quasi-length $\ell$ in a tube of rank $r$. 
Since all tubes in  $\Gam_\tLam$ are standard,
the endomorphism ring of $M$ is non-trivial if and only if
either  $\ell>r$, or $\ell=r$ and $M$ has a non-trivial 
endomorphism that factors over a 
projective module in the same tube. 

\begin{Cor} \label{cor:triven}
A module
$M\in\cHe$  has trivial endomorphism ring if and only if $[M]=0$, or
$[M]\in R_S^+$ and $\dimv{M}$ does not belong to the following list:
\[
\dva{h}_\infty +\dvp{i},\quad \dva{h}_0+\dvp{j}
\]
with $i=1,3,5$ and $j=2,4$. 
\end{Cor}
\proof
This follows from~\eqref{eqn:ql} together with our description 
of $\Tub_0$ and $\Tub_\infty$.
\qed

%%%%%%%%%%%%%%%%%%%%%%%%%%%%%%%%%%%%%%%%%%%%%%%%%%%%%%%%%%%%

\section{Case $\A_5$: parametrization of the indecomposable 
irreducible components}\label{category}

%%%%%%%%%%%%%%%%%%%%%%%%%%%%%%%%%%%%%%%%%%%%%%%%%%%%%%%%%%%%

We shall now explain how $R^+$ parametrizes 
(i) the indecomposable $\Lam$-modules and (ii) the set
of dimension vectors of indecomposable $\tLam$-modules modulo
the Galois group action.
From (ii) we shall deduce the main result of this section,
namely the parametrization of the indecomposable
irreducible components of $\Lam$-modules by $R^+_S$.

%%%%%%%%%%%%%%%%%%%%%%%%%%%%%%%%%%%%%%%%%%%%%%%
\subsection{} % AR-quiver of \Lam
%%%%%%%%%%%%%%%%%%%%%%%%%%%%%%%%%%%%%%%%%%%%%%%
By \ref{ARtLam}, the indecomposable $\Lam$-modules are in one-to-one
correspondence with the $\tLam$-modules of the class $\cHe$.
This class decomposes into tubular families $\Tub_\lam$ as shown
in \ref{expldec}.
The projective modules $P_1$ and $P_5$ appear at the mouth of
the non-homogeneous tube $T_{\infty,x_1}$,
the module $P_3$ at the mouth of
the non-homogeneous tube $T_{\infty,x_2}$,
and the modules $P_2$ and $P_4$ at the mouth of
the non-homogeneous tube $T_{0,x_1}$.
As a result we obtain the following parametrization
of the indecomposable $\Lam$-modules by $R^+$.

\begin{Prop}\label{rootsdimvectors}
Let $\lam \in \Q_\infty$, $i \in \{ 1,2,3,6 \}$ and 
$\ell \in \ENp$.
Then the following hold:
\begin{itemize}

\item[(a)]
If $i$ does not divide $\ell$
then there exists a one-to-one correspondence 
between $R_\ell^\lam(i)$ and the set of indecomposable $\Lam$-modules
of quasi-length $\ell$ in the
non-homogeneous tube~$T_{\lam,x_j}$ with $m_j = i$.
This correspondence maps $\cla{d} \in R^\lam_\ell(i)$ to $M=F(N)$
with $[N] = \cla{d}$.

\item[(b)]
If $i = 1$ then
$R_\ell^\lam(1) = \{ \cla{h}_{a,b} \}$ where $b/a = \lam$ and
$\gcd(a,b) = \ell$.
There is an infinite set of 
indecomposable $\Lam$-modules $M=F(N)$ with $[N] = \cla{h}_{a,b}$,
parametrized by the weighted projective line $\X$.
More precisely, for each ordinary point $x\in\X$ there is an
indecomposable module of quasi-length $\ell$ in the tube $T_{\lam,x}$,
and for each exceptional point $x_j$, there are $m_j$ indecomposable
modules of quasi-length $\ell$ on $T_{\lam,x_j}$. 

\item[(c)]
The only indecomposable $\Lam$-modules not appearing in the above
lists are the five indecomposable projective modules. 

\end{itemize}
\end{Prop}

%%%%%%%%%%%%%%%%%%%%%%%%%%%%%%%%
\subsection{}
%%%%%%%%%%%%%%%%%%%%%%%%%%%%%%%%
Recall from Corollary~\ref{cor:eta} that we have a canonical map
$\xi\colon H \to R^+\cup \{0\}$.
We are going to define a `right inverse' $\delta$ of $\xi$.
Define
$\udel: R^+ \lra K_0(\tLam)$
by
\[
\udel(\cla{r}) =
\begin{cases}
\cla{r} -\min\{0,r_{2_2}\}\dvp{2} -\min\{0, r_{4_1}\}\dvp{4} &
\text{ if $\dva{r}\in R^{\infty}$},\\
\cla{r} &\text{ if $\cla{r}\in R^{++}$},\\
\cla{r}-\min\{0,r_{1_1}\}\dvp{1}-\min\{0,r_{3_0}\}\dvp{3}-\min\{0,r_{5_{-1}}\}\dvp{5}&
\text{ if $\cla{r}\in R^0$},\\
\cla{r}- r_{2_2}\dvp{2}- r_{4_1}\dvp{4} &\text{ if $\dva{r}\in R^{+-}$}.
\end{cases}
\]

\begin{Prop} \label{deltares}
With the above definition of $\udel$ we have:
\begin{itemize}

\item[(a)] 
$\udel$ induces a well-defined and injective map
$\delta\colon{R^+}\lra H \cong \Idv(\tLam)/\!\sim$;

\item[(b)] 
The only elements of $H$ not in the image of $\delta$
are
\[
\ENn\,\dva{h}_\infty +\dvp{i},\quad \ENn\,\dva{h}_0+\dvp{j}
\]
with $i=1,3,5$ and $j=2,4$;

\item[(c)] 
The map $\delta$ restricts to a map
$\delta_S: R_S^+ \lra H_S=H\cap R^+_S$. 
The only elements of $H_S$ not in the image of $\delta_S$ are
$\dvp{1},\ldots,\dvp{5}$.

\end{itemize}

\end{Prop}

\begin{proof} 
First we have to show that $\udel(\cla{r})$
is the dimension vector of an indecomposable $\tLam$-module. This is clear for
$\cla{r}\in R^{++}$. 
For $\cla{r}\in R^{+-}$ we notice that 
$\cla{x}\mapsto \cla{x}- x_{2_2}\dvp{2}- x_{4_1}\dvp{4}$ gives an isometry
$K_0(\Da)\lra K_0(\Daa{-1,0})$, so that $\udel(R^{+-})$ consists of 
dimension vectors of objects in $\cMa{-1,0}$. The remaining two cases are
treated directly. It is easy to calculate the map $\xi$ 
(Corollary~\ref{cor:eta})
explicitly with~\eqref{eqn:decK}. 
It follows that $\xi\udel=\1_{R^+}$
and the rest of (a) follows since $H$ is a transversal for the action of
$\Z$ on $\ind({\tLam})$, see Proposition~\ref{prp:transv}.
Now we obtain $(b)$ from the description of the fibres of $\xi$ 
in Corollary~\ref{cor:eta}, and (c) follows from Corollary~\ref{cor:triven}.
\end{proof}

%%%%%%%%%%%%%%%%%%%%%%%%%%%%%%%%%%%%%%
\subsection{}
%%%%%%%%%%%%%%%%%%%%%%%%%%%%%%%%%%%%%%
Collecting the results of Theorem~\ref{ThmGS}, Theorem~\ref{Thm1}
and Proposition~\ref{deltares} we can now state the following
parametrization of the indecomposable irreducible components
of varieties of $\Lam$-modules and of the corresponding 
multisegments.
Let $C_j\ (1\le j \le 5)$ be the irreducible components containing 
the five indecomposable projective $\Lam$-modules.
Let $\msp_j$ denote the corresponding multisegments, namely
\[
\msp_j = \sum_{i=1}^{5-j+1} [i,i+j-1],\qquad (j=1,\ldots,5).
\]

\begin{Thm}\label{paramet}
\begin{itemize}
\item[(a)] The map 
\[
\cla{d} \mapsto \eta(Z_{\delta(\cla{d})})
\] 
is a one-to-one
correspondence from the set $R_S^+$ of Schur roots of the Ringel
form $\bil{-,-}$ in $\Z^{10}$ to the set 
$\ind(\irr(\Lam))-\{C_1,\ldots ,C_5\}$.
\item[(b)] 
The map 
\[
\cla{d} \mapsto \psi(\delta(\cla{d}))
\]
is a one-to-one
correspondence from $R_S^+$ to 
$\ind(\MM) - \{\msp_1,\ldots ,\msp_5\}$.
\end{itemize}
\end{Thm}
Note that the descriptions of $R_S^+, \delta, \eta$ and $\psi$
are completely explicit, so that we get a very concrete
parametrization of the factors arising in the canonical factorization
of the elements of ${\mathcal S}^*$.

\begin{example}{\rm
(i)\quad Let 
$\cla{d}=\left( 
\bsm &1&&0\\2&&3&&1\\&3&&3\\1&&2&&1 \esm \right)
\in R^{++}$.
Then $\delta(\cla{d}) = \dva{d}$ and 
\[
\psi(\delta(\cla{d}))=[1,1]+[1,2]+[1,3]+2\,[2,3]+[3,4]+[3,5]+[4,4]+[5,5].
\]
(ii)\quad  Let 
$\cla{d}=\left( 
\bsm &-1&&0\\0&&0&&0\\&1&&0\\0&&0&&0 \esm \right)
\in R^0$.
Then $\delta(\cla{d})=\left( 
\bsm &0&&0\\1&&1&&0\\&2&&1\\0&&1&&1\\
&0&&1 \esm \right)$
and 
\[
\psi(\delta(\cla{d}))=[1,1]+[2,2]+[2,3]+[3,4]+[4,5].
\]
}
\end{example}

%%%%%%%%%%%%%%%%%%%%%%%%%%%%%%%%%%%%%%%%%%%%%%%%%%%

\section{Case $\A_5$: the component graph}\label{extensions}

%%%%%%%%%%%%%%%%%%%%%%%%%%%%%%%%%%%%%%%%%%%%%%%%%%%

Recall that the component graph ${\mathcal C}(\Lam)$ has for vertices 
the indecomposable irreducible components of the 
varieties of $\Lam$-modules, and two vertices $Z_1$ and $Z_2$
are connected by an edge if and only if
$\overline{Z_1\oplus Z_2}$ is an irreducible component,
or equivalently $\ext_\Lam^1(Z_1,Z_2)=0$.
There are edges from the irreducible components 
$C_i\ (1\le i\le 5)$ to every other vertex.
The following theorem describes all remaining edges.
In agreement with Theorem~\ref{paramet} (a), we shall
label the vertices other than $C_i$ by the elements of $R^+_S$.

In order to state the theorem we introduce the following 
definition. We call a pair of Schur roots $(\cla{d},\cla{e})$ {\em critical}
if the following three conditions hold:
\begin{itemize}
\item $\{\cla{d},\cla{e}\}\subset R_\ell^\mu(6)$ for some 
$\mu\in\Q_\infty$ and $\ell \in \Z^+$,
\item $\bil{\cla{d},\cla{e}}=0=\bil{\cla{e},\cla{d}}$,
\item $\ql(\cla{d})+\ql(\cla{e})\geq 7$.
\end{itemize}

\begin{Thm}\label{Thm4}
Two Schur roots 
$\cla{d}$ and $\cla{e}$ are connected by an edge in 
${\mathcal C}(\Lam)$ if and only if the following two conditions hold:
\begin{itemize}

\item[(i)]
$\bil{\cla{d},\cla{e}}\geq 0$ and $\bil{\cla{e},\cla{d}}\geq 0$,

\item[(ii)] 
$(\cla{d},\cla{e})$ is not critical, or
$\bil{\cla{d},\Phi^i(\cla{e})} < 0$ where
$i = \min \{ j \geq 1 \mid \bil{\cla{d},\Phi^j(\cla{e})} \not= 0 \}$.

\end{itemize}
\end{Thm}
\proof This follows from an adaptation of \cite[Theorem 1.3, Lemma
6.4]{GSc}. \qed

Thus, the edges of the component graph of $\Lam$ are completely 
determined by the bilinear form $\bil{-,-}$ and the Coxeter
matrix $\Phi$.
Moreover there is an edge between $d$ and $e$ if and only if
there is an edge between $\Phi(\cla{d})$ and $\Phi(\cla{e})$.

%%%%%%%%%%%%%%%%%%%%%%%%%%%%%%%%%%%%%%%%%%%%%%%%

\section{Proof of Theorem~\ref{Thm1}}\label{proof1}

%%%%%%%%%%%%%%%%%%%%%%%%%%%%%%%%%%%%%%%%%%%%%%%%

It was already proved in \cite{GSc} that
$\theta$ and $\eta$ are well-defined and bijective.
The map $\phi$ is bijective by definition.
It remains to explicitly construct the map
$\psi = \phi^{-1} \eta \theta$.
It is enough to prove Theorem~\ref{Thm1} for $\Lam = \Lam_5$, since
$\Lam_n$ ($n=2,3,4$) are full convex subalgebras of $\Lam_5$.

We will use the following result from \cite[Theorem 1, 2]{BS}: 
Let $A$ be a tame quasi-tilted basic algebra, and let $\cla{d}$ be a 
dimension vector of an indecomposable $A$-module.
Then $\md(A,\cla{d})$ has at most two irreducible components, and 
$\md(A,\cla{d})$ is irreducible if and only if
$\cla{d}$ is not of one of the following forms:
\begin{itemize}

\item[(a)] $\cla{d} = \cla{h}+\cla{z}$ where $\cla{h}$ and $\cla{z}$ are 
connected positive vectors with disjoint
support, $\bil{\cla{h},\cla{h}}_A = 0$, $\bil{\cla{z},\cla{z}}_A = 1$ 
and $z_i \leq 1$ for all entries
$z_i$ of $\cla{z}$,

\item[(b)] $\cla{d} = \cla{h}+\cla{h}'$ where $\cla{h}$ and $\cla{h}'$ are 
connected positive vectors with
$\bil{\cla{h},\cla{h}}_A = \bil{\cla{h}',\cla{h}'}_A = 0$, 
$\bil{\cla{h},\cla{h}'}_A = 1$ and $\bil{\cla{h}',\cla{h}}_A = 0$.
\end{itemize}
The algebras $\Da=\Daa{0,1}$ and $\Da^*=\Daa{-1,0}$ are both tubular algebras, 
in particular they are tame quasi-tilted algebras.
All connected positive vectors $\cla{h}$ for $\Daa{i,i+1}$  such that
$\bia{i,i+1}{\cla{h},\cla{h}}= 0$ 
are of the form 
$a\bha{i} + b\bha{i+1}$ where 
$(0,0) \not= (a,b) \in \ENn \times \ENn$.
We have $\bia{i,i+1}{\bha{i},\bha{i+1}}= 6$ and 
$\bia{i,i+1}{\bha{i+1},\bha{i}} = -6$.
It follows that the case $(b)$ above cannot occur for 
$A = \Daa{i,i+1}$, $i=-1,0$.

For $A = \Da$ we are in case $(a)$ precisely when
$\cla{d} \in \{ \dva{e}_1(n), \ldots, \dva{e}_5(n) \mid n \geq 1 \}$
where
\[
\dva{e}_1(n) = \left( \bsm &1&&0\\n&&2n&&n\\&3n&&3n\\n&&2n&&n \esm \right),\,
\dva{e}_2(n) = \left( \bsm &0&&1\\n&&2n&&n\\&3n&&3n\\n&&2n&&n \esm \right),
\]
and
\[
\dva{e}_3(n) = \left( \bsm &n&&n\\n&&2n&&n\\&n&&n\\1&&0&&0 \esm \right),\,
\dva{e}_4(n) = \left( \bsm &n&&n\\n&&2n&&n\\&n&&n\\0&&1&&0 \esm \right),\,
\dva{e}_5(n) = \left( \bsm &n&&n\\n&&2n&&n\\&n&&n\\0&&0&&1 \esm \right).
\]
For $A = \Da^*$ we are in case $(a)$ precisely when 
$\dva{d} \in \{ (\dva{e}_1(n))^*, \ldots, (\dva{e}_5(n))^* \mid n \geq 1 \}$
where
\[
\dva{e}_1(n)^* = \left( \bsm n&&2n&&n\\&3n&&3n\\n&&2n&&n\\&0&&1 \esm \right),\,
\dva{e}_2(n)^* = \left( \bsm n&&2n&&n\\&3n&&3n\\n&&2n&&n\\&1&&0 \esm \right),
\]
and
\[
\dva{e}_3(n)^* = \left( \bsm 0&&0&&1\\&n&&n\\n&&2n&&n\\&n&&n \esm \right),\,
\dva{e}_4(n)^* = \left( \bsm 0&&1&&0\\&n&&n\\n&&2n&&n\\&n&&n \esm \right),\,
\dva{e}_5(n)^* = \left( \bsm 1&&0&&0\\&n&&n\\n&&2n&&n\\&n&&n \esm \right).
\]
Let $M$ be an indecomposable $\tLam$-module, and let
$\dimv{M}$ be its dimension vector.
If  
\[
\dimv{M} \in \{ \dva{e}_i(n),\, (\dva{e}_i(n))^* \mid i = 1,2, n \geq 1 \},
\]
then $F(M) \in \Tub_0$, and if
\[
\dimv{M} \in \{ \dva{e}_i(n),\, (\dva{e}_i(n))^* \mid i = 3,4,5, n \geq 1 \},
\]
then $F(M) \in \Tub_\infty$.
One easily checks that $\dva{e}_i(n)$ is a Schur root if and only if
$n = 1$.
Set $\dva{e}_i=\dva{e}_i(1)$ and $\dva{e}_i^*=(\dva{e}_i(1))^*$, $1 \leq i \leq 5$.

Assume that $\dva{d} \in \Idv_S(\tLam)$ is a Schur root whose support
lies (up to shift) in $\Da$ or $\Da^*$, and assume
$\dva{d} \not\in \{ \dva{e}_i,\, \dva{e}_i^* \mid 1 \leq i \leq 5 \}$.
By the result mentioned above this implies that $\md(\tLam,\dva{d})$ 
is irreducible.
Thus $Z_\dva{d} = Z_{\rm max}(\dva{d})$ and
\[
\psi(\dva{d}) = \mu(Z_\dva{d}) = \msm_{\rm max}(\dva{d}).
\]

Next, assume $\dva{d} \in \{ \dva{e}_i,\, \dva{e}_i^* \mid 1 \leq i \leq 5 \}$.
Thus $\md(\tLam,\dva{d})$ has exactly two irreducible components.
Furthermore, we know that 
\[
Z_\dva{d} = \overline{\orb(M_\dva{d})}
\]
for some indecomposable $\tLam$-module $M_\dva{d}$.

For any vertex $i_j$ of $\tQ$ let $S_{i_j}$ be the
corresponding simple $\tLam$-module, and let $Z_{i_j}$ 
be the irreducible component consisting of the single point
corresponding to $S_{i_j}$.
Then 
\[
\ext_{\tLam}^1(Z_{\dva{h}_0},Z_{i_j}) = \ext_{\tLam}^1(Z_{i_j},Z_{\dva{h}_0}) = 0
\]
for $i = 2,4$ and all $j \in \Z$, and
\[
\ext_{\tLam}^1(Z_{\dva{h}_\infty},Z_{i_j}) = 
\ext_{\tLam}^1(Z_{i_j},Z_{\dva{h}_\infty}) = 0
\]
for $i = 1,3,5$ and all $j \in \Z$. 
Thus \cite[Theorem 1.2]{CBSc} implies that 
$\overline{Z_{\dva{h}_0} \oplus Z_{i_j}}$ (resp. 
$\overline{Z_{\dva{h}_\infty} \oplus Z_{i_j}}$) are irreducible components provided
$i = 2,4$ (resp. $i = 1,3,5$).
Exactly one of these irreducible components lies
in $\md(\tLam,\dva{d})$, we denote this component by $Z_{\rm dec}(\dva{d})$.
Thus $\md(\tLam,\dva{d})$ contains
exactly one indecomposable irreducible component, namely $Z_\dva{d}$, 
and exactly one 
decomposable irreducible component, namely $Z_{\rm dec}(\dva{d})$.

If $\overline{Z' \oplus Z''}$ is an irreducible component, then 
\[
\mu(\overline{Z' \oplus Z''}) = \mu(Z') + \mu(Z'').
\]
Thus, if $\dva{d} \not\in \{ \dva{e}_3^*,\, \dva{e}_5 \}$, then 
$Z_{\rm dec}(\dva{d}) \not= Z_{\rm max}(\dva{d})$.
This yields $Z_\dva{d} = Z_{\rm max}(\dva{d})$ and
\[
\psi(\dva{d}) = \mu(Z_\dva{d}) = \msm_{\rm max}(\dva{d}).
\]
For $\dva{d} \in \{ \dva{e}_3^*,\, \dva{e}_5 \}$ it is not difficult to construct
the module $M_\dva{d}$ explicitly.
We get
\[
\psi(\dva{e}_3^*) = 2\,[1,1] + [2,2] + [2,4] + [3,3] + [4,5]
\]
and
\[
\psi(\dva{e}_5) = [1,2] + [2,4] + [3,3] + [4,4] + 2\,[5,5].
\]

The only Schur roots in $\Idv_S(\tLam)$ whose support is (up to shift) 
not contained in $\Da$ or $\Da^*$ are 
\[
\dva{f}_1 = \left( \bsm &0&&1\\1&&2&&1\\&3&&2\\1&&2&&1\\&0&&1 \esm \right),\,
\dva{f}_2 = \left( \bsm &1&&0\\1&&2&&1\\&2&&3\\1&&2&&1\\&1&&0 \esm \right),\,
\dvp{2},\, \dvp{4} 
\] 
(these are all in $\Tub_0$), and
\[
\dva{g}_1 = \left( \bsm 1&&0&&0\\&1&&1\\0&&2&&1\\&1&&1\\1&&0&&0 \esm \right),\,
\dva{g}_2 = \left( \bsm 0&&0&&1\\&1&&1\\1&&2&&0\\&1&&1\\0&&0&&1 \esm \right),\,
\dvp{1},\, \dvp{3},\,\dvp{5} 
\]
(these are all in $\Tub_\infty$).
Here $\dvp{1}, \ldots, \dvp{5}$ are (up to shift) the dimension vectors of the
indecomposable projective $\tLam$-modules as displayed in~\ref{sus:isomt}.

For each 
$\dva{d}\in\{\dvp{1},\ldots,\dvp{5},\dva{f}_1,\dva{f}_2,\dva{g}_1,\dva{g}_2 \}$, 
there exists an
indecomposable $\tLam$-module $M_\dva{d}$ such that 
$Z_\dva{d} = \overline{\orb(M_\dva{d})}$.
For $\dva{d} \in \{ \dvp{1}, \ldots, \dvp{5},\dva{g}_1,\dva{g}_2 \}$ it is easy 
to construct $M_\dva{d}$ explicitly.
In these cases, it follows that $Z_\dva{d} = Z_{\rm max}(\dva{d})$, thus
\[
\psi(\dva{d}) = \mu(Z_\dva{d}) = \msm_{\rm max}(\dva{d}).
\]

For any element $a \in \tLam$ and any $\tLam$-module $M$,
let $f_{a,M}: M \to M$ be the linear map defined by the $\tLam$-module action 
of $a$ on $M$, i.e. $f_{a,M}(m) = am$.
If $M$ is a submodule of a module $N$, then we have
\[
\rk(f_{a,M}) \geq \rk(f_{a,N})
\]
for all $a$.
Now let 
\[
\dva{d} = \dva{f}_1 = 
\left( \bsm &0&&1\\1&&2&&1\\&3&&2\\1&&2&&1\\&0&&1 \esm \right).
\]
In Section~\ref{Tnull} we can see that $M_\dva{d}$ contains two 
indecomposable submodules $M_{\dva{f}_1'}$ and $M_{\dva{f}_1''}$
with dimension vectors
\[
\dva{f}_1' = \left( \bsm &0&&0\\1&&2&&1\\&3&&2\\1&&2&&1\\&0&&1 \esm \right)  
\text{ and }
\dva{f}_1'' = \left( \bsm &0&&1\\0&&1&&1\\&2&&1\\1&&1&&0\\&0&&0 \esm \right),
\]
respectively.
We get
\[
Z_{\dva{f}_1'} = \overline{\orb(M_{\dva{f}_1'})} 
\text{ and }
Z_{\dva{f}_1''} = \overline{\orb(M_{\dva{f}_1''})}.
\]
Since the support of $\dva{f}_1'$ and $\dva{f}_1''$ lies in $\Da^*$ and $\Da$, 
respectively,
we get
\[
\psi(\dva{f}_1') = \mu(Z_{\dva{f}_1'}) = \msm_{\rm max}(\dva{f}_1') 
\text{ and } 
\psi(\dva{f}_1'') = \mu(Z_{\dva{f}_1''}) = \msm_{\rm max}(\dva{f}_1'').
\]
This enables us to compute the ranks of the maps $f_{p,M_{\dva{f}_1'}}$ and
$f_{p,M_{\dva{f}_1''}}$ for all paths $p$ in $\tQ$ of the form
$\alpha_{lj} \alpha_{l+1,j} \cdots \alpha_{mj}$, $1 \leq l \leq m \leq 4$.
Then we can use the above rank inequality, and get the rank of
$f_{p,M_{\dva{f}_1}}$ for any path $p$.
It turns out that
\[
\psi(\dva{f}_1) = \mu(Z_{\dva{f}_1}) = \msm_{\rm max}(\dva{f}_1).
\]
The case $\dva{d} = \dva{f}_2$ is done in a similar way, and we get again
\[
\psi(\dva{f}_2) = \mu(Z_{\dva{f}_2}) = \msm_{\rm max}(\dva{f}_2).
\]

%%%%%%%%%%%%%%%%%%%%%%%%%%%%%%%%%%%%%%%%%%%%%%%%%%%%%%%%

\section{Concluding remarks} \label{concl}

%%%%%%%%%%%%%%%%%%%%%%%%%%%%%%%%%%%%%%%%%%%%%%%%%%%%%%%%
%%%%%%%%%%%%%%%%%%%%%%%%%%%%
\subsection{} %type D4
%%%%%%%%%%%%%%%%%%%%%%%%%%%%

By Proposition~\ref{infinite_type}, the preprojective algebra 
$\Lam=P(Q)$ is tame 
if and only if the quiver $Q$ is of Dynkin type $\A_5$
or $\D_4$.
Using the same methods as in this paper it is possible
to obtain a complete analogue of Theorem~\ref{ThmC} for 
type $\D_4$. 
In this case, $\Lam$ has a Galois covering which is isomorphic
to the repetitive algebra of a tubular algebra of type
$(3,3,3)$, and the corresponding root system is an elliptic
root system of type $\E_6^{(1,1)}$. 
We plan to give a detailed analysis of this case in a 
forthcoming publication.

%%%%%%%%%%%%%%%%%%%%%%%%%%%%%%%
\subsection{} %cluster type 
%%%%%%%%%%%%%%%%%%%%%%%%%%%%%%%

It is shown in \cite{BFZ} that for all Dynkin types, 
the algebra $\C[N]$ has the structure of an (upper) cluster algebra, 
and that it has finite type as a cluster algebra if and only
$\g$ is of Lie type $\A_n \ (n\le 4)$.
In that case one can associate to $\C[N]$ a root system ${\mathcal R}$
called its {\em cluster type}, which controls the combinatorics
of the cluster variables and of the cluster sets.
More precisely, the cluster variables are parametrized
by the set ${\mathcal R}_{\ge -1}$ of almost positive roots 
of ${\mathcal R}$, and the
pairs of cluster variables which can occur simultaneously
in a cluster set can be explicitly described in terms of
${\mathcal R}$ and a piecewise linear Coxeter transformation 
acting on ${\mathcal R}_{\ge -1}$. 
The Cartan matrix $A$ of ${\mathcal R}$ can be obtained by
a certain symmetrization procedure from the principal part 
$B(t)$ of the exchange matrix of $\C[N]$ 
at certain vertices $t$ of its exchange graph (see \cite{FZ3}).

For $\g$ of Lie type $\A_2, \A_3, \A_4$, the algebra $\C[N]$
has cluster type $\A_1, \A_3, \D_6$ respectively \cite{BFZ}.
As mentioned in \ref{reDatla}, these root systems also occur 
in our setting in the following way.
For $n=2, 3, 4$, we have 
$\smd(\tLam_n)\cong D^b(\md(\C{\mathcal Q}_n))$
where ${\mathcal Q}_n$ is a quiver of type  $\A_1, \A_3, \D_6$ respectively.

%%%%%%%%%%%%%%%%%%%%%%%%%%%%%%
\subsection{} %type A5, D4
%%%%%%%%%%%%%%%%%%%%%%%%%%%%%%
\begin{figure}[t]
\begin{center}
\leavevmode
\epsfxsize =8cm
\epsffile{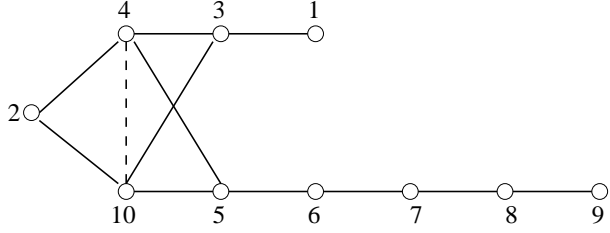}
\end{center}
\caption{\label{figu2} The Dynkin diagram of $\E_8^{(1,1)}$}
\end{figure}

At the moment, there is no notion of cluster type for
the algebras $\C[N]$ which are not of finite type as
cluster algebras. 
The results of this paper strongly suggest that if such
a cluster type exists for $\g$ of Lie type $\A_5$ ({\it resp.} $\D_4$),
then it should be the elliptic type $\E_8^{(1,1)}$ 
({\it resp.} $\E_6^{(1,1)}$) in the notation of Saito, or the
tubular type $(6,3,2)$ ({\it resp.} $(3,3,3)$) in the language of Ringel.
Remember in particular that 
$\smd(\tLam_5)\cong D^b(\coh(\X))$
where $\X$ is a weighted projective line of type $(6,3,2)$.

\subsection{}
Here is another remark supporting that guess.
For $\g$ of Lie type $\A_5$ one can find a vertex $t$
of the exchange graph of $\C[N]$ at which the 
principal part of the exchange matrix is 
(up to simultaneous permutations of rows and columns)
\[
B(t)=
\left(\begin{matrix} 
                 0  &  0 &  1  & 0  &  0 &  0  &  0  &  0  &  0  &  0\cr
                 0  &  0 &  0  & 1  &  0 &  0  &  0  &  0  &  0  &  -1\cr
                -1  &  0 &  0  & 1  &  0 &  0  &  0  &  0  &  0  &  -1\cr
                 0  &  -1&  -1 & 0  &  -1&  0  &  0  &  0  &  0  &  2\cr
                 0  &  0 &  0  & 1  &  0 & -1  &  0  &  0  &  0  &  -1\cr
                 0  &  0 &  0  & 0  &  1 &  0  &  -1 &  0  &  0  &  0\cr
                 0  &  0 &  0  & 0  &  0 &  1  &  0  & -1  &  0  &  0\cr
                 0  &  0 &  0  & 0  &  0 &  0  &  1  &  0  & -1  &  0\cr
                 0  &  0 &  0  & 0  &  0 &  0  &  0  &  1  &  0  &  0\cr
                 0  &  1 &  1  &-2  &  1 &  0  &  0  &  0  &  0  &  0
\end{matrix}\right) 
\]
We suggest to take as symmetric counterpart of $B(t)$ the matrix
\[
A=
\left(\begin{matrix} 
                 2  &  0 &  -1  & 0  &  0 &  0  &  0  &  0  &  0  &  0\cr
                 0  &  2 &  0  & -1  &  0 &  0  &  0  &  0  &  0  &  -1\cr
                -1  &  0 &  2  & -1  &  0 &  0  &  0  &  0  &  0  &  -1\cr
                 0  & -1 & -1  & 2  & -1 &  0  &  0  &  0  &  0  &  2\cr
                 0  &  0 &  0  & -1 &  2 &  -1  &  0   &  0  &  0  & -1\cr
                 0  &  0 &  0  & 0  & -1 &  2  & -1  &  0  &  0  &  0\cr
                 0  &  0 &  0  & 0  &  0 & -1  &  2  & -1  &  0  &  0\cr
                 0  &  0 &  0  & 0  &  0 &  0  & -1  &  2  &  -1 &  0\cr
                 0  &  0 &  0  & 0  &  0 &  0  &  0  & -1  &  2  &  0\cr
                 0  & -1 & -1  & 2  & -1 &  0  &   0  &  0  &  0  & 2
\end{matrix}\right). 
\]
Note that $A$ has two positive entries off the diagonal, so it differs
from the matrix obtained from $B(t)$ by the symmetrizing procedure 
of \cite{FZ3}.
It turns out that $A$ is exactly the Cartan matrix of the root system 
$\E_8^{(1,1)}$,
that is, the matrix of scalar products of a basis of simple roots
in the sense of Saito \cite{Sa}. It can be visualized with the
help of the Dynkin type diagram of Figure~\ref{figu2} 
in which an ordinary edge between $i$ and $j$ means that $a_{ij}=a_{ji}=-1$
and the dashed line between $4$ and $10$ means that $a_{4,10}=a_{10,4}=2$
(see \cite{Sa}). 

There is a similar fact for type $\D_4$ and $\E_6^{(1,1)}$.

%\newpage
%%%%%%%%%%%%%%%%%%%%%%%%%%%%%%%%%%%%%%%%%%%%%%%%%%%%%%%%%

\section{Pictures and tables}\label{pictures}

%%%%%%%%%%%%%%%%%%%%%%%%%%%%%%%%%%%%%%%%%%%%%%%%%%%%%%%%%

%%%%%%%%%%%%%%%%%%%%%%%%%%%%%%%%%%%%%%%%%%%%%%%%%%%%%
\subsection{The Auslander-Reiten quivers of $\Lam_2$, $\Lam_3$, 
$\Lam_4$}\label{AR}
%%%%%%%%%%%%%%%%%%%%%%%%%%%%%%%%%%%%%%%%%%%%%%%%%%%%%

\[
\xymatrix@-1.3pc{
&&&\\
{\bsm&0\\0\\&1\\1\esm} \ar[dr]\ar@{.}[dd] &&&& {\bsm&1\\1\\&0\\0\esm}
\ar@{.}[dd]\\
&{\bsm&0\\0\\&1\\0\esm} \ar[dr] && {\bsm&0\\1\\&0\\0\esm} \ar[ur]\\
&&{\bsm&0\\1\\&1\\0\esm} \ar[ur]&&\\
}
\]
\begin{center}
The Auslander-Reiten quiver of $\Lam_2$
\end{center}

%%%%%%%%%%%%%%%%%%%%%%%%%%%%%%%%%

\[
\xymatrix@-1.3pc{
&&&\\
{\bsm0&&0\\&0\\1&&0\\&1\\0&&1\esm} \ar[dr]\ar@{.}[dd] &&&&&&
{\bsm0&&1\\&1\\1&&0\\&0\\0&&0\esm}\ar@{.}[dd]\\
&{\bsm0&&0\\&0\\1&&0\\&1\\0&&0\esm} \ar[dr] &&
{\bsm0&&0\\&0\\0&&1\\&0\\0&&0\esm} \ar[dr] &&
{\bsm0&&0\\&1\\1&&0\\&0\\0&&0\esm} \ar[dr]\ar[ur]\\
{\bsm0&&0\\&0\\0&&0\\&1\\0&&0\esm} \ar[dr]\ar[ur]\ar@{.}[dd] &&
{\bsm0&&0\\&0\\1&&1\\&1\\0&&0\esm} \ar[dr]\ar[r]\ar[ur] &
{\bsm0&&0\\&1\\1&&1\\&1\\0&&0\esm} \ar[r] &
{\bsm0&&0\\&1\\1&&1\\&0\\0&&0\esm} \ar[dr]\ar[ur] &&
{\bsm0&&0\\&1\\0&&0\\&0\\0&&0\esm}\ar@{.}[dd]\\
&{\bsm0&&0\\&0\\0&&1\\&1\\0&&0\esm} \ar[ur] &&
{\bsm0&&0\\&0\\1&&0\\&0\\0&&0\esm} \ar[ur] &&
{\bsm0&&0\\&1\\0&&1\\&0\\0&&0\esm} \ar[dr]\ar[ur]\\
{\bsm0&&0\\&0\\0&&1\\&1\\1&&0\esm} \ar[ur] &&&&&&
{\bsm1&&0\\&1\\0&&1\\&0\\0&&0\esm}\\ 
}
\]
\begin{center}
The Auslander-Reiten quiver of $\Lam_3$
\end{center}

%%%%%%%%%%%%%%%%%%%%%%%%%%%%%%%%%%%%%%%%%%%%%%%%%%%%%%%%%%%%%%%%
\newpage
\begin{landscape}
\[
\xymatrix@-1.3pc{
\ar@{.}[d]&&&{\bsm0&&0\\&0&&1\\0&&1\\&1&&0\\1&&0\esm}\ar[dr] &&&&&& {\bsm1&&0\\&1&&0\\0&&1\\&0&&1\\0&&0\esm}\ar[dr] &&&\ar@{.}[d]\\
{\bsm0&&0\\&0&&0\\0&&0\\&0&&1\\0&&0\esm}\ar[dr]\ar@{.}[dd] && 
{\bsm0&&0\\&0&&0\\0&&1\\&1&&0\\1&&0\esm}\ar[ur]\ar[dr] && 
{\bsm0&&0\\&0&&1\\0&&1\\&1&&0\\0&&0\esm}\ar[dr] && 
{\bsm0&&0\\&0&&0\\1&&0\\&0&&0\\0&&0\esm}\ar[dr] &&
{\bsm0&&0\\&1&&0\\0&&1\\&0&&1\\0&&0\esm}\ar[dr]\ar[ur] && 
{\bsm1&&0\\&1&&0\\0&&1\\&0&&0\\0&&0\esm}\ar[dr] &&
{\bsm0&&0\\&0&&1\\0&&0\\&0&&0\\0&&0\esm} \ar@{.}[dd]\\
&{\bsm0&&0\\&0&&0\\0&&1\\&1&&1\\1&&0\esm}\ar[ur]\ar[dr]&& 
{\bsm0&&0\\&0&&0\\0&&1\\&1&&0\\0&&0\esm}\ar[ur]\ar[dr] && 
{\bsm0&&0\\&0&&1\\1&&1\\&1&&0\\0&&0\esm}\ar[ur]\ar[dr] && 
{\bsm0&&0\\&1&&0\\1&&1\\&0&&1\\0&&0\esm}\ar[ur]\ar[dr] &&
{\bsm0&&0\\&1&&0\\0&&1\\&0&&0\\0&&0\esm}\ar[ur]\ar[dr]&& 
{\bsm1&&0\\&1&&1\\0&&1\\&0&&0\\0&&0\esm}\ar[ur]\ar[dr]&\\
{\bsm0&&0\\&0&&0\\1&&1\\&2&&1\\1&&1\esm} \ar[dr] \ar[ur]\ar@{.}[d] && 
{\bsm0&&0\\&0&&0\\0&&1\\&1&&1\\0&&0\esm} \ar[dr] \ar[ur] &&
{\bsm0&&0\\&0&&0\\1&&1\\&1&&0\\0&&0\esm} \ar[dr] \ar[ur] && 
{\bsm0&&0\\&1&&1\\1&&2\\&1&&1\\0&&0\esm} \ar[dr] \ar[ur] &&
{\bsm0&&0\\&1&&0\\1&&1\\&0&&0\\0&&0\esm} \ar[dr] \ar[ur] && 
{\bsm0&&0\\&1&&1\\0&&1\\&0&&0\\0&&0\esm} \ar[dr] \ar[ur] &&
{\bsm1&&1\\&2&&1\\1&&1\\&0&&0\\0&&0\esm}\ar@{.}[d]\\
{\bsm0&&0\\&0&&0\\0&&0\\&1&&0\\0&&0\esm} \ar[r]\ar@{.}[d] & 
{\bsm0&&0\\&0&&0\\1&&1\\&2&&1\\0&&1\esm}\ar[r]\ar[dr]\ar[ur]&
{\bsm0&&0\\&0&&0\\1&&1\\&1&&1\\0&&1\esm} \ar[r]\ar[ddr] & 
{\bsm0&&0\\&0&&0\\1&&1\\&1&&1\\0&&0\esm}\ar[r]\ar[dr]\ar[ur]&
{\bsm0&&0\\&1&&0\\1&&1\\&1&&1\\0&&0\esm}\ar[r] & 
{\bsm0&&0\\&1&&0\\1&&2\\&1&&1\\0&&0\esm}\ar[r]\ar[dr]\ar[ur] &
{\bsm0&&0\\&0&&0\\0&&1\\&0&&0\\0&&0\esm} \ar[r] & 
{\bsm0&&0\\&1&&1\\1&&2\\&1&&0\\0&&0\esm}\ar[r]\ar[dr]\ar[ur]&
{\bsm0&&0\\&1&&1\\1&&1\\&1&&0\\0&&0\esm}\ar[r]\ar[ddr] & 
{\bsm0&&0\\&1&&1\\1&&1\\&0&&0\\0&&0\esm}\ar[r]\ar[dr]\ar[ur] &
{\bsm0&&1\\&1&&1\\1&&1\\&0&&0\\0&&0\esm} \ar[r] & 
{\bsm0&&1\\&2&&1\\1&&1\\&0&&0\\0&&0\esm}\ar[r]\ar[dr]\ar[ur]&
{\bsm0&&0\\&1&&0\\0&&0\\&0&&0\\0&&0\esm}\ar@{.}[d]\\
\ar@{.}[d]{\bsm0&&0\\&0&&0\\0&&1\\&1&&1\\0&&1\esm} \ar[ur] && 
{\bsm0&&0\\&0&&0\\1&&0\\&1&&0\\0&&0\esm} \ar[ur] &&
{\bsm0&&0\\&0&&0\\0&&1\\&0&&1\\0&&0\esm} \ar[ur] && 
{\bsm0&&0\\&1&&0\\1&&1\\&1&&0\\0&&0\esm} \ar[ur]&&
{\bsm0&&0\\&0&&1\\0&&1\\&0&&0\\0&&0\esm}\ar[ur] && 
{\bsm0&&0\\&1&&0\\1&&0\\&0&&0\\0&&0\esm}\ar[ur] &&
{\bsm0&&1\\&1&&1\\0&&1\\&0&&0\\0&&0\esm}\ar@{.}[d]\\
&&&{\bsm0&&0\\&1&&0\\1&&1\\&1&&1\\0&&1\esm} \ar[uur] &&&&&&
{\bsm0&&1\\&1&&1\\1&&1\\&1&&0\\0&&0\esm} \ar[uur] &&&\\ 
}
\]
\begin{center}
The Auslander-Reiten quiver of $\Lam_4$ 
\end{center}
\end{landscape}

\newpage
%%%%%%%%%%%%%%%%%%%%%%%%%%%%%%%%%%%%%%%%%%%%%%%%%%%%%%%%%%%%%%%%%%%%
\subsection{The indecomposable multisegments in $\MM(4)$}
\label{A4multi}
%%%%%%%%%%%%%%%%%%%%%%%%%%%%%%%%%%%%%%%%%%%%%%%%%%%%%%%%%%%%%%%%%%%%

\begin{align*}
&&\\
\msm_1 &= [1,1]&\msm_{21} &= [1,3]+[4,4]\\
\msm_2 &= [2,2]&\msm_{22} &= [1,2]+[3,3]+[4,4]\\
\msm_3 &= [3,3]&\msm_{23} &= [1,1]+[2,2]+[3,4]\\
\msm_4 &= [4,4]&\msm_{24} &= [1,1]+[2,4]\\
\msm_5 &= [1,2]&\msm_{25} &= [1,2]+[3,4]\\
\msm_6 &= [1,1]+[2,2]&\msm_{26} &= [1,1]+[2,3]+[4,4]\\
\msm_7 &= [2,3]&\msm_{27} &= [1,2]+[2,3]+[4,4]\\  
\msm_8 &= [2,2]+[3,3]& \msm_{28} &= [1,2]+[2,4]\\ 
\msm_9 &= [3,4]&\msm_{29} &= [1,3]+[3,4]\\
\msm_{10} &= [3,3]+[4,4]&\msm_{30} &= [1,1]+[2,3]+[3,4]\\
\msm_{11} &= [1,3]&\msm_{31} &= [1,2]+[2,3]+[3,3]+[4,4]\\
\msm_{12} &= [1,1]+[2,2]+[3,3]&\msm_{32} &= [1,2]+[2,4]+[3,3]\\
\msm_{13} &= [1,2]+[3,3]&\msm_{33} &= [1,3]+[2,2]+[3,4]\\
\msm_{14} &= [1,1]+[2,3]&\msm_{34} &= [1,1]+[2,2]+[2,3]+[3,4]\\
\msm_{15} &= [1,2]+[2,3]&\msm_{35} &= [1,1]+[1,3]+[2,2]+[3,4]\\
\msm_{16} &= [2,4]&\msm_{36} &= [1,2]+[2,4]+[3,3]+[4,4]\\
\msm_{17} &= [2,2]+[3,3]+[4,4]&\msm_{37} &= [1,4]\\
\msm_{18} &= [2,3]+[4,4]&\msm_{38} &= [1,3]+[2,4]\\
\msm_{19} &= [2,2]+[3,4]&\msm_{39} &= [1,2]+[2,3]+[3,4]\\
\msm_{20} &= [2,3]+[3,4]&\msm_{40} &= [1,1]+[2,2]+[3,3]+[4,4]\\
\end{align*}
\[
\xymatrix@-1.4pc{
&&&&\\
\ar@{.}[d]&&&{\msm_{37}}\ar[dr] &&&&&& {\msm_{40}}\ar[dr] &&&\ar@{.}[d]\\
{\msm_4}\ar[dr]\ar@{.}[dd] && 
{\msm_{11}}\ar[ur]\ar[dr] && 
{\msm_{16}}\ar[dr] && 
{\msm_1}\ar[dr] &&
{\msm_{17}}\ar[dr]\ar[ur] && 
{\msm_{12}}\ar[dr] &&
{\msm_4} \ar@{.}[dd]\\
&{\msm_{21}}\ar[ur]\ar[dr]&& 
{\msm_7}\ar[ur]\ar[dr] && 
{\msm_{24}}\ar[ur]\ar[dr] && 
{\msm_{22}}\ar[ur]\ar[dr] &&
{\msm_8}\ar[ur]\ar[dr]&& 
{\msm_{23}}\ar[ur]\ar[dr]&\\
{\msm_{35}} \ar[dr] \ar[ur]\ar@{.}[d] && 
{\msm_{18}} \ar[dr] \ar[ur] &&
{\msm_{14}} \ar[dr] \ar[ur] && 
{\msm_{36}} \ar[dr] \ar[ur] &&
{\msm_{13}} \ar[dr] \ar[ur] && 
{\msm_{19}} \ar[dr] \ar[ur] &&
{\msm_{35}}\ar@{.}[d]\\
{\msm_2} \ar[r]\ar@{.}[d] & 
{\msm_{34}}\ar[r]\ar[dr]\ar[ur]&
{\msm_{30}} \ar[r]\ar[ddr] & 
{\msm_{26}}\ar[r]\ar[dr]\ar[ur]&
{\msm_{27}}\ar[r] & 
{\msm_{31}}\ar[r]\ar[dr]\ar[ur] &
{\msm_3} \ar[r] & 
{\msm_{32}}\ar[r]\ar[dr]\ar[ur]&
{\msm_{28}}\ar[r]\ar[ddr] & 
{\msm_{25}}\ar[r]\ar[dr]\ar[ur] &
{\msm_{29}} \ar[r] & 
{\msm_{33}}\ar[r]\ar[dr]\ar[ur]&
{\msm_2}\ar@{.}[d]\\
\ar@{.}[d]{\msm_{20}} \ar[ur] && 
{\msm_6} \ar[ur] &&
{\msm_{10}} \ar[ur] && 
{\msm_{15}} \ar[ur]&&
{\msm_9}\ar[ur] && 
{\msm_5}\ar[ur] &&
{\msm_{20}}\ar@{.}[d]\\
&&&{\msm_{39}} \ar[uur] &&&&&&
{\msm_{38}} \ar[uur] &&&\\ 
}
\]
\begin{center}
The Auslander-Reiten quiver of $\Lam_4$ in terms of multisegments
\end{center}
\vspace{0.5cm}

\newpage
%%%%%%%%%%%%%%%%%%%%%%%%%%%%%%%%%%%%%%%%%%%%%%%%%%%%%%%%%%%%%%%%
\subsection{Extensions between indecomposable $\Lam_3$-modules}\label{extA3}
%%%%%%%%%%%%%%%%%%%%%%%%%%%%%%%%%%%%%%%%%%%%%%%%%%%%%%%%%%%%%%%%

\[
\xymatrix@-1.4pc{
\ar@{.}[d]&\cdot \ar[dr] && \bullet \ar[dr] &&\cdot \ar[dr]&\ar@{.}[d]&&&
\ar@{.}[d]&X_2 \ar[dr] && \bullet \ar[dr] &&\cdot \ar[dr]&\ar@{.}[d]\\ 
X_1 \ar[ur]\ar[dr]\ar@{.}[d] && \bullet \ar[ur]\ar[dr] && 
\bullet \ar[ur]\ar[dr] && \cdot\ar@{.}[d]&&&
\cdot \ar[ur]\ar[dr]\ar@{.}[d] && \cdot \ar[ur]\ar[dr] && 
\bullet \ar[ur]\ar[dr] && \cdot \ar@{.}[d]\\
&\cdot \ar[ur] && \bullet \ar[ur] && \cdot \ar[ur]&&&&
&\cdot \ar[ur] && \cdot \ar[ur] && \bullet \ar[ur]&
}
\]

%%%%%%%%%%%%%%%%%%%%%%%%%%%%%%%%%%%%%%%%%%%%%%%%%%%%%%%%%%%%%%%%
\subsection{Extensions between indecomposable $\Lam_4$-modules}\label{extA4}
%%%%%%%%%%%%%%%%%%%%%%%%%%%%%%%%%%%%%%%%%%%%%%%%%%%%%%%%%%%%%%%%

\[
\xymatrix@-1.4pc{
X_1 \ar[dr]\ar@{.}[dd] && 
{\bullet} \ar[dr] && {\cdot} \ar[dr] && 
{\cdot} \ar[dr] && {\cdot} \ar[dr] && 
{\bullet} \ar[dr] && {\cdot}\ar@{.}[dd] \\
&{\cdot} \ar[dr]\ar[ur] &&
{\bullet} \ar[dr]\ar[ur] && {\cdot} \ar[dr]\ar[ur] &&
{\cdot} \ar[dr]\ar[ur] && {\bullet} \ar[dr]\ar[ur] &&
{\cdot} \ar[dr]\ar[ur] &&\\
{\cdot} \ar[dr]\ar[ur]\ar@{.}[d] &&
{\cdot} \ar[dr]\ar[ur] && {\bullet} \ar[dr]\ar[ur] &&
{\cdot} \ar[dr]\ar[ur] && {\bullet} \ar[dr]\ar[ur] &&
{\cdot} \ar[dr]\ar[ur] && {\cdot}\ar@{.}[d] \\
{\cdot} \ar[r]\ar@{.}[d] & {\cdot} \ar[r]\ar[dr]\ar[ur] & 
{\cdot} \ar[r] & {\cdot} \ar[r]\ar[dr]\ar[ur] & 
{\cdot} \ar[r] & {\bullet} \ar[r]\ar[dr]\ar[ur] & 
{\bullet} \ar[r] & {\bullet} \ar[r]\ar[dr]\ar[ur] & 
{\cdot} \ar[r] & {\cdot} \ar[r]\ar[dr]\ar[ur] & 
{\cdot} \ar[r] & {\cdot} \ar[r]\ar[dr]\ar[ur] & {\cdot}\ar@{.}[d]\\
{\cdot} \ar[ur] && 
{\cdot}\ar[ur] && {\cdot} \ar[ur] && 
{\bullet}\ar[ur] && {\cdot} \ar[ur] &&
{\cdot}\ar[ur] && {\cdot}
}
\]
\[
\xymatrix@-1.4pc{
{\cdot} \ar[dr]\ar@{.}[dd] && 
{\cdot} \ar[dr] && {\bullet} \ar[dr] && 
{\cdot} \ar[dr] && {\bullet} \ar[dr] && 
{\cdot} \ar[dr] && {\cdot}\ar@{.}[dd] \\
&{\cdot} \ar[dr]\ar[ur] &&
{\bullet} \ar[dr]\ar[ur] && {\bullet} \ar[dr]\ar[ur] &&
{\bullet} \ar[dr]\ar[ur] && {\bullet} \ar[dr]\ar[ur] &&
{\cdot} \ar[dr]\ar[ur] &&\\
X_2 \ar[dr]\ar@{.}[d]\ar[ur] &&
{\bullet} \ar[dr]\ar[ur] && {\bullet} \ar[dr]\ar[ur] &&
{\bullet} \ar[dr]\ar[ur] && {\bullet} \ar[dr]\ar[ur] &&
{\bullet} \ar[dr]\ar[ur] && {\cdot}\ar@{.}[d] \\
{\cdot} \ar[r]\ar@{.}[d] & {\cdot} \ar[r]\ar[dr]\ar[ur] & 
{\cdot} \ar[r] & {\bullet} \ar[r]\ar[dr]\ar[ur] & 
{\bullet} \ar[r] & {\bullet} \ar[r]\ar[dr]\ar[ur] & 
{\bullet} \ar[r] & {\bullet} \ar[r]\ar[dr]\ar[ur] & 
{\bullet} \ar[r] & {\bullet} \ar[r]\ar[dr]\ar[ur] & 
{\cdot} \ar[r] & {\cdot} \ar[r]\ar[dr]\ar[ur] & {\cdot}\ar@{.}[d]\\
{\cdot} \ar[ur] && 
{\cdot}\ar[ur] && {\bullet} \ar[ur] && 
{\bullet}\ar[ur] && {\bullet} \ar[ur] &&
{\cdot}\ar[ur] && {\cdot}
}
\]
\[
\xymatrix@-1.4pc{
{\cdot} \ar[dr]\ar@{.}[dd] && 
{\cdot} \ar[dr] && {\cdot} \ar[dr] && 
{\bullet} \ar[dr] && {\cdot} \ar[dr] && 
{\cdot} \ar[dr] && {\cdot}\ar@{.}[dd] \\
&{\cdot} \ar[dr]\ar[ur] &&
{\cdot} \ar[dr]\ar[ur] && {\bullet} \ar[dr]\ar[ur] &&
{\bullet} \ar[dr]\ar[ur] && {\cdot} \ar[dr]\ar[ur] &&
{\cdot} \ar[dr]\ar[ur] &&\\
{\cdot} \ar[dr]\ar[ur]\ar@{.}[d] &&
{\cdot} \ar[dr]\ar[ur] && {\bullet} \ar[dr]\ar[ur] &&
{\bullet} \ar[dr]\ar[ur] && {\bullet} \ar[dr]\ar[ur] &&
{\cdot} \ar[dr]\ar[ur] && {\cdot}\ar@{.}[d] \\
X_3 \ar[r]\ar@{.}[d] & {\cdot} \ar[r]\ar[dr]\ar[ur] & 
{\bullet} \ar[r] & {\bullet} \ar[r]\ar[dr]\ar[ur] & 
{\cdot} \ar[r] & {\bullet} \ar[r]\ar[dr]\ar[ur] & 
{\bullet} \ar[r] & {\bullet} \ar[r]\ar[dr]\ar[ur] & 
{\cdot} \ar[r] & {\bullet} \ar[r]\ar[dr]\ar[ur] & 
{\bullet} \ar[r] & {\cdot} \ar[r]\ar[dr]\ar[ur] & {\cdot}\ar@{.}[d]\\
{\cdot} \ar[ur] && 
{\cdot}\ar[ur] && {\bullet} \ar[ur] && 
{\cdot}\ar[ur] && {\bullet} \ar[ur] &&
{\cdot}\ar[ur] && {\cdot}
}
\]
\[
\xymatrix@-1.4pc{
{\cdot} \ar[dr]\ar@{.}[dd] && 
{\cdot} \ar[dr] && {\cdot} \ar[dr] && 
{\bullet} \ar[dr] && {\cdot} \ar[dr] && 
{\cdot} \ar[dr] && {\cdot}\ar@{.}[dd] \\
&{\cdot} \ar[dr]\ar[ur] &&
{\cdot} \ar[dr]\ar[ur] && {\bullet} \ar[dr]\ar[ur] &&
{\bullet} \ar[dr]\ar[ur] && {\cdot} \ar[dr]\ar[ur] &&
{\cdot} \ar[dr]\ar[ur] &&\\
{\cdot} \ar[dr]\ar[ur]\ar@{.}[d] &&
{\cdot} \ar[dr]\ar[ur] && {\bullet} \ar[dr]\ar[ur] &&
{\bullet} \ar[dr]\ar[ur] && {\bullet} \ar[dr]\ar[ur] &&
{\cdot} \ar[dr]\ar[ur] && {\cdot}\ar@{.}[d] \\
{\cdot} \ar[r]\ar@{.}[d] & {\cdot} \ar[r]\ar[dr]\ar[ur] & 
{\cdot} \ar[r] & {\bullet} \ar[r]\ar[dr]\ar[ur] & 
{\bullet} \ar[r] & {\bullet} \ar[r]\ar[dr]\ar[ur] & 
{\bullet} \ar[r] & {\bullet} \ar[r]\ar[dr]\ar[ur] & 
{\bullet} \ar[r] & {\bullet} \ar[r]\ar[dr]\ar[ur] & 
{\cdot} \ar[r] & {\cdot} \ar[r]\ar[dr]\ar[ur] & {\cdot}\ar@{.}[d]\\
X_4 \ar[ur] && 
{\bullet}\ar[ur] && {\cdot} \ar[ur] && 
{\bullet}\ar[ur] && {\cdot} \ar[ur] &&
{\bullet}\ar[ur] && {\cdot}
}
\]
\[
\xymatrix@-1.4pc{
{\cdot} \ar[dr]\ar@{.}[dd] && 
{\cdot} \ar[dr] && {\bullet} \ar[dr] && 
{\cdot} \ar[dr] && {\cdot} \ar[dr] && 
{\bullet} \ar[dr] && {\cdot}\ar@{.}[dd] \\
&X_5 \ar[dr]\ar[ur] &&
{\bullet} \ar[dr]\ar[ur] && {\bullet} \ar[dr]\ar[ur] &&
{\cdot} \ar[dr]\ar[ur] && {\bullet} \ar[dr]\ar[ur] &&
{\bullet} \ar[dr]\ar[ur] &&\\
{\cdot} \ar[dr]\ar[ur]\ar@{.}[d] &&
{\cdot} \ar[dr]\ar[ur] && {\bullet} \ar[dr]\ar[ur] &&
{\bullet} \ar[dr]\ar[ur] && {\bullet} \ar[dr]\ar[ur] &&
{\bullet} \ar[dr]\ar[ur] && {\cdot}\ar@{.}[d] \\
{\cdot} \ar[r]\ar@{.}[d] & {\cdot} \ar[r]\ar[dr]\ar[ur] & 
{\cdot} \ar[r] & {\cdot} \ar[r]\ar[dr]\ar[ur] & 
{\cdot} \ar[r] & {\bullet} \ar[r]\ar[dr]\ar[ur] & 
{\bullet} \ar[r] & {\bullet} \ar[r]\ar[dr]\ar[ur] & 
{\bullet} \ar[r] & {\bullet} \ar[r]\ar[dr]\ar[ur] & 
{\cdot} \ar[r] & {\cdot} \ar[r]\ar[dr]\ar[ur] & {\cdot}\ar@{.}[d]\\
{\cdot} \ar[ur] && 
{\cdot}\ar[ur] && {\cdot} \ar[ur] && 
{\bullet}\ar[ur] && {\bullet} \ar[ur] &&
{\cdot}\ar[ur] && {\cdot}
}
\]
\[
\xymatrix@-1.4pc{
{\cdot} \ar[dr]\ar@{.}[dd] && 
{\cdot} \ar[dr] && {\cdot} \ar[dr] && 
{\bullet} \ar[dr] && {\bullet} \ar[dr] && 
{\cdot} \ar[dr] && {\cdot}\ar@{.}[dd] \\
&{\cdot} \ar[dr]\ar[ur] &&
{\cdot} \ar[dr]\ar[ur] && {\bullet} \ar[dr]\ar[ur] &&
{\bullet} \ar[dr]\ar[ur] && {\bullet} \ar[dr]\ar[ur] &&
{\cdot} \ar[dr]\ar[ur] &&\\
{\cdot} \ar[dr]\ar[ur]\ar@{.}[d] &&
{\cdot} \ar[dr]\ar[ur] && {\bullet} \ar[dr]\ar[ur] &&
{\bullet} \ar[dr]\ar[ur] && {\bullet} \ar[dr]\ar[ur] &&
{\bullet} \ar[dr]\ar[ur] && {\cdot}\ar@{.}[d] \\
{\cdot} \ar[r]\ar@{.}[d] & X_6 \ar[r]\ar[dr]\ar[ur] & 
{\cdot} \ar[r] & {\bullet} \ar[r]\ar[dr]\ar[ur] & 
{\bullet} \ar[r] & {\bullet} \ar[r]\ar[dr]\ar[ur] & 
{\bullet} \ar[r] & {\bullet} \ar[r]\ar[dr]\ar[ur] & 
{\bullet} \ar[r] & {\bullet} \ar[r]\ar[dr]\ar[ur] & 
{\bullet} \ar[r] & {\bullet} \ar[r]\ar[dr]\ar[ur] & {\cdot}\ar@{.}[d]\\
{\cdot} \ar[ur] && 
{\cdot}\ar[ur] && {\bullet} \ar[ur] && 
{\bullet}\ar[ur] && {\bullet} \ar[ur] &&
{\bullet}\ar[ur] && {\cdot}
}
\]

\newpage
%%%%%%%%%%%%%%%%%%%%%%%%%%%%%%%%%%%%%%%%%%%%%%%%%%%%%%%%%%%%%%%%
\subsection{The graph ${\mathcal G}^\circ_4$}\label{graphA4}
%%%%%%%%%%%%%%%%%%%%%%%%%%%%%%%%%%%%%%%%%%%%%%%%%%%%%%%%%%%%%%%%

\begin{small}
\[
\begin{tabular}{llllllllll}
&&&&&&&&&\\
(1,3) & (1,4) & (1,5) & (1,6) & (1,9) & (1,10) & (1,11) & (1,12)&(1,13) & (1,14)\\ 
(1,15) & (1,21) & (1,22) & (1,23) & (1,24) & (1,25)&(1,26) & (1,27) & (1,28) & (1,29)\\ 
(1,30) & (1,31) & (1,32) & (1,35)&(1,36) & (2,4) & (2,5) & (2,6) & (2,7) & (2,8)\\ 
(2,11) & (2,12)&(2,15) & (2,16) & (2,17) & (2,18) & (2,19) & (2,20) & (2,21) & (2,23)\\
(2,27) & (2,28) & (2,33) & (2,34) & (2,35) & (3,7) & (3,8) & (3,9)&(3,10) & (3,11)\\ 
(3,12) & (3,13) & (3,14) & (3,15) & (3,16) & (3,17)&(3,20) & (3,22) & (3,24) & (3,29)\\ 
(3,30) & (3,31) & (3,32) & (3,36)&(4,5) & (4,6) & (4,9) & (4,10) & (4,16) & (4,17)\\ 
(4,18) & (4,19)&(4,20) & (4,21) & (4,22) & (4,23) & (4,24) & (4,25) & (4,26) & (4,27)\\
(4,28) & (4,29) & (4,30) & (4,33) & (4,34) & (4,35) & (4,36) & (5,6)&(5,8) & (5,11)\\ 
(5,12) & (5,13) & (5,15) &(5,17) & (5,19) &(5,21) &(5,22) &(5,23) &(5,25) &(5,27)\\
(5,28) & (5,29) &(5,33) &(5,35) &(6,7) & (6,11) &(6,12) &(6,14) &(6,15)&(6,16)\\ 
(6,18) &(6,21) &(6,23) &(6,24) &(6,26) &(6,27) &(6,28) & (6,30)&(6,34) & (6,35)\\
(7,8) &(7,10) &(7,11) &(7,12) &(7,14) &(7,15) & (7,16)&(7,17)&(7,18) &(7,20)\\
(7,26) &(7,27) &(7,30) &(7,31) &(7,34) &(8,9)&(8,11) &(8,12) &(8,13) &(8,15)\\
(8,16) &(8,17) &(8,19) &(8,20) &(8,25)&(8,28) &(8,29) &(8,32) &(8,33) &(9,10)\\
(9,12) & (9,13)&(9,16) &(9,17)&(9,19) &(9,20) &(9,22) &(9,23) &(9,24) &(9,25)\\
(9,28) &(9,29)&(9,30)&(9,32) &(9,36) &(10,11) &(10,14) &(10,16) &(10,17) &(10,18) \\
(10,20)&(10,21) &(10,22)&(10,24) &(10,26) &(10,27) &(10,29) &(10,30) &(10,31)&(10,36)\\
(11,12) &(11,13)&(11,14) &(11,15) &(11,17) &(11,18) &(11,20) &(11,21) &(11,22) &(11,26)\\
(11,27) &(11,29) &(11,30) &(11,31) &(11,33)&(11,34) &(11,35) &(12,13) &(12,14) &(12,15)\\
(12,16) &(12,19)&(12,20)&(12,23) &(12,24) &(12,25) &(12,28) &(12,29) &(12,30) &(12,32) \\
(12,33)&(12,34)&(12,35) &(13,15) &(13,17) &(13,22) &(13,25) &(13,28) &(13,29)&(13,32)\\
(14,15) &(14,16) &(14,24) &(14,26) &(14,27) &(14,30) &(14,31)&(15,16) &(15,17) &(15,22)\\
(15,24) &(15,27) &(15,28) &(15,31) &(15,32)&(15,36) &(16,17) &(16,18) &(16,19) &(16,20)\\
(16,23) &(16,24)&(16,26)&(16,27) &(16,28) &(16,30) &(16,31) &(16,32) &(16,34) &(16,36) \\
(17,18)&(17,19) &(17,20) &(17,21) &(17,22) &(17,25) &(17,27) &(17,28) &(17,29)&(17,31)\\
(17,32) &(17,33) &(17,36) &(18,20) &(18,21) &(18,26) &(18,27)&(18,30) &(18,34) &(19,20) \\
(19,23) &(19,25) &(19,28) &(19,29) &(19,33)&(20,21) &(20,23) &(20,29) &(20,30) &(20,33)\\
(20,34) &(20,35) &(21,22)&(21,26) &(21,27) &(21,29) &(21,30) &(21,33) &(21,34) &(21,35) \\
(22,25)&(22,27) &(22,28) &(22,29) &(22,31) &(22,32) &(22,36) &(23,24) &(23,25)&(23,28) \\
(23,29) &(23,30) &(23,33) &(23,34) &(23,35) &(24,26) &(24,27)&(24,28) &(24,30) &(24,31) \\
(24,32) &(24,36) &(25,28) &(25,29) &(26,27)&(26,30) &(27,28) &(27,31) &(27,36) &(28,32) \\
(28,36) &(29,30)&(29,33)&(29,35) &(30,34) & (30,35) &(31,36) &(32,36) &(33,35) &(34,35)
\end{tabular}
\]
\end{small}

\newpage
%%%%%%%%%%%%%%%%%%%%%%%%%%%%%%%%%%%%%%%%%%%%%%
\subsection{The non-homogeneous tubes in $\Tub_0$}\label{Tnull}
%%%%%%%%%%%%%%%%%%%%%%%%%%%%%%%%%%%%%%%%%%%%%%

\[
\xymatrix@-1.5pc{
&&&&&&&&&&&&&&\\
&&&&&&&&&&&&&&\\
&&&&&&&&{}\save[]*{\bsm 1&&2&&1\\&3&&3\\1&&2&&1\esm}\restore\phZ\ar[rd]\ar@{.}[u]&&
 {}\save[]*{\bsm 1&&2&&1\\&3&&3\\1&&2&&1\esm}\restore\phZ\ar[rd]&&
 {}\save[]*{\bsm 1&&2&&1\\&3&&3\\1&&2&&1\esm}\restore\phZ\ar[rd]&&
 {}\save[]*{\bsm 1&&2&&1\\&3&&3\\1&&2&&1\esm}\restore\phZ\ar@{.}[u]\\
&{}\save[]*{\bsm 1&&2&&1\\&3&&3\\1&&2&&1\esm}\restore\phZ\ar[rd]&&
 {}\save[]*{\bsm 1&&2&&1\\&3&&3\\1&&2&&1\esm}\restore\phZ\ar[rd]& &&&&
&{}\save[]*{\bsm 1&&1&&1\\&2&&2\\1&&1&&1\esm}\restore\phZ\ar[rd]\ar[ru]&&
 {}\save[]*{\bsm 0&&2&&0\\&2&&2\\1&&1&&1\esm}\restore\phZ\ar[rd]\ar[ru]&&
 {}\save[]*{\bsm 1&&1&&1\\&2&&2\\0&&2&&0\esm}\restore\phZ\ar[rd]\ar[ru]&\\
 {}\save[]*{\bsm 1&&1&&0\\&2&&1\\1&&1&&0\esm}\restore\phZ\ar[ru]\ar@{.}[uu]&&
 {}\save[]*{\bsm 0&&1&&1\\&1&&2\\0&&1&&1\esm}\restore\phZ\ar[ru]&&
 {}\save[]*{\bsm 1&&1&&0\\&2&&1\\1&&1&&0\esm}\restore\phZ\ar@{.}[uu] &&&&
 {}\save[]*{\bsm 1&&0&&1\\&1&&1\\0&&1&&0\esm}\restore\phZ\ar[ru]\ar@{.}[uu]&&
 {}\save[]*{\bsm 0&&1&&0\\&1&&1\\1&&0&&1\esm}\restore\phZ\ar[ru]&&
 {}\save[]*{\bsm 0&&1&&0\\&1&&1\\0&&1&&0\esm}\restore\phZ\ar[ru]&&
 {}\save[]*{\bsm 1&&0&&1\\&1&&1\\0&&1&&0\esm}\restore\phZ\ar@{.}[uu]
}
\]
\[
\xymatrix@-1.5pc{
&&&&&&&&&&&&\\
&&&&&&&&&&&&\\
&&&&&&&&&&&&\\
&{}\save[]*{\bsm&0&&0\\1&&2&&1\\&3&&3\\1&&2&&1\\&0&&0\esm}\restore\phX\ar[rd]&&
 {}\save[]*{\bsm&0&&1\\1&&3&&2\\&4&&4\\2&&3&&1\\&1&&0\esm}\restore\phX\ar[rd]&&
 {}\save[]*{\bsm&0&&0\\1&&2&&1\\&3&&3\\1&&2&&1\\&0&&0\esm}\restore\phX\ar[rd]&&
 {}\save[]*{\bsm&0&&0\\1&&2&&1\\&3&&3\\1&&2&&1\\&0&&0\esm}\restore\phX\ar[rd]&&
 {}\save[]*{\bsm&1&&0\\2&&3&&1\\&4&&4\\1&&3&&2\\&0&&1\esm}\restore\phX\ar[rd]&&
 {}\save[]*{\bsm&0&&0\\1&&2&&1\\&3&&3\\1&&2&&1\\&0&&0\esm}\restore\phX\ar[rd]&\\
 {}\save[]*{\bsm&0&&0\\1&&2&&1\\&3&&2\\1&&2&&1\\&0&&0\esm}\restore\phX\ar[rd]\ar[ru]\ar@{.}[dd]\ar@{.}[uu]&&
 {}\save[]*{\bsm&0&&1\\1&&2&&1\\&3&&3\\1&&2&&1\\&0&&0\esm}\restore\phX\ar[rd]\ar[ru]&&
 {}\save[]*{\bsm&0&&0\\1&&2&&1\\&3&&3\\1&&2&&1\\&1&&0\esm}\restore\phX\ar[rd]\ar[ru]&&
 {}\save[]*{\bsm&0&&0\\1&&2&&1\\&2&&3\\1&&2&&1\\&0&&0\esm}\restore\phX\ar[rd]\ar[ru]&&
 {}\save[]*{\bsm&1&&0\\1&&2&&1\\&3&&3\\1&&2&&1\\&0&&0\esm}\restore\phX\ar[rd]\ar[ru]&&
 {}\save[]*{\bsm&0&&0\\1&&2&&1\\&3&&3\\1&&2&&1\\&0&&1\esm}\restore\phX\ar[rd]\ar[ru]&&
 {}\save[]*{\bsm&0&&0\\1&&2&&1\\&3&&2\\1&&2&&1\\&0&&0\esm}\restore\phX\ar@{.}[dd]\ar@{.}[uu]\\
&{}\save[]*{\bsm&0&&1\\1&&2&&1\\&3&&2\\1&&2&&1\\&0&&0\esm}\restore\phX\ar[rd]\ar[ru]&&
 {}\save[]*{\bsm&0&&0\\1&&1&&0\\&2&&2\\0&&1&&1\\&0&&0\esm}\restore\phX\ar[rd]\ar[ru]&&
 {}\save[]*{\bsm&0&&0\\1&&2&&1\\&2&&3\\1&&2&&1\\&1&&0\esm}\restore\phX\ar[rd]\ar[ru]&&
 {}\save[]*{\bsm&1&&0\\1&&2&&1\\&2&&3\\1&&2&&1\\&0&&0\esm}\restore\phX\ar[rd]\ar[ru]&&
 {}\save[]*{\bsm&0&&0\\0&&1&&1\\&2&&2\\1&&1&&0\\&0&&0\esm}\restore\phX\ar[rd]\ar[ru]&&
 {}\save[]*{\bsm&0&&0\\1&&2&&1\\&3&&2\\1&&2&&1\\&0&&1\esm}\restore\phX\ar[rd]\ar[ru]&\\
 {}\save[]*{\bsm&0&&1\\1&&2&&1\\&3&&2\\1&&2&&1\\&0&&1\esm}\restore\phX\ar[ru]\ar[rd]&&
 {}\save[]*{\bsm&0&&0\\1&&1&&0\\&2&&1\\0&&1&&1\\&0&&0\esm}\restore\phX\ar[ru]\ar[rd]&&
 {}\save[]*{\bsm&0&&0\\1&&1&&0\\&1&&2\\0&&1&&1\\&0&&0\esm}\restore\phX\ar[ru]\ar[rd]&&
 {}\save[]*{\bsm&1&&0\\1&&2&&1\\&2&&3\\1&&2&&1\\&1&&0\esm}\restore\phX\ar[ru]\ar[rd]&&
 {}\save[]*{\bsm&0&&0\\0&&1&&1\\&1&&2\\1&&1&&0\\&0&&0\esm}\restore\phX\ar[ru]\ar[rd]&&
 {}\save[]*{\bsm&0&&0\\0&&1&&1\\&2&&1\\1&&1&&0\\&0&&0\esm}\restore\phX\ar[ru]\ar[rd]&&
 {}\save[]*{\bsm&0&&1\\1&&2&&1\\&3&&2\\1&&2&&1\\&0&&1\esm}\restore\phX\\
&{}\save[]*{\bsm&0&&0\\1&&1&&0\\&2&&1\\0&&1&&1\\&0&&1\esm}\restore\phX\ar[rd]\ar[ru]&&
 {}\save[]*{\bsm&0&&0\\1&&1&&0\\&1&&1\\0&&1&&1\\&0&&0\esm}\restore\phX\ar[rd]\ar[ru]&&
 {}\save[]*{\bsm&1&&0\\1&&1&&0\\&1&&2\\0&&1&&1\\&0&&0\esm}\restore\phX\ar[rd]\ar[ru]&&
 {}\save[]*{\bsm&0&&0\\0&&1&&1\\&1&&2\\1&&1&&0\\&1&&0\esm}\restore\phX\ar[rd]\ar[ru]&&
 {}\save[]*{\bsm&0&&0\\0&&1&&1\\&1&&1\\1&&1&&0\\&0&&0\esm}\restore\phX\ar[rd]\ar[ru]&&
 {}\save[]*{\bsm&0&&1\\0&&1&&1\\&2&&1\\1&&1&&0\\&0&&0\esm}\restore\phX\ar[rd]\ar[ru]&\\
 {}\save[]*{\bsm&0&&0\\0&&0&&0\\&1&&0\\0&&0&&0\\&0&&0\esm}\restore\phX\ar[ru]\ar@{.}[uu]&&
 {}\save[]*{\bsm&0&&0\\1&&1&&0\\&1&&1\\0&&1&&1\\&0&&1\esm}\restore\phX\ar[ru]\ar[rd]&&
 {}\save[]*{\bsm&1&&0\\1&&1&&0\\&1&&1\\0&&1&&1\\&0&&0\esm}\restore\phX\ar[ru]&&
 {}\save[]*{\bsm&0&&0\\0&&0&&0\\&0&&1\\0&&0&&0\\&0&&0\esm}\restore\phX\ar[ru]&&
 {}\save[]*{\bsm&0&&0\\0&&1&&1\\&1&&1\\1&&1&&0\\&1&&0\esm}\restore\phX\ar[ru]\ar[rd]&&
 {}\save[]*{\bsm&0&&1\\0&&1&&1\\&1&&1\\1&&1&&0\\&0&&0\esm}\restore\phX\ar[ru]&&
 {}\save[]*{\bsm&0&&0\\0&&0&&0\\&1&&0\\0&&0&&0\\&0&&0\esm}\restore\phX\ar@{.}[uu]\\
&&&{}\save[]*{\bsm&1&&0\\1&&1&&0\\&1&&1\\0&&1&&1\\&0&&1\esm}\restore\phX\ar[ru]&&&&&&
 {}\save[]*{\bsm&0&&1\\0&&1&&1\\&1&&1\\1&&1&&0\\&1&&0\esm}\restore\phX\ar[ru]\\
}
\]

\newpage
%%%%%%%%%%%%%%%%%%%%%%%%%%%%%%%%%%%%%%%%%%%%%%%%%%%
\subsection{The non-homogeneous tubes in $\Tub_\infty$}\label{Tinfty}
%%%%%%%%%%%%%%%%%%%%%%%%%%%%%%%%%%%%%%%%%%%%%%%%%%%

\[
\xymatrix@-1.5pc{
&&&&&&&&&&&&&&&\\
&&&&&&&&&&&&&&&\\
&&&&&&&&{}\save[]*{\bsm&&1\\&2&&2\\2&&3&&2\\&2&&2\\&&1\esm}\restore\phZ\ar[dr]\ar@{.}[uu]&&
 {}\save[]*{\bsm&&0\\&1&&1\\1&&2&&1\\&1&&1\\&&0\esm}\restore\phZ\ar[rd]&&
 {}\save[]*{\bsm&&0\\&1&&1\\1&&2&&1\\&1&&1\\&&0\esm}\restore\phZ\ar[rd]&&
 {}\save[]*{\bsm&&1\\&2&&2\\2&&3&&2\\&2&&2\\&&1\esm}\restore\phZ\ar@{.}[uu]&&&\\
&{}\save[]*{\bsm&1&&1\\1&&2&&1\\&1&&1\esm}\restore\phZ\ar[rd]&&
 {}\save[]*{\bsm&1&&1\\1&&2&&1\\&1&&1\esm}\restore\phZ\ar[rd]& &&&&
&{}\save[]*{\bsm&&0\\&1&&1\\1&&2&&1\\&1&&1\\&&1\esm}\restore\phZ\ar[rd]\ar[ur]&&
 {}\save[]*{\bsm&&0\\&1&&1\\1&&1&&1\\&1&&1\\&&0\esm}\restore\phZ\ar[rd]\ar[ur]&&
 {}\save[]*{\bsm&&1\\&1&&1\\1&&2&&1\\&1&&1\\&&0\esm}\restore\phZ\ar[rd]\ar[ur]&\\
 {}\save[]*{\bsm&1&&0\\1&&1&&0\\&1&&0\esm}\restore\phZ\ar[ru]\ar@{.}[uu]&&
 {}\save[]*{\bsm&0&&1\\0&&1&&1\\&0&&1\esm}\restore\phZ\ar[ru]&&
 {}\save[]*{\bsm&1&&0\\1&&1&&0\\&1&&0\esm}\restore\phZ\ar@{.}[uu] &&&&
 {}\save[]*{\bsm&&0\\&0&&0\\0&&1&&0\\&0&&0\\&&0\esm}\restore\phZ\ar[ru]\ar@{.}[uu]&&
 {}\save[]*{\bsm&&0\\&1&&1\\1&&1&&1\\&1&&1\\&&1\esm}\restore\phZ\ar[ru]\ar[rd]&&
 {}\save[]*{\bsm&&1\\&1&&1\\1&&1&&1\\&1&&1\\&&0\esm}\restore\phZ\ar[ru]&&
 {}\save[]*{\bsm&&0\\&0&&0\\0&&1&&0\\&0&&0\\&&0\esm}\restore\phZ\ar@{.}[uu]\\
&&&&&&&& 
&&&{}\save[]*{\bsm&&1\\&1&&1\\1&&1&&1\\&1&&1\\&&1\esm}\restore\phZ\ar[ru]
}
\]
\[
\xymatrix@-1.5pc{
&&&&&&&&&&&&\\
&&&&&&&&&&&&\\
&&&&&&&&&&&&\\
&{}\save[]*{\bsm0&&&&0\\&1&&1\\1&&2&&1\\&1&&1\\0&&&&0\esm}\restore\phX\ar[rd]&&
 {}\save[]*{\bsm0&&&&1\\&1&&2\\1&&3&&1\\&2&&1\\1&&&&0\esm}\restore\phX\ar[rd]&&
 {}\save[]*{\bsm0&&&&0\\&1&&1\\1&&2&&1\\&1&&1\\0&&&&0\esm}\restore\phX\ar[rd]&&
 {}\save[]*{\bsm0&&&&0\\&1&&1\\1&&2&&1\\&1&&1\\0&&&&0\esm}\restore\phX\ar[rd]&&
 {}\save[]*{\bsm1&&&&0\\&2&&1\\1&&3&&1\\&1&&2\\0&&&&1\esm}\restore\phX\ar[rd]&&
 {}\save[]*{\bsm0&&&&0\\&1&&1\\1&&2&&1\\&1&&1\\0&&&&0\esm}\restore\phX\ar[rd]&\\
 {}\save[]*{\bsm0&&&&0\\&1&&1\\1&&2&&0\\&1&&1\\0&&&&0\esm}\restore\phX\ar[rd]\ar[ur]\ar@{.}[uu]\ar@{.}[dd]&&
 {}\save[]*{\bsm0&&&&1\\&1&&1\\1&&2&&1\\&1&&1\\0&&&&0\esm}\restore\phX\ar[rd]\ar[ur]&&
 {}\save[]*{\bsm0&&&&0\\&1&&1\\1&&2&&1\\&1&&1\\1&&&&0\esm}\restore\phX\ar[rd]\ar[ur]&&
 {}\save[]*{\bsm0&&&&0\\&1&&1\\0&&2&&1\\&1&&1\\0&&&&0\esm}\restore\phX\ar[rd]\ar[ur]&&
 {}\save[]*{\bsm1&&&&0\\&1&&1\\1&&2&&1\\&1&&1\\0&&&&0\esm}\restore\phX\ar[rd]\ar[ur]&&
 {}\save[]*{\bsm0&&&&0\\&1&&1\\1&&2&&1\\&1&&1\\0&&&&1\esm}\restore\phX\ar[rd]\ar[ur]&&
 {}\save[]*{\bsm0&&&&0\\&1&&1\\1&&2&&0\\&1&&1\\0&&&&0\esm}\restore\phX\ar@{.}[uu]\ar@{.}[dd]\\
&{}\save[]*{\bsm0&&&&1\\&1&&1\\1&&2&&0\\&1&&1\\0&&&&0\esm}\restore\phX\ar[rd]\ar[ru]&&
 {}\save[]*{\bsm0&&&&0\\&1&&0\\1&&1&&1\\&0&&1\\0&&&&0\esm}\restore\phX\ar[rd]\ar[ru]&&
 {}\save[]*{\bsm0&&&&0\\&1&&1\\0&&2&&1\\&1&&1\\1&&&&0\esm}\restore\phX\ar[rd]\ar[ru]&&
 {}\save[]*{\bsm1&&&&0\\&1&&1\\0&&2&&1\\&1&&1\\0&&&&0\esm}\restore\phX\ar[rd]\ar[ru]&&
 {}\save[]*{\bsm0&&&&0\\&0&&1\\1&&1&&1\\&1&&0\\0&&&&0\esm}\restore\phX\ar[rd]\ar[ru]&&
 {}\save[]*{\bsm0&&&&0\\&1&&1\\1&&2&&0\\&1&&1\\0&&&&1\esm}\restore\phX\ar[rd]\ar[ru]&\\
 {}\save[]*{\bsm0&&&&1\\&1&&1\\1&&2&&0\\&1&&1\\0&&&&1\esm}\restore\phX\ar[rd]\ar[ru]&&
 {}\save[]*{\bsm0&&&&0\\&1&&0\\1&&1&&0\\&0&&1\\0&&&&0\esm}\restore\phX\ar[rd]\ar[ru]&&
 {}\save[]*{\bsm0&&&&0\\&1&&0\\0&&1&&1\\&0&&1\\0&&&&0\esm}\restore\phX\ar[rd]\ar[ru]&&
 {}\save[]*{\bsm1&&&&0\\&1&&1\\0&&2&&1\\&1&&1\\1&&&&0\esm}\restore\phX\ar[rd]\ar[ru]&&
 {}\save[]*{\bsm0&&&&0\\&0&&1\\0&&1&&1\\&1&&0\\0&&&&0\esm}\restore\phX\ar[rd]\ar[ru]&&
 {}\save[]*{\bsm0&&&&0\\&0&&1\\1&&1&&0\\&1&&0\\0&&&&0\esm}\restore\phX\ar[rd]\ar[ru]&&
 {}\save[]*{\bsm0&&&&1\\&1&&1\\1&&2&&0\\&1&&1\\0&&&&1\esm}\restore\phX\\
&{}\save[]*{\bsm0&&&&0\\&1&&0\\1&&1&&0\\&0&&1\\0&&&&1\esm}\restore\phX\ar[rd]\ar[ru]&&
 {}\save[]*{\bsm0&&&&0\\&1&&0\\0&&1&&0\\&0&&1\\0&&&&0\esm}\restore\phX\ar[rd]\ar[ru]&&
 {}\save[]*{\bsm1&&&&0\\&1&&0\\0&&1&&1\\&0&&1\\0&&&&0\esm}\restore\phX\ar[rd]\ar[ru]&&
 {}\save[]*{\bsm0&&&&0\\&0&&1\\0&&1&&1\\&1&&0\\1&&&&0\esm}\restore\phX\ar[rd]\ar[ru]&&
 {}\save[]*{\bsm0&&&&0\\&0&&1\\0&&1&&0\\&1&&0\\0&&&&0\esm}\restore\phX\ar[rd]\ar[ru]&&
 {}\save[]*{\bsm0&&&&1\\&0&&1\\1&&1&&0\\&1&&0\\0&&&&0\esm}\restore\phX\ar[rd]\ar[ru]&\\
 {}\save[]*{\bsm0&&&&0\\&0&&0\\1&&0&&0\\&0&&0\\0&&&&0\esm}\restore\phX\ar[ru]\ar@{.}[uu]&&
 {}\save[]*{\bsm0&&&&0\\&1&&0\\0&&1&&0\\&0&&1\\0&&&&1\esm}\restore\phX\ar[ru]\ar[rd]&&
 {}\save[]*{\bsm1&&&&0\\&1&&0\\0&&1&&0\\&0&&1\\0&&&&0\esm}\restore\phX\ar[ru]&&
 {}\save[]*{\bsm0&&&&0\\&0&&0\\0&&0&&1\\&0&&0\\0&&&&0\esm}\restore\phX\ar[ru]&&
 {}\save[]*{\bsm0&&&&0\\&0&&1\\0&&1&&0\\&1&&0\\1&&&&0\esm}\restore\phX\ar[ru]\ar[rd]&&
 {}\save[]*{\bsm0&&&&1\\&0&&1\\0&&1&&0\\&1&&0\\0&&&&0\esm}\restore\phX\ar[ru]&&
 {}\save[]*{\bsm0&&&&0\\&0&&0\\1&&0&&0\\&0&&0\\0&&&&0\esm}\restore\phX\ar@{.}[uu]\\
 &&&{}\save[]*{\bsm1&&&&0\\&1&&0\\0&&1&&0\\&0&&1\\0&&&&1\esm}\restore\phX\ar[ru]&&&&&&
 {}\save[]*{\bsm0&&&&1\\&0&&1\\0&&1&&0\\&1&&0\\1&&&&0\esm}\restore\phX\ar[ru]\\
}
\]

\newpage
%%%%%%%%%%%%%%%%%%%%%%%%%%%%%%%%%%%%
\subsection{The sets $R_{[m,n]}(i)$}\label{rootlist}
%%%%%%%%%%%%%%%%%%%%%%%%%%%%%%%%%%%%

\vspace{1.0cm}
\begin{center}
{\bf 
The sets $R_{[m,n]}(2)$, $[m,n] \in L_2$
}
\end{center}

{\begin{small}
\[
[[1,0], \left[ \begin {smallmatrix} &0&&0&
\\\noalign{\medskip}0&&1&&1\\\noalign{\medskip}&1&&2&
\\\noalign{\medskip}0&&1&&1\end {smallmatrix} \right] , \left[ 
\begin {smallmatrix} &0&&0&\\\noalign{\medskip}1&&1&&0
\\\noalign{\medskip}&2&&1&\\\noalign{\medskip}1&&1&&0
\end {smallmatrix} \right] ]
\] \[
[[1,1], \left[ \begin {smallmatrix} &1&&0&
\\\noalign{\medskip}1&&2&&1\\\noalign{\medskip}&2&&2&
\\\noalign{\medskip}0&&1&&1\end {smallmatrix} \right] , \left[ 
\begin {smallmatrix} &0&&1&\\\noalign{\medskip}1&&2&&1
\\\noalign{\medskip}&2&&2&\\\noalign{\medskip}1&&1&&0
\end {smallmatrix} \right] ]
\] \[
[[0,1], \left[ \begin {smallmatrix} &0&&1&
\\\noalign{\medskip}0&&1&&1\\\noalign{\medskip}&0&&1&
\\\noalign{\medskip}0&&0&&0\end {smallmatrix} \right] , \left[ 
\begin {smallmatrix} &1&&0&\\\noalign{\medskip}1&&1&&0
\\\noalign{\medskip}&1&&0&\\\noalign{\medskip}0&&0&&0
\end {smallmatrix} \right] ]
\] 
\end{small}}

\vspace{1.0cm}
\begin{center}
{\bf 
The sets $R_{[m,n]}(3)$, $[m,n] \in L_3$
}
\end{center}

{\begin{small}
\[
[[1,0], \left[ \begin {smallmatrix} &0&&0&
\\\noalign{\medskip}0&&1&&0\\\noalign{\medskip}&1&&1&
\\\noalign{\medskip}1&&0&&1\end {smallmatrix} \right] , \left[ 
\begin {smallmatrix} &0&&0&\\\noalign{\medskip}0&&1&&0
\\\noalign{\medskip}&1&&1&\\\noalign{\medskip}0&&1&&0
\end {smallmatrix} \right] , \left[ \begin {smallmatrix} &0&&0&
\\\noalign{\medskip}1&&0&&1\\\noalign{\medskip}&1&&1&
\\\noalign{\medskip}0&&1&&0\end {smallmatrix} \right] ]
\] \[
[[2,0], \left[ \begin {smallmatrix} &0&&0&
\\\noalign{\medskip}1&&1&&1\\\noalign{\medskip}&2&&2&
\\\noalign{\medskip}0&&2&&0\end {smallmatrix} \right] , \left[ 
\begin {smallmatrix} &0&&0&\\\noalign{\medskip}1&&1&&1
\\\noalign{\medskip}&2&&2&\\\noalign{\medskip}1&&1&&1
\end {smallmatrix} \right] , \left[ \begin {smallmatrix} &0&&0&
\\\noalign{\medskip}0&&2&&0\\\noalign{\medskip}&2&&2&
\\\noalign{\medskip}1&&1&&1\end {smallmatrix} \right] ]
\] \[
[[2,1], \left[ \begin {smallmatrix} &0&&0&
\\\noalign{\medskip}1&&2&&1\\\noalign{\medskip}&2&&2&
\\\noalign{\medskip}1&&1&&1\end {smallmatrix} \right] , \left[ 
\begin {smallmatrix} &1&&1&\\\noalign{\medskip}1&&3&&1
\\\noalign{\medskip}&3&&3&\\\noalign{\medskip}1&&2&&1
\end {smallmatrix} \right] , \left[ \begin {smallmatrix} &0&&0&
\\\noalign{\medskip}1&&1&&1\\\noalign{\medskip}&2&&2&
\\\noalign{\medskip}0&&1&&0\end {smallmatrix} \right] ]
\] \[
[[1,1], \left[ \begin {smallmatrix} &0&&0&
\\\noalign{\medskip}0&&1&&0\\\noalign{\medskip}&1&&1&
\\\noalign{\medskip}0&&0&&0\end {smallmatrix} \right] , \left[ 
\begin {smallmatrix} &0&&0&\\\noalign{\medskip}1&&1&&1
\\\noalign{\medskip}&1&&1&\\\noalign{\medskip}0&&1&&0
\end {smallmatrix} \right] , \left[ \begin {smallmatrix} &1&&1&
\\\noalign{\medskip}1&&2&&1\\\noalign{\medskip}&2&&2&
\\\noalign{\medskip}1&&1&&1\end {smallmatrix} \right] ]
\] \[
[[2,2], \left[ \begin {smallmatrix} &1&&1&
\\\noalign{\medskip}2&&3&&2\\\noalign{\medskip}&3&&3&
\\\noalign{\medskip}1&&2&&1\end {smallmatrix} \right] , \left[ 
\begin {smallmatrix} &1&&1&\\\noalign{\medskip}1&&3&&1
\\\noalign{\medskip}&3&&3&\\\noalign{\medskip}1&&1&&1
\end {smallmatrix} \right] , \left[ \begin {smallmatrix} &0&&0&
\\\noalign{\medskip}1&&2&&1\\\noalign{\medskip}&2&&2&
\\\noalign{\medskip}0&&1&&0\end {smallmatrix} \right] ]
\] \[
[[1,2], \left[ \begin {smallmatrix} &1&&1&
\\\noalign{\medskip}1&&2&&1\\\noalign{\medskip}&2&&2&
\\\noalign{\medskip}0&&1&&0\end {smallmatrix} \right] , \left[ 
\begin {smallmatrix} &0&&0&\\\noalign{\medskip}1&&1&&1
\\\noalign{\medskip}&1&&1&\\\noalign{\medskip}0&&0&&0
\end {smallmatrix} \right] , \left[ \begin {smallmatrix} &1&&1&
\\\noalign{\medskip}1&&3&&1\\\noalign{\medskip}&2&&2&
\\\noalign{\medskip}1&&1&&1\end {smallmatrix} \right] ]
\] \[
[[0,1], \left[ \begin {smallmatrix} &0&&0&
\\\noalign{\medskip}0&&1&&0\\\noalign{\medskip}&0&&0&
\\\noalign{\medskip}0&&0&&0\end {smallmatrix} \right] , \left[ 
\begin {smallmatrix} &1&&1&\\\noalign{\medskip}1&&1&&1
\\\noalign{\medskip}&1&&1&\\\noalign{\medskip}0&&1&&0
\end {smallmatrix} \right] , \left[ \begin {smallmatrix} &0&&0&
\\\noalign{\medskip}0&&0&&0\\\noalign{\medskip}&0&&0&
\\\noalign{\medskip}0&&-1&&0\end {smallmatrix} \right] ]
\] \[
[[0,2], \left[ \begin {smallmatrix} &1&&1&
\\\noalign{\medskip}1&&1&&1\\\noalign{\medskip}&1&&1&
\\\noalign{\medskip}0&&0&&0\end {smallmatrix} \right] , \left[ 
\begin {smallmatrix} &0&&0&\\\noalign{\medskip}0&&1&&0
\\\noalign{\medskip}&0&&0&\\\noalign{\medskip}0&&-1&&0
\end {smallmatrix} \right] , \left[ \begin {smallmatrix} &1&&1&
\\\noalign{\medskip}1&&2&&1\\\noalign{\medskip}&1&&1&
\\\noalign{\medskip}0&&1&&0\end {smallmatrix} \right] ]
\] 
\end{small}}

\newpage
\begin{center}
{\bf 
The sets $R_{[m,n]}(6)$, $[m,n] \in L_6$
}
\end{center}

{\begin{small}
\[
[[1,0], \left[ \begin {smallmatrix} &-1&&0&
\\\noalign{\medskip}0&&0&&0\\\noalign{\medskip}&0&&0&
\\\noalign{\medskip}0&&0&&0\end {smallmatrix} \right] , \left[ 
\begin {smallmatrix} &1&&0&\\\noalign{\medskip}1&&1&&0
\\\noalign{\medskip}&1&&1&\\\noalign{\medskip}0&&1&&1
\end {smallmatrix} \right] , \left[ \begin {smallmatrix} &0&&0&
\\\noalign{\medskip}0&&0&&0\\\noalign{\medskip}&0&&1&
\\\noalign{\medskip}0&&0&&0\end {smallmatrix} \right] , \left[ 
\begin {smallmatrix} &0&&-1&\\\noalign{\medskip}0&&0&&0
\\\noalign{\medskip}&0&&0&\\\noalign{\medskip}0&&0&&0
\end {smallmatrix} \right] , \left[ \begin {smallmatrix} &0&&1&
\\\noalign{\medskip}0&&1&&1\\\noalign{\medskip}&1&&1&
\\\noalign{\medskip}1&&1&&0\end {smallmatrix} \right] , \left[ 
\begin {smallmatrix} &0&&0&\\\noalign{\medskip}0&&0&&0
\\\noalign{\medskip}&1&&0&\\\noalign{\medskip}0&&0&&0
\end {smallmatrix} \right] ]
\] \[
[[2,0], \left[ \begin {smallmatrix} &-1&&0&
\\\noalign{\medskip}0&&0&&0\\\noalign{\medskip}&1&&0&
\\\noalign{\medskip}0&&0&&0\end {smallmatrix} \right] , \left[ 
\begin {smallmatrix} &0&&0&\\\noalign{\medskip}1&&1&&0
\\\noalign{\medskip}&1&&1&\\\noalign{\medskip}0&&1&&1
\end {smallmatrix} \right] , \left[ \begin {smallmatrix} &1&&0&
\\\noalign{\medskip}1&&1&&0\\\noalign{\medskip}&1&&2&
\\\noalign{\medskip}0&&1&&1\end {smallmatrix} \right] , \left[ 
\begin {smallmatrix} &0&&-1&\\\noalign{\medskip}0&&0&&0
\\\noalign{\medskip}&0&&1&\\\noalign{\medskip}0&&0&&0
\end {smallmatrix} \right] , \left[ \begin {smallmatrix} &0&&0&
\\\noalign{\medskip}0&&1&&1\\\noalign{\medskip}&1&&1&
\\\noalign{\medskip}1&&1&&0\end {smallmatrix} \right] , \left[ 
\begin {smallmatrix} &0&&1&\\\noalign{\medskip}0&&1&&1
\\\noalign{\medskip}&2&&1&\\\noalign{\medskip}1&&1&&0
\end {smallmatrix} \right] ]
\] \[
[[3,0], \left[ \begin {smallmatrix} &-1&&1&
\\\noalign{\medskip}0&&1&&1\\\noalign{\medskip}&2&&1&
\\\noalign{\medskip}1&&1&&0\end {smallmatrix} \right] , \left[ 
\begin {smallmatrix} &0&&0&\\\noalign{\medskip}1&&1&&0
\\\noalign{\medskip}&2&&1&\\\noalign{\medskip}0&&1&&1
\end {smallmatrix} \right] , \left[ \begin {smallmatrix} &0&&0&
\\\noalign{\medskip}1&&1&&0\\\noalign{\medskip}&1&&2&
\\\noalign{\medskip}0&&1&&1\end {smallmatrix} \right] , \left[ 
\begin {smallmatrix} &1&&-1&\\\noalign{\medskip}1&&1&&0
\\\noalign{\medskip}&1&&2&\\\noalign{\medskip}0&&1&&1
\end {smallmatrix} \right] , \left[ \begin {smallmatrix} &0&&0&
\\\noalign{\medskip}0&&1&&1\\\noalign{\medskip}&1&&2&
\\\noalign{\medskip}1&&1&&0\end {smallmatrix} \right] , \left[ 
\begin {smallmatrix} &0&&0&\\\noalign{\medskip}0&&1&&1
\\\noalign{\medskip}&2&&1&\\\noalign{\medskip}1&&1&&0
\end {smallmatrix} \right] ]
\] \[
[[4,0], \left[ \begin {smallmatrix} &1&&0&
\\\noalign{\medskip}1&&2&&1\\\noalign{\medskip}&2&&3&
\\\noalign{\medskip}1&&2&&1\end {smallmatrix} \right] , \left[ 
\begin {smallmatrix} &0&&0&\\\noalign{\medskip}0&&1&&1
\\\noalign{\medskip}&2&&2&\\\noalign{\medskip}1&&1&&0
\end {smallmatrix} \right] , \left[ \begin {smallmatrix} &-1&&0&
\\\noalign{\medskip}0&&1&&1\\\noalign{\medskip}&2&&1&
\\\noalign{\medskip}1&&1&&0\end {smallmatrix} \right] , \left[ 
\begin {smallmatrix} &0&&1&\\\noalign{\medskip}1&&2&&1
\\\noalign{\medskip}&3&&2&\\\noalign{\medskip}1&&2&&1
\end {smallmatrix} \right] , \left[ \begin {smallmatrix} &0&&0&
\\\noalign{\medskip}1&&1&&0\\\noalign{\medskip}&2&&2&
\\\noalign{\medskip}0&&1&&1\end {smallmatrix} \right] , \left[ 
\begin {smallmatrix} &0&&-1&\\\noalign{\medskip}1&&1&&0
\\\noalign{\medskip}&1&&2&\\\noalign{\medskip}0&&1&&1
\end {smallmatrix} \right] ]
\] \[
[[5,0], \left[ \begin {smallmatrix} &1&&0&
\\\noalign{\medskip}1&&2&&1\\\noalign{\medskip}&3&&3&
\\\noalign{\medskip}1&&2&&1\end {smallmatrix} \right] , \left[ 
\begin {smallmatrix} &-1&&0&\\\noalign{\medskip}0&&1&&
1\\\noalign{\medskip}&2&&2&\\\noalign{\medskip}1&&1&&0
\end {smallmatrix} \right] , \left[ \begin {smallmatrix} &0&&0&
\\\noalign{\medskip}1&&2&&1\\\noalign{\medskip}&3&&2&
\\\noalign{\medskip}1&&2&&1\end {smallmatrix} \right] , \left[ 
\begin {smallmatrix} &0&&1&\\\noalign{\medskip}1&&2&&1
\\\noalign{\medskip}&3&&3&\\\noalign{\medskip}1&&2&&1
\end {smallmatrix} \right] , \left[ \begin {smallmatrix} &0&&-1&
\\\noalign{\medskip}1&&1&&0\\\noalign{\medskip}&2&&2&
\\\noalign{\medskip}0&&1&&1\end {smallmatrix} \right] , \left[ 
\begin {smallmatrix} &0&&0&\\\noalign{\medskip}1&&2&&1
\\\noalign{\medskip}&2&&3&\\\noalign{\medskip}1&&2&&1
\end {smallmatrix} \right] ]
\] \[
[[5,1], \left[ \begin {smallmatrix} &0&&0&
\\\noalign{\medskip}1&&2&&1\\\noalign{\medskip}&2&&3&
\\\noalign{\medskip}0&&2&&1\end {smallmatrix} \right] , \left[ 
\begin {smallmatrix} &1&&0&\\\noalign{\medskip}2&&2&&1
\\\noalign{\medskip}&3&&3&\\\noalign{\medskip}1&&2&&1
\end {smallmatrix} \right] , \left[ \begin {smallmatrix} &0&&0&
\\\noalign{\medskip}0&&2&&1\\\noalign{\medskip}&2&&3&
\\\noalign{\medskip}1&&1&&1\end {smallmatrix} \right] , \left[ 
\begin {smallmatrix} &0&&0&\\\noalign{\medskip}1&&2&&1
\\\noalign{\medskip}&3&&2&\\\noalign{\medskip}1&&2&&0
\end {smallmatrix} \right] , \left[ \begin {smallmatrix} &0&&1&
\\\noalign{\medskip}1&&2&&2\\\noalign{\medskip}&3&&3&
\\\noalign{\medskip}1&&2&&1\end {smallmatrix} \right] , \left[ 
\begin {smallmatrix} &0&&0&\\\noalign{\medskip}1&&2&&0
\\\noalign{\medskip}&3&&2&\\\noalign{\medskip}1&&1&&1
\end {smallmatrix} \right] ]
\] \[
[[4,1], \left[ \begin {smallmatrix} &0&&0&
\\\noalign{\medskip}1&&1&&1\\\noalign{\medskip}&2&&2&
\\\noalign{\medskip}0&&1&&1\end {smallmatrix} \right] , \left[ 
\begin {smallmatrix} &0&&0&\\\noalign{\medskip}1&&2&&0
\\\noalign{\medskip}&2&&2&\\\noalign{\medskip}1&&1&&1
\end {smallmatrix} \right] , \left[ \begin {smallmatrix} &1&&0&
\\\noalign{\medskip}1&&2&&1\\\noalign{\medskip}&2&&3&
\\\noalign{\medskip}0&&2&&1\end {smallmatrix} \right] , \left[ 
\begin {smallmatrix} &0&&0&\\\noalign{\medskip}1&&1&&1
\\\noalign{\medskip}&2&&2&\\\noalign{\medskip}1&&1&&0
\end {smallmatrix} \right] , \left[ \begin {smallmatrix} &0&&0&
\\\noalign{\medskip}0&&2&&1\\\noalign{\medskip}&2&&2&
\\\noalign{\medskip}1&&1&&1\end {smallmatrix} \right] , \left[ 
\begin {smallmatrix} &0&&1&\\\noalign{\medskip}1&&2&&1
\\\noalign{\medskip}&3&&2&\\\noalign{\medskip}1&&2&&0
\end {smallmatrix} \right] ]
\] \[
[[3,1], \left[ \begin {smallmatrix} &0&&0&
\\\noalign{\medskip}1&&1&&0\\\noalign{\medskip}&2&&1&
\\\noalign{\medskip}0&&1&&0\end {smallmatrix} \right] , \left[ 
\begin {smallmatrix} &0&&0&\\\noalign{\medskip}1&&1&&1
\\\noalign{\medskip}&1&&2&\\\noalign{\medskip}0&&1&&1
\end {smallmatrix} \right] , \left[ \begin {smallmatrix} &1&&0&
\\\noalign{\medskip}1&&2&&0\\\noalign{\medskip}&2&&2&
\\\noalign{\medskip}1&&1&&1\end {smallmatrix} \right] , \left[ 
\begin {smallmatrix} &0&&0&\\\noalign{\medskip}0&&1&&1
\\\noalign{\medskip}&1&&2&\\\noalign{\medskip}0&&1&&0
\end {smallmatrix} \right] , \left[ \begin {smallmatrix} &0&&0&
\\\noalign{\medskip}1&&1&&1\\\noalign{\medskip}&2&&1&
\\\noalign{\medskip}1&&1&&0\end {smallmatrix} \right] , \left[ 
\begin {smallmatrix} &0&&1&\\\noalign{\medskip}0&&2&&1
\\\noalign{\medskip}&2&&2&\\\noalign{\medskip}1&&1&&1
\end {smallmatrix} \right] ]
\] \[
[[5,2], \left[ \begin {smallmatrix} &0&&0&
\\\noalign{\medskip}1&&2&&1\\\noalign{\medskip}&2&&3&
\\\noalign{\medskip}1&&1&&1\end {smallmatrix} \right] , \left[ 
\begin {smallmatrix} &1&&0&\\\noalign{\medskip}1&&3&&1
\\\noalign{\medskip}&3&&3&\\\noalign{\medskip}1&&2&&1
\end {smallmatrix} \right] , \left[ \begin {smallmatrix} &0&&1&
\\\noalign{\medskip}1&&2&&2\\\noalign{\medskip}&3&&3&
\\\noalign{\medskip}1&&2&&0\end {smallmatrix} \right] , \left[ 
\begin {smallmatrix} &0&&0&\\\noalign{\medskip}1&&2&&1
\\\noalign{\medskip}&3&&2&\\\noalign{\medskip}1&&1&&1
\end {smallmatrix} \right] , \left[ \begin {smallmatrix} &0&&1&
\\\noalign{\medskip}1&&3&&1\\\noalign{\medskip}&3&&3&
\\\noalign{\medskip}1&&2&&1\end {smallmatrix} \right] , \left[ 
\begin {smallmatrix} &1&&0&\\\noalign{\medskip}2&&2&&1
\\\noalign{\medskip}&3&&3&\\\noalign{\medskip}0&&2&&1
\end {smallmatrix} \right] ]
\] \[
[[2,1], \left[ \begin {smallmatrix} &0&&0&
\\\noalign{\medskip}0&&1&&0\\\noalign{\medskip}&1&&1&
\\\noalign{\medskip}0&&0&&1\end {smallmatrix} \right] , \left[ 
\begin {smallmatrix} &0&&0&\\\noalign{\medskip}1&&1&&0
\\\noalign{\medskip}&1&&1&\\\noalign{\medskip}0&&1&&0
\end {smallmatrix} \right] , \left[ \begin {smallmatrix} &1&&0&
\\\noalign{\medskip}1&&1&&1\\\noalign{\medskip}&1&&2&
\\\noalign{\medskip}0&&1&&1\end {smallmatrix} \right] , \left[ 
\begin {smallmatrix} &0&&0&\\\noalign{\medskip}0&&1&&0
\\\noalign{\medskip}&1&&1&\\\noalign{\medskip}1&&0&&0
\end {smallmatrix} \right] , \left[ \begin {smallmatrix} &0&&0&
\\\noalign{\medskip}0&&1&&1\\\noalign{\medskip}&1&&1&
\\\noalign{\medskip}0&&1&&0\end {smallmatrix} \right] , \left[ 
\begin {smallmatrix} &0&&1&\\\noalign{\medskip}1&&1&&1
\\\noalign{\medskip}&2&&1&\\\noalign{\medskip}1&&1&&0
\end {smallmatrix} \right] ]
\] \[
[[4,2], \left[ \begin {smallmatrix} &0&&0&
\\\noalign{\medskip}0&&2&&1\\\noalign{\medskip}&2&&2&
\\\noalign{\medskip}1&&1&&0\end {smallmatrix} \right] , \left[ 
\begin {smallmatrix} &0&&1&\\\noalign{\medskip}1&&2&&2
\\\noalign{\medskip}&3&&2&\\\noalign{\medskip}1&&2&&0
\end {smallmatrix} \right] , \left[ \begin {smallmatrix} &0&&1&
\\\noalign{\medskip}1&&2&&1\\\noalign{\medskip}&3&&2&
\\\noalign{\medskip}1&&1&&1\end {smallmatrix} \right] , \left[ 
\begin {smallmatrix} &0&&0&\\\noalign{\medskip}1&&2&&0
\\\noalign{\medskip}&2&&2&\\\noalign{\medskip}0&&1&&1
\end {smallmatrix} \right] , \left[ \begin {smallmatrix} &1&&0&
\\\noalign{\medskip}2&&2&&1\\\noalign{\medskip}&2&&3&
\\\noalign{\medskip}0&&2&&1\end {smallmatrix} \right] , \left[ 
\begin {smallmatrix} &1&&0&\\\noalign{\medskip}1&&2&&1
\\\noalign{\medskip}&2&&3&\\\noalign{\medskip}1&&1&&1
\end {smallmatrix} \right] ]
\] \[
[[5,3], \left[ \begin {smallmatrix} &1&&0&
\\\noalign{\medskip}2&&3&&1\\\noalign{\medskip}&3&&3&
\\\noalign{\medskip}1&&2&&1\end {smallmatrix} \right] , \left[ 
\begin {smallmatrix} &1&&1&\\\noalign{\medskip}1&&3&&2
\\\noalign{\medskip}&3&&4&\\\noalign{\medskip}1&&2&&1
\end {smallmatrix} \right] , \left[ \begin {smallmatrix} &0&&0&
\\\noalign{\medskip}1&&2&&1\\\noalign{\medskip}&3&&2&
\\\noalign{\medskip}1&&1&&0\end {smallmatrix} \right] , \left[ 
\begin {smallmatrix} &0&&1&\\\noalign{\medskip}1&&3&&2
\\\noalign{\medskip}&3&&3&\\\noalign{\medskip}1&&2&&1
\end {smallmatrix} \right] , \left[ \begin {smallmatrix} &1&&1&
\\\noalign{\medskip}2&&3&&1\\\noalign{\medskip}&4&&3&
\\\noalign{\medskip}1&&2&&1\end {smallmatrix} \right] , \left[ 
\begin {smallmatrix} &0&&0&\\\noalign{\medskip}1&&2&&1
\\\noalign{\medskip}&2&&3&\\\noalign{\medskip}0&&1&&1
\end {smallmatrix} \right] ]
\] \[
[[3,2], \left[ \begin {smallmatrix} &0&&0&
\\\noalign{\medskip}1&&1&&1\\\noalign{\medskip}&2&&1&
\\\noalign{\medskip}0&&1&&0\end {smallmatrix} \right] , \left[ 
\begin {smallmatrix} &0&&1&\\\noalign{\medskip}1&&2&&1
\\\noalign{\medskip}&2&&2&\\\noalign{\medskip}1&&1&&1
\end {smallmatrix} \right] , \left[ \begin {smallmatrix} &1&&0&
\\\noalign{\medskip}1&&2&&0\\\noalign{\medskip}&2&&2&
\\\noalign{\medskip}0&&1&&1\end {smallmatrix} \right] , \left[ 
\begin {smallmatrix} &0&&0&\\\noalign{\medskip}1&&1&&1
\\\noalign{\medskip}&1&&2&\\\noalign{\medskip}0&&1&&0
\end {smallmatrix} \right] , \left[ \begin {smallmatrix} &1&&0&
\\\noalign{\medskip}1&&2&&1\\\noalign{\medskip}&2&&2&
\\\noalign{\medskip}1&&1&&1\end {smallmatrix} \right] , \left[ 
\begin {smallmatrix} &0&&1&\\\noalign{\medskip}0&&2&&1
\\\noalign{\medskip}&2&&2&\\\noalign{\medskip}1&&1&&0
\end {smallmatrix} \right] ]
\] \[
[[4,3], \left[ \begin {smallmatrix} &0&&0&
\\\noalign{\medskip}1&&2&&1\\\noalign{\medskip}&2&&2&
\\\noalign{\medskip}1&&1&&0\end {smallmatrix} \right] , \left[ 
\begin {smallmatrix} &1&&1&\\\noalign{\medskip}1&&3&&2
\\\noalign{\medskip}&3&&3&\\\noalign{\medskip}1&&2&&1
\end {smallmatrix} \right] , \left[ \begin {smallmatrix} &0&&1&
\\\noalign{\medskip}1&&2&&1\\\noalign{\medskip}&3&&2&
\\\noalign{\medskip}1&&1&&0\end {smallmatrix} \right] , \left[ 
\begin {smallmatrix} &0&&0&\\\noalign{\medskip}1&&2&&1
\\\noalign{\medskip}&2&&2&\\\noalign{\medskip}0&&1&&1
\end {smallmatrix} \right] , \left[ \begin {smallmatrix} &1&&1&
\\\noalign{\medskip}2&&3&&1\\\noalign{\medskip}&3&&3&
\\\noalign{\medskip}1&&2&&1\end {smallmatrix} \right] , \left[ 
\begin {smallmatrix} &1&&0&\\\noalign{\medskip}1&&2&&1
\\\noalign{\medskip}&2&&3&\\\noalign{\medskip}0&&1&&1
\end {smallmatrix} \right] ]
\] \[
[[5,4], \left[ \begin {smallmatrix} &0&&1&
\\\noalign{\medskip}1&&3&&2\\\noalign{\medskip}&3&&3&
\\\noalign{\medskip}1&&2&&0\end {smallmatrix} \right] , \left[ 
\begin {smallmatrix} &1&&1&\\\noalign{\medskip}2&&3&&2
\\\noalign{\medskip}&4&&3&\\\noalign{\medskip}1&&2&&1
\end {smallmatrix} \right] , \left[ \begin {smallmatrix} &0&&1&
\\\noalign{\medskip}1&&3&&1\\\noalign{\medskip}&3&&3&
\\\noalign{\medskip}1&&1&&1\end {smallmatrix} \right] , \left[ 
\begin {smallmatrix} &1&&0&\\\noalign{\medskip}2&&3&&1
\\\noalign{\medskip}&3&&3&\\\noalign{\medskip}0&&2&&1
\end {smallmatrix} \right] , \left[ \begin {smallmatrix} &1&&1&
\\\noalign{\medskip}2&&3&&2\\\noalign{\medskip}&3&&4&
\\\noalign{\medskip}1&&2&&1\end {smallmatrix} \right] , \left[ 
\begin {smallmatrix} &1&&0&\\\noalign{\medskip}1&&3&&1
\\\noalign{\medskip}&3&&3&\\\noalign{\medskip}1&&1&&1
\end {smallmatrix} \right] ]
\] \[
[[1,1], \left[ \begin {smallmatrix} &0&&0&
\\\noalign{\medskip}1&&0&&0\\\noalign{\medskip}&1&&0&
\\\noalign{\medskip}0&&0&&0\end {smallmatrix} \right] , \left[ 
\begin {smallmatrix} &0&&0&\\\noalign{\medskip}0&&1&&0
\\\noalign{\medskip}&0&&1&\\\noalign{\medskip}0&&0&&1
\end {smallmatrix} \right] , \left[ \begin {smallmatrix} &1&&0&
\\\noalign{\medskip}1&&1&&0\\\noalign{\medskip}&1&&1&
\\\noalign{\medskip}0&&1&&0\end {smallmatrix} \right] , \left[ 
\begin {smallmatrix} &0&&0&\\\noalign{\medskip}0&&0&&1
\\\noalign{\medskip}&0&&1&\\\noalign{\medskip}0&&0&&0
\end {smallmatrix} \right] , \left[ \begin {smallmatrix} &0&&0&
\\\noalign{\medskip}0&&1&&0\\\noalign{\medskip}&1&&0&
\\\noalign{\medskip}1&&0&&0\end {smallmatrix} \right] , \left[ 
\begin {smallmatrix} &0&&1&\\\noalign{\medskip}0&&1&&1
\\\noalign{\medskip}&1&&1&\\\noalign{\medskip}0&&1&&0
\end {smallmatrix} \right] ]
\] \[
[[2,2], \left[ \begin {smallmatrix} &0&&1&
\\\noalign{\medskip}1&&1&&1\\\noalign{\medskip}&2&&1&
\\\noalign{\medskip}0&&1&&0\end {smallmatrix} \right] , \left[ 
\begin {smallmatrix} &0&&0&\\\noalign{\medskip}1&&1&&0
\\\noalign{\medskip}&1&&1&\\\noalign{\medskip}0&&0&&1
\end {smallmatrix} \right] , \left[ \begin {smallmatrix} &1&&0&
\\\noalign{\medskip}1&&2&&0\\\noalign{\medskip}&1&&2&
\\\noalign{\medskip}0&&1&&1\end {smallmatrix} \right] , \left[ 
\begin {smallmatrix} &1&&0&\\\noalign{\medskip}1&&1&&1
\\\noalign{\medskip}&1&&2&\\\noalign{\medskip}0&&1&&0
\end {smallmatrix} \right] , \left[ \begin {smallmatrix} &0&&0&
\\\noalign{\medskip}0&&1&&1\\\noalign{\medskip}&1&&1&
\\\noalign{\medskip}1&&0&&0\end {smallmatrix} \right] , \left[ 
\begin {smallmatrix} &0&&1&\\\noalign{\medskip}0&&2&&1
\\\noalign{\medskip}&2&&1&\\\noalign{\medskip}1&&1&&0
\end {smallmatrix} \right] ]
\] \[
[[3,3], \left[ \begin {smallmatrix} &0&&1&
\\\noalign{\medskip}1&&2&&1\\\noalign{\medskip}&3&&1&
\\\noalign{\medskip}1&&1&&0\end {smallmatrix} \right] , \left[ 
\begin {smallmatrix} &0&&1&\\\noalign{\medskip}1&&2&&1
\\\noalign{\medskip}&2&&2&\\\noalign{\medskip}0&&1&&1
\end {smallmatrix} \right] , \left[ \begin {smallmatrix} &1&&0&
\\\noalign{\medskip}2&&2&&0\\\noalign{\medskip}&2&&2&
\\\noalign{\medskip}0&&1&&1\end {smallmatrix} \right] , \left[ 
\begin {smallmatrix} &1&&0&\\\noalign{\medskip}1&&2&&1
\\\noalign{\medskip}&1&&3&\\\noalign{\medskip}0&&1&&1
\end {smallmatrix} \right] , \left[ \begin {smallmatrix} &1&&0&
\\\noalign{\medskip}1&&2&&1\\\noalign{\medskip}&2&&2&
\\\noalign{\medskip}1&&1&&0\end {smallmatrix} \right] , \left[ 
\begin {smallmatrix} &0&&1&\\\noalign{\medskip}0&&2&&2
\\\noalign{\medskip}&2&&2&\\\noalign{\medskip}1&&1&&0
\end {smallmatrix} \right] ]
\] \[
[[4,4], \left[ \begin {smallmatrix} &1&&0&
\\\noalign{\medskip}1&&3&&1\\\noalign{\medskip}&2&&3&
\\\noalign{\medskip}1&&1&&1\end {smallmatrix} \right] , \left[ 
\begin {smallmatrix} &1&&1&\\\noalign{\medskip}1&&3&&2
\\\noalign{\medskip}&3&&3&\\\noalign{\medskip}1&&2&&0
\end {smallmatrix} \right] , \left[ \begin {smallmatrix} &0&&1&
\\\noalign{\medskip}1&&2&&2\\\noalign{\medskip}&3&&2&
\\\noalign{\medskip}1&&1&&0\end {smallmatrix} \right] , \left[ 
\begin {smallmatrix} &0&&1&\\\noalign{\medskip}1&&3&&1
\\\noalign{\medskip}&3&&2&\\\noalign{\medskip}1&&1&&1
\end {smallmatrix} \right] , \left[ \begin {smallmatrix} &1&&1&
\\\noalign{\medskip}2&&3&&1\\\noalign{\medskip}&3&&3&
\\\noalign{\medskip}0&&2&&1\end {smallmatrix} \right] , \left[ 
\begin {smallmatrix} &1&&0&\\\noalign{\medskip}2&&2&&1
\\\noalign{\medskip}&2&&3&\\\noalign{\medskip}0&&1&&1
\end {smallmatrix} \right] ]
\] \[
[[5,5], \left[ \begin {smallmatrix} &1&&1&
\\\noalign{\medskip}1&&4&&2\\\noalign{\medskip}&3&&4&
\\\noalign{\medskip}1&&2&&1\end {smallmatrix} \right] , \left[ 
\begin {smallmatrix} &1&&1&\\\noalign{\medskip}2&&3&&2
\\\noalign{\medskip}&4&&3&\\\noalign{\medskip}1&&2&&0
\end {smallmatrix} \right] , \left[ \begin {smallmatrix} &0&&1&
\\\noalign{\medskip}1&&3&&2\\\noalign{\medskip}&3&&3&
\\\noalign{\medskip}1&&1&&1\end {smallmatrix} \right] , \left[ 
\begin {smallmatrix} &1&&1&\\\noalign{\medskip}2&&4&&1
\\\noalign{\medskip}&4&&3&\\\noalign{\medskip}1&&2&&1
\end {smallmatrix} \right] , \left[ \begin {smallmatrix} &1&&1&
\\\noalign{\medskip}2&&3&&2\\\noalign{\medskip}&3&&4&
\\\noalign{\medskip}0&&2&&1\end {smallmatrix} \right] , \left[ 
\begin {smallmatrix} &1&&0&\\\noalign{\medskip}2&&3&&1
\\\noalign{\medskip}&3&&3&\\\noalign{\medskip}1&&1&&1
\end {smallmatrix} \right] ]
\] \[
[[4,5], \left[ \begin {smallmatrix} &1&&1&
\\\noalign{\medskip}2&&3&&2\\\noalign{\medskip}&3&&3&
\\\noalign{\medskip}1&&2&&0\end {smallmatrix} \right] , \left[ 
\begin {smallmatrix} &1&&1&\\\noalign{\medskip}1&&3&&2
\\\noalign{\medskip}&3&&3&\\\noalign{\medskip}1&&1&&1
\end {smallmatrix} \right] , \left[ \begin {smallmatrix} &0&&1&
\\\noalign{\medskip}1&&3&&1\\\noalign{\medskip}&3&&2&
\\\noalign{\medskip}1&&1&&0\end {smallmatrix} \right] , \left[ 
\begin {smallmatrix} &1&&1&\\\noalign{\medskip}2&&3&&2
\\\noalign{\medskip}&3&&3&\\\noalign{\medskip}0&&2&&1
\end {smallmatrix} \right] , \left[ \begin {smallmatrix} &1&&1&
\\\noalign{\medskip}2&&3&&1\\\noalign{\medskip}&3&&3&
\\\noalign{\medskip}1&&1&&1\end {smallmatrix} \right] , \left[ 
\begin {smallmatrix} &1&&0&\\\noalign{\medskip}1&&3&&1
\\\noalign{\medskip}&2&&3&\\\noalign{\medskip}0&&1&&1
\end {smallmatrix} \right] ]
\] \[
[[3,4], \left[ \begin {smallmatrix} &1&&1&
\\\noalign{\medskip}1&&3&&1\\\noalign{\medskip}&2&&3&
\\\noalign{\medskip}1&&1&&1\end {smallmatrix} \right] , \left[ 
\begin {smallmatrix} &1&&0&\\\noalign{\medskip}1&&2&&1
\\\noalign{\medskip}&2&&2&\\\noalign{\medskip}0&&1&&0
\end {smallmatrix} \right] , \left[ \begin {smallmatrix} &0&&1&
\\\noalign{\medskip}1&&2&&2\\\noalign{\medskip}&2&&2&
\\\noalign{\medskip}1&&1&&0\end {smallmatrix} \right] , \left[ 
\begin {smallmatrix} &1&&1&\\\noalign{\medskip}1&&3&&1
\\\noalign{\medskip}&3&&2&\\\noalign{\medskip}1&&1&&1
\end {smallmatrix} \right] , \left[ \begin {smallmatrix} &0&&1&
\\\noalign{\medskip}1&&2&&1\\\noalign{\medskip}&2&&2&
\\\noalign{\medskip}0&&1&&0\end {smallmatrix} \right] , \left[ 
\begin {smallmatrix} &1&&0&\\\noalign{\medskip}2&&2&&1
\\\noalign{\medskip}&2&&2&\\\noalign{\medskip}0&&1&&1
\end {smallmatrix} \right] ]
\] \[
[[2,3], \left[ \begin {smallmatrix} &1&&1&
\\\noalign{\medskip}1&&2&&1\\\noalign{\medskip}&2&&2&
\\\noalign{\medskip}0&&1&&1\end {smallmatrix} \right] , \left[ 
\begin {smallmatrix} &0&&0&\\\noalign{\medskip}1&&1&&0
\\\noalign{\medskip}&1&&1&\\\noalign{\medskip}0&&0&&0
\end {smallmatrix} \right] , \left[ \begin {smallmatrix} &1&&0&
\\\noalign{\medskip}1&&2&&1\\\noalign{\medskip}&1&&2&
\\\noalign{\medskip}0&&1&&1\end {smallmatrix} \right] , \left[ 
\begin {smallmatrix} &1&&1&\\\noalign{\medskip}1&&2&&1
\\\noalign{\medskip}&2&&2&\\\noalign{\medskip}1&&1&&0
\end {smallmatrix} \right] , \left[ \begin {smallmatrix} &0&&0&
\\\noalign{\medskip}0&&1&&1\\\noalign{\medskip}&1&&1&
\\\noalign{\medskip}0&&0&&0\end {smallmatrix} \right] , \left[ 
\begin {smallmatrix} &0&&1&\\\noalign{\medskip}1&&2&&1
\\\noalign{\medskip}&2&&1&\\\noalign{\medskip}1&&1&&0
\end {smallmatrix} \right] ]
\] \[
[[3,5], \left[ \begin {smallmatrix} &1&&1&
\\\noalign{\medskip}1&&3&&2\\\noalign{\medskip}&2&&3&
\\\noalign{\medskip}1&&1&&1\end {smallmatrix} \right] , \left[ 
\begin {smallmatrix} &1&&1&\\\noalign{\medskip}1&&3&&1
\\\noalign{\medskip}&3&&2&\\\noalign{\medskip}1&&1&&0
\end {smallmatrix} \right] , \left[ \begin {smallmatrix} &0&&1&
\\\noalign{\medskip}1&&2&&2\\\noalign{\medskip}&2&&2&
\\\noalign{\medskip}0&&1&&0\end {smallmatrix} \right] , \left[ 
\begin {smallmatrix} &1&&1&\\\noalign{\medskip}2&&3&&1
\\\noalign{\medskip}&3&&2&\\\noalign{\medskip}1&&1&&1
\end {smallmatrix} \right] , \left[ \begin {smallmatrix} &1&&1&
\\\noalign{\medskip}1&&3&&1\\\noalign{\medskip}&2&&3&
\\\noalign{\medskip}0&&1&&1\end {smallmatrix} \right] , \left[ 
\begin {smallmatrix} &1&&0&\\\noalign{\medskip}2&&2&&1
\\\noalign{\medskip}&2&&2&\\\noalign{\medskip}0&&1&&0
\end {smallmatrix} \right] ]
\] \[
[[1,2], \left[ \begin {smallmatrix} &1&&0&
\\\noalign{\medskip}1&&1&&0\\\noalign{\medskip}&1&&1&
\\\noalign{\medskip}0&&0&&1\end {smallmatrix} \right] , \left[ 
\begin {smallmatrix} &0&&0&\\\noalign{\medskip}0&&1&&0
\\\noalign{\medskip}&0&&1&\\\noalign{\medskip}0&&0&&0
\end {smallmatrix} \right] , \left[ \begin {smallmatrix} &1&&0&
\\\noalign{\medskip}1&&1&&1\\\noalign{\medskip}&1&&1&
\\\noalign{\medskip}0&&1&&0\end {smallmatrix} \right] , \left[ 
\begin {smallmatrix} &0&&1&\\\noalign{\medskip}0&&1&&1
\\\noalign{\medskip}&1&&1&\\\noalign{\medskip}1&&0&&0
\end {smallmatrix} \right] , \left[ \begin {smallmatrix} &0&&0&
\\\noalign{\medskip}0&&1&&0\\\noalign{\medskip}&1&&0&
\\\noalign{\medskip}0&&0&&0\end {smallmatrix} \right] , \left[ 
\begin {smallmatrix} &0&&1&\\\noalign{\medskip}1&&1&&1
\\\noalign{\medskip}&1&&1&\\\noalign{\medskip}0&&1&&0
\end {smallmatrix} \right] ]
\] \[
[[2,4], \left[ \begin {smallmatrix} &1&&1&
\\\noalign{\medskip}2&&2&&1\\\noalign{\medskip}&2&&2&
\\\noalign{\medskip}0&&1&&1\end {smallmatrix} \right] , \left[ 
\begin {smallmatrix} &1&&0&\\\noalign{\medskip}1&&2&&0
\\\noalign{\medskip}&1&&2&\\\noalign{\medskip}0&&0&&1
\end {smallmatrix} \right] , \left[ \begin {smallmatrix} &1&&0&
\\\noalign{\medskip}1&&2&&1\\\noalign{\medskip}&1&&2&
\\\noalign{\medskip}0&&1&&0\end {smallmatrix} \right] , \left[ 
\begin {smallmatrix} &1&&1&\\\noalign{\medskip}1&&2&&2
\\\noalign{\medskip}&2&&2&\\\noalign{\medskip}1&&1&&0
\end {smallmatrix} \right] , \left[ \begin {smallmatrix} &0&&1&
\\\noalign{\medskip}0&&2&&1\\\noalign{\medskip}&2&&1&
\\\noalign{\medskip}1&&0&&0\end {smallmatrix} \right] , \left[ 
\begin {smallmatrix} &0&&1&\\\noalign{\medskip}1&&2&&1
\\\noalign{\medskip}&2&&1&\\\noalign{\medskip}0&&1&&0
\end {smallmatrix} \right] ]
\] \[
[[2,5], \left[ \begin {smallmatrix} &1&&1&
\\\noalign{\medskip}1&&3&&1\\\noalign{\medskip}&2&&2&
\\\noalign{\medskip}1&&1&&0\end {smallmatrix} \right] , \left[ 
\begin {smallmatrix} &1&&1&\\\noalign{\medskip}1&&2&&2
\\\noalign{\medskip}&2&&2&\\\noalign{\medskip}0&&1&&0
\end {smallmatrix} \right] , \left[ \begin {smallmatrix} &0&&1&
\\\noalign{\medskip}1&&2&&1\\\noalign{\medskip}&2&&1&
\\\noalign{\medskip}1&&0&&0\end {smallmatrix} \right] , \left[ 
\begin {smallmatrix} &1&&1&\\\noalign{\medskip}1&&3&&1
\\\noalign{\medskip}&2&&2&\\\noalign{\medskip}0&&1&&1
\end {smallmatrix} \right] , \left[ \begin {smallmatrix} &1&&1&
\\\noalign{\medskip}2&&2&&1\\\noalign{\medskip}&2&&2&
\\\noalign{\medskip}0&&1&&0\end {smallmatrix} \right] , \left[ 
\begin {smallmatrix} &1&&0&\\\noalign{\medskip}1&&2&&1
\\\noalign{\medskip}&1&&2&\\\noalign{\medskip}0&&0&&1
\end {smallmatrix} \right] ]
\] \[
[[1,3], \left[ \begin {smallmatrix} &0&&1&
\\\noalign{\medskip}0&&1&&1\\\noalign{\medskip}&1&&1&
\\\noalign{\medskip}0&&0&&0\end {smallmatrix} \right] , \left[ 
\begin {smallmatrix} &0&&0&\\\noalign{\medskip}1&&1&&0
\\\noalign{\medskip}&1&&0&\\\noalign{\medskip}0&&0&&0
\end {smallmatrix} \right] , \left[ \begin {smallmatrix} &1&&1&
\\\noalign{\medskip}1&&2&&1\\\noalign{\medskip}&1&&2&
\\\noalign{\medskip}0&&1&&1\end {smallmatrix} \right] , \left[ 
\begin {smallmatrix} &1&&0&\\\noalign{\medskip}1&&1&&0
\\\noalign{\medskip}&1&&1&\\\noalign{\medskip}0&&0&&0
\end {smallmatrix} \right] , \left[ \begin {smallmatrix} &0&&0&
\\\noalign{\medskip}0&&1&&1\\\noalign{\medskip}&0&&1&
\\\noalign{\medskip}0&&0&&0\end {smallmatrix} \right] , \left[ 
\begin {smallmatrix} &1&&1&\\\noalign{\medskip}1&&2&&1
\\\noalign{\medskip}&2&&1&\\\noalign{\medskip}1&&1&&0
\end {smallmatrix} \right] ]
\] \[
[[1,4], \left[ \begin {smallmatrix} &1&&1&
\\\noalign{\medskip}1&&2&&1\\\noalign{\medskip}&2&&1&
\\\noalign{\medskip}0&&1&&0\end {smallmatrix} \right] , \left[ 
\begin {smallmatrix} &0&&1&\\\noalign{\medskip}1&&1&&1
\\\noalign{\medskip}&1&&1&\\\noalign{\medskip}0&&0&&0
\end {smallmatrix} \right] , \left[ \begin {smallmatrix} &1&&0&
\\\noalign{\medskip}1&&2&&0\\\noalign{\medskip}&1&&1&
\\\noalign{\medskip}0&&0&&1\end {smallmatrix} \right] , \left[ 
\begin {smallmatrix} &1&&1&\\\noalign{\medskip}1&&2&&1
\\\noalign{\medskip}&1&&2&\\\noalign{\medskip}0&&1&&0
\end {smallmatrix} \right] , \left[ \begin {smallmatrix} &1&&0&
\\\noalign{\medskip}1&&1&&1\\\noalign{\medskip}&1&&1&
\\\noalign{\medskip}0&&0&&0\end {smallmatrix} \right] , \left[ 
\begin {smallmatrix} &0&&1&\\\noalign{\medskip}0&&2&&1
\\\noalign{\medskip}&1&&1&\\\noalign{\medskip}1&&0&&0
\end {smallmatrix} \right] ]
\] \[
[[1,5], \left[ \begin {smallmatrix} &1&&1&
\\\noalign{\medskip}1&&2&&1\\\noalign{\medskip}&2&&1&
\\\noalign{\medskip}1&&0&&0\end {smallmatrix} \right] , \left[ 
\begin {smallmatrix} &0&&1&\\\noalign{\medskip}0&&2&&1
\\\noalign{\medskip}&1&&1&\\\noalign{\medskip}0&&0&&0
\end {smallmatrix} \right] , \left[ \begin {smallmatrix} &1&&1&
\\\noalign{\medskip}2&&2&&1\\\noalign{\medskip}&2&&1&
\\\noalign{\medskip}0&&1&&0\end {smallmatrix} \right] , \left[ 
\begin {smallmatrix} &1&&1&\\\noalign{\medskip}1&&2&&1
\\\noalign{\medskip}&1&&2&\\\noalign{\medskip}0&&0&&1
\end {smallmatrix} \right] , \left[ \begin {smallmatrix} &1&&0&
\\\noalign{\medskip}1&&2&&0\\\noalign{\medskip}&1&&1&
\\\noalign{\medskip}0&&0&&0\end {smallmatrix} \right] , \left[ 
\begin {smallmatrix} &1&&1&\\\noalign{\medskip}1&&2&&2
\\\noalign{\medskip}&1&&2&\\\noalign{\medskip}0&&1&&0
\end {smallmatrix} \right] ]
\] \[
[[0,1], \left[ \begin {smallmatrix} &0&&0&
\\\noalign{\medskip}1&&0&&0\\\noalign{\medskip}&0&&0&
\\\noalign{\medskip}0&&0&&0\end {smallmatrix} \right] , \left[ 
\begin {smallmatrix} &1&&0&\\\noalign{\medskip}0&&1&&0
\\\noalign{\medskip}&0&&1&\\\noalign{\medskip}0&&0&&1
\end {smallmatrix} \right] , \left[ \begin {smallmatrix} &0&&0&
\\\noalign{\medskip}0&&0&&0\\\noalign{\medskip}&0&&0&
\\\noalign{\medskip}0&&0&&-1\end {smallmatrix} \right] , \left[ 
\begin {smallmatrix} &0&&0&\\\noalign{\medskip}0&&0&&1
\\\noalign{\medskip}&0&&0&\\\noalign{\medskip}0&&0&&0
\end {smallmatrix} \right] , \left[ \begin {smallmatrix} &0&&1&
\\\noalign{\medskip}0&&1&&0\\\noalign{\medskip}&1&&0&
\\\noalign{\medskip}1&&0&&0\end {smallmatrix} \right] , \left[ 
\begin {smallmatrix} &0&&0&\\\noalign{\medskip}0&&0&&0
\\\noalign{\medskip}&0&&0&\\\noalign{\medskip}-1&&0&&0
\end {smallmatrix} \right] ]
\] \[
[[0,2], \left[ \begin {smallmatrix} &0&&0&
\\\noalign{\medskip}1&&0&&0\\\noalign{\medskip}&0&&0&
\\\noalign{\medskip}-1&&0&&0\end {smallmatrix} \right] , \left[ 
\begin {smallmatrix} &1&&0&\\\noalign{\medskip}1&&1&&0
\\\noalign{\medskip}&0&&1&\\\noalign{\medskip}0&&0&&1
\end {smallmatrix} \right] , \left[ \begin {smallmatrix} &1&&0&
\\\noalign{\medskip}0&&1&&0\\\noalign{\medskip}&0&&1&
\\\noalign{\medskip}0&&0&&0\end {smallmatrix} \right] , \left[ 
\begin {smallmatrix} &0&&0&\\\noalign{\medskip}0&&0&&1
\\\noalign{\medskip}&0&&0&\\\noalign{\medskip}0&&0&&-1
\end {smallmatrix} \right] , \left[ \begin {smallmatrix} &0&&1&
\\\noalign{\medskip}0&&1&&1\\\noalign{\medskip}&1&&0&
\\\noalign{\medskip}1&&0&&0\end {smallmatrix} \right] , \left[ 
\begin {smallmatrix} &0&&1&\\\noalign{\medskip}0&&1&&0
\\\noalign{\medskip}&1&&0&\\\noalign{\medskip}0&&0&&0
\end {smallmatrix} \right] ]
\] \[
[[0,3], \left[ \begin {smallmatrix} &0&&1&
\\\noalign{\medskip}1&&1&&0\\\noalign{\medskip}&1&&0&
\\\noalign{\medskip}0&&0&&0\end {smallmatrix} \right] , \left[ 
\begin {smallmatrix} &1&&0&\\\noalign{\medskip}1&&1&&0
\\\noalign{\medskip}&0&&1&\\\noalign{\medskip}-1&&0&&1
\end {smallmatrix} \right] , \left[ \begin {smallmatrix} &1&&0&
\\\noalign{\medskip}1&&1&&0\\\noalign{\medskip}&0&&1&
\\\noalign{\medskip}0&&0&&0\end {smallmatrix} \right] , \left[ 
\begin {smallmatrix} &1&&0&\\\noalign{\medskip}0&&1&&1
\\\noalign{\medskip}&0&&1&\\\noalign{\medskip}0&&0&&0
\end {smallmatrix} \right] , \left[ \begin {smallmatrix} &0&&1&
\\\noalign{\medskip}0&&1&&1\\\noalign{\medskip}&1&&0&
\\\noalign{\medskip}1&&0&&-1\end {smallmatrix} \right] , \left[ 
\begin {smallmatrix} &0&&1&\\\noalign{\medskip}0&&1&&1
\\\noalign{\medskip}&1&&0&\\\noalign{\medskip}0&&0&&0
\end {smallmatrix} \right] ]
\] \[
[[0,4], \left[ \begin {smallmatrix} &1&&1&
\\\noalign{\medskip}0&&2&&1\\\noalign{\medskip}&1&&1&
\\\noalign{\medskip}1&&0&&0\end {smallmatrix} \right] , \left[ 
\begin {smallmatrix} &0&&1&\\\noalign{\medskip}0&&1&&1
\\\noalign{\medskip}&1&&0&\\\noalign{\medskip}0&&0&&-1
\end {smallmatrix} \right] , \left[ \begin {smallmatrix} &0&&1&
\\\noalign{\medskip}1&&1&&1\\\noalign{\medskip}&1&&0&
\\\noalign{\medskip}0&&0&&0\end {smallmatrix} \right] , \left[ 
\begin {smallmatrix} &1&&1&\\\noalign{\medskip}1&&2&&0
\\\noalign{\medskip}&1&&1&\\\noalign{\medskip}0&&0&&1
\end {smallmatrix} \right] , \left[ \begin {smallmatrix} &1&&0&
\\\noalign{\medskip}1&&1&&0\\\noalign{\medskip}&0&&1&
\\\noalign{\medskip}-1&&0&&0\end {smallmatrix} \right] , \left[ 
\begin {smallmatrix} &1&&0&\\\noalign{\medskip}1&&1&&1
\\\noalign{\medskip}&0&&1&\\\noalign{\medskip}0&&0&&0
\end {smallmatrix} \right] ]
\] \[
[[0,5], \left[ \begin {smallmatrix} &1&&1&
\\\noalign{\medskip}0&&2&&1\\\noalign{\medskip}&1&&1&
\\\noalign{\medskip}0&&0&&0\end {smallmatrix} \right] , \left[ 
\begin {smallmatrix} &0&&1&\\\noalign{\medskip}1&&1&&1
\\\noalign{\medskip}&1&&0&\\\noalign{\medskip}0&&0&&-1
\end {smallmatrix} \right] , \left[ \begin {smallmatrix} &1&&1&
\\\noalign{\medskip}1&&2&&1\\\noalign{\medskip}&1&&1&
\\\noalign{\medskip}0&&0&&1\end {smallmatrix} \right] , \left[ 
\begin {smallmatrix} &1&&1&\\\noalign{\medskip}1&&2&&0
\\\noalign{\medskip}&1&&1&\\\noalign{\medskip}0&&0&&0
\end {smallmatrix} \right] , \left[ \begin {smallmatrix} &1&&0&
\\\noalign{\medskip}1&&1&&1\\\noalign{\medskip}&0&&1&
\\\noalign{\medskip}-1&&0&&0\end {smallmatrix} \right] , \left[      
\begin {smallmatrix} &1&&1&\\\noalign{\medskip}1&&2&&1
\\\noalign{\medskip}&1&&1&\\\noalign{\medskip}1&&0&&0
\end {smallmatrix} \right] ]
\]
\end{small}}

%%%%%%%%%%%%%%%%%%%%%%%%%%%%%%%%%%%%%%%%%%%%%%%%%%%%%%%%%

\end{document}